\documentclass{amsart}

\usepackage[english]{babel}

\usepackage[letterpaper,top=2cm,bottom=2cm,left=3cm,right=3cm,marginparwidth=1.75cm]{geometry}

\usepackage{amsmath}
\usepackage{graphicx}
\usepackage[colorlinks=true, allcolors=blue]{hyperref}
\usepackage{amsthm}
\usepackage{amsmath}
\usepackage{amsfonts}
\usepackage{amssymb}
\usepackage{xcolor}
\usepackage{cleveref}
\usepackage{mathrsfs}
\usepackage{enumitem}
\setlist[description]{
    style=multiline,
    leftmargin=2.5cm,
    labelwidth=2.2cm,
    labelsep=0.3cm
}

\newcommand{\ip}[2]{\langle {#1}, {#2} \rangle}

\newtheorem{maintheorem}{Theorem}

\newtheorem{maincorollary}[maintheorem]{Corollary}

\newtheorem{theorem}{Theorem}[section]

\newtheorem{definition}[theorem]{Definition}
\newtheorem{prop}[theorem]{Proposition}
\newtheorem{corollary}[theorem]{Corollary}
\newtheorem{lemma}[theorem]{Lemma}

\theoremstyle{definition}
\newtheorem{remark}[theorem]{Remark}
\newtheorem*{notation}{Notation}

\theoremstyle{remark}

\title{A splitting theorem for Einstein manifolds with symmetries}
\author{Marcos Martinez Wagner}

\begin{document}
\begin{abstract}
We prove that every negative-Einstein manifold with cocompact symmetries splits isometrically as a Riemannian product $H \times K$ with $H$ homogeneous and $K$ compact.
This confirms a conjecture due to Böhm and Lafuente.
As a consequence, we obtain a new proof of the 1975 Alekseevskii conjecture.
\end{abstract}
\maketitle

\tableofcontents

\section{Introduction}\label{sec:introduction}

A Riemannian manifold $M=(M,\mathrm{g})$ is called \emph{Einstein} if its Ricci tensor satisfies $\mathrm{ric}_{\mathrm{g}}=\lambda\,\mathrm{g}$ for some $\lambda\in\mathbb{R}$.
When $\lambda<0$, we say that $M$ is \emph{negative-Einstein}.
In this article we study negative-Einstein manifolds with \emph{cocompact symmetries}, meaning that a connected Lie group $\mathsf{G}$ acts properly and isometrically on $M$ with compact orbit space $M/\mathsf{G}$. Our main result confirms a conjecture of B\"ohm and Lafuente \cite[Conjecture~B]{bohm2025non}.

\begin{maintheorem}\label{thm:main-splitting}
Any connected negative-Einstein manifold with cocompact symmetries is isometric to a Riemannian product $H\times K$, with $H$ a homogeneous manifold and $K$ compact.
\end{maintheorem}

The factor $H$ is a $\mathsf{G}$-homogeneous manifold and $M$ is equivariantly isometric to $H\times K$.
In particular, all $\mathsf{G}$-orbits are equivariantly isometric to $H$ and the orbit space $M/\mathsf{G}$ is a smooth manifold isometric to $K$.
Each factor is negative-Einstein with the same Einstein constant as $M$.

B\"ohm and Lafuente established this splitting for semisimple symmetry groups in \cite{bohm2025non}; in this case, the homogeneous factor is a symmetric space of non-compact type.
More generally, they proved that a closed, connected, unimodular subgroup of $\mathsf{Isom}(M)$ acting cocompactly on a negative-Einstein manifold must be semisimple.
Theorem~\ref{thm:main-splitting} establishes the product splitting without imposing unimodularity on the symmetry group.
Our proof builds partially on their earlier structural study of negative-Einstein manifolds with cocompact symmetries \cite{BL23}.

\subsubsection*{The role of the hypotheses:}

The connectedness of $\mathsf{G}$ is essential.
Let $S_2$ be a closed hyperbolic surface of genus two and let $M\to S_2$ be the infinite cyclic covering associated with a surjection $\pi_1(S_2)\to\mathbb{Z}$; see, for example, \cite[Section~3]{agol2016pseudo}.
With the lifted metric, $M$ is not homogeneous, since it is not simply connected \cite{BL23}.
Its deck group $\mathbb{Z}$ acts isometrically and cocompactly.
Since $M$ is noncompact and nonhomogeneous, an isometric product splitting as in Theorem~\ref{thm:main-splitting} would have two positive-dimensional factors.
Both would therefore be one-dimensional, making their product flat, a contradiction.

Compactness of the orbit space is also essential.
The complete hyperbolic cylinder
$$
M=\mathbb{R}\times\mathbb{S}^1,
\qquad
\mathrm{g}=dt^2+\cosh^2(t)\,d\theta^2
$$
is $\mathbb{S}^1$-invariant and noncompact.
Moreover, every isometry preserves the closed geodesic $\{0\}\times\mathbb{S}^1$, so $M$ is not homogeneous.
The same dimension argument excludes a product splitting as in Theorem \ref{thm:main-splitting}.

The negativity of the Einstein constant is essential for the equivariant refinement, as shown by the Hopf circle action on the round $\mathbb{S}^3$.

\subsubsection*{The homogeneous case:}

Theorem~\ref{thm:main-splitting} yields a new proof of the Alekseevskii conjecture, including its solvmanifold formulation.

\begin{maincorollary}[Alekseevskii Conjecture {\cite{BL23}}]\label{cor:alekseevskii}
Every connected homogeneous negative-Einstein manifold is isometric to a simply connected solvable Lie group endowed with a left-invariant Riemannian metric.
In particular, it is diffeomorphic to $\mathbb{R}^n$.
\end{maincorollary}

To prove this, choose a closed connected subgroup $\mathsf{G}\subseteq\mathsf{Isom}(M)$ acting transitively on $M$.
One can show that $\mathsf{G}$ is centerless (see Remark~\ref{rem:trivial-center}).
It therefore admits a \emph{generalized Iwasawa decomposition} $\mathsf{G}=\mathsf{K}\mathsf{S}$ (see Theorem~\ref{thm:generalized-iwasawa}), where $\mathsf{K}$ is maximal compact and $\mathsf{S}$ is closed, connected, simply connected and solvable.
The quotient $\mathsf{G}/\mathsf{S}$ is compact, hence $M/\mathsf{S}$ is compact by Lemma~\ref{lem:cocompact-subgroup}.
The equivariant refinement of Theorem~\ref{thm:main-splitting} gives an isometric splitting into negative-Einstein factors,
$$
M\equiv(\mathsf{S}\cdot x_0)\times M/\mathsf{S}.
$$
The factors of a homogeneous Riemannian product are homogeneous; see \cite[Chapter~VI, Theorem~3.5]{KobayashiNomizu63}.
Thus, $M/\mathsf{S}$ is compact, homogeneous and negative-Einstein, so Bochner's theorem \cite{bochner1946vector} forces $M/\mathsf{S}=\{*\}$.
Consequently, $\mathsf{S}$ acts transitively on $M$.
Its action is proper because $\mathsf{S}$ is closed, and hence free since solvable simply connected groups have no non-trivial compact subgroups (see Lemma~\ref{lem:simply-connected-solvable-free-action}).
The orbit map therefore identifies $M$ isometrically with $\mathsf{S}$ endowed with a left-invariant Riemannian metric. The same argument shows that a closed connected subgroup $\mathsf{H}\subseteq \mathsf{G}$ acts transitively on $M$ if and only if $\mathsf{G}/\mathsf{H}$ is compact.

The structure theory for Einstein solvmanifolds was developed by Heber \cite{Heber98} under the standardness assumption, which Lauret subsequently proved for all Einstein solvmanifolds \cite{standard}.
Lauret's work on nilsolitons \cite{Lauret01,lauret2002finding} relates this theory to Ricci solitons on nilpotent Lie groups and geometric invariant theory.
The algebraic structure and symmetries of Einstein solvmanifolds have been studied in \cite{MR936815,Nik11,Jbl2015,Jab15,GJ19}.
Beyond the solvable setting, structural results of Lafuente--Lauret \cite{lafuente2014structure} and Jablonski--Petersen \cite{jablonski2017step} preceded the proof of the Alekseevskii conjecture by B\"ohm and Lafuente \cite{BL23}.
See also the surveys \cite{Lauret09survey,Jablonski23survey} and the references therein.

In dimensions at most eight, homogeneous negative-Einstein manifolds are classified up to isometry and scaling by this theory.
Heber's structure theory and Lauret's standardness theorem describe the metrics by nilsolitons and their solvable extensions.
The nilradical has dimension at most seven in this range, and the required nilsolitons are classified through dimension six in \cite{lauret2002finding,Will03} and in dimension seven in \cite{FernandezCulma12,FernandezCulma14}.

\subsubsection*{Low cohomogeneity:}

Any two closed connected subgroups of $\mathsf{Isom}(M)$ acting cocompactly on $M$ have the same orbits and hence the same cohomogeneity; see Corollary~\ref{cor:canonical-orbit-decomposition}.

For any such subgroup $\mathsf{G}$, the equivariant splitting and Corollary~\ref{cor:alekseevskii} give a diffeomorphism
$$
M\cong\mathbb{R}^{d}\times M/\mathsf{G},
$$
where $d$ is the dimension of the $\mathsf{G}$-orbits.
The projection $M\rightarrow M/\mathsf{G}$ is therefore a homotopy equivalence, and the quotient carries a negative-Einstein metric with the same Einstein constant as $M$.
This has particularly strong consequences in low cohomogeneity.
A one-dimensional quotient would be a flat circle, which is incompatible with its negative Einstein constant.

\begin{maincorollary}\label{cor:no-cocompact-cohomogeneity-one}
Let $M$ be a connected negative-Einstein manifold.
Then there is no closed connected subgroup $\mathsf{G}\subseteq\mathsf{Isom}(M)$ acting cocompactly and with cohomogeneity one.
\end{maincorollary}
Compare with \cite[Theorem~D]{bohm2022euclidean}, which establishes a rigidity result for cohomogeneity-one actions with orbit space $\mathbb{S}^1$, imposing the negative-Einstein equation only in vertical directions.
Such actions do exist: for example, any Lie group $G$ acting on the first factor of $H \times \mathbb{S}^1$ for $H$ a $G$-homogeneous negative-Einstein manifold.

If the cohomogeneity is two or three, then the quotient
$M/\mathsf{G}$ has constant negative sectional curvature
\cite{Besse87}.
Together with the splitting, this implies that the universal cover
of $M$ is a product of homogeneous spaces.

\begin{maincorollary}\label{cor:low-cohomogeneity}
Let $M$ be a connected negative-Einstein manifold, and suppose that a closed connected subgroup $\mathsf{G}\subseteq\mathsf{Isom}(M)$ acts cocompactly with cohomogeneity two or three.
Then $M$ is locally homogeneous, and its universal cover is diffeomorphic to $\mathbb{R}^n$.
In particular, $M$ is aspherical: $\pi_i(M)=0$ for every $i\geq2$.
\end{maincorollary}

\subsubsection*{Obstructions to invariant negative-Einstein metrics:}

Suppose that a connected Lie group $\mathsf{G}$ acts smoothly, effectively, properly and cocompactly on a connected smooth manifold $M$, with cohomogeneity $k$.
By the equivariant refinement of Theorem~\ref{thm:main-splitting}, if this action has a nonprincipal orbit, then $M$ carries no $\mathsf{G}$-invariant negative-Einstein metric.
There are topological obstructions as well: the existence of such a metric would force
$$
H_i(M;\mathbb{Z})=0\quad(i>k),
\qquad
H_k(M;\mathbb{Z}_2)\cong\mathbb{Z}_2.
$$
If $M$ is orientable, then also $H_k(M;\mathbb{Z})\cong\mathbb{Z}$; if $k$ is odd, then $\chi(M)=0$.
These restrictions follow from the homotopy equivalence $M\simeq M/\mathsf{G}$ and Poincar\'e duality on the compact quotient \cite{Hatcher02}.

In cohomogeneity two, an orientable manifold with odd first Betti number cannot admit a $\mathsf{G}$-invariant negative-Einstein metric.
Indeed, such a metric would make $M/\mathsf{G}$ a closed orientable surface of genus $g$, and the homotopy equivalence with $M$ would give $b_1(M)=2g$.
Moreover, the quotient metric has constant negative sectional curvature, so Gauss--Bonnet forces $g\geq2$, and hence $b_1(M)\geq4$.

\subsubsection*{Manifolds with parallel Ricci tensor:}
The product conclusion also extends beyond the Einstein condition.
A complete connected simply connected manifold with parallel Ricci tensor splits into Einstein factors by de Rham's theorem \cite{derham1952reductibilite}.
Applying Theorem~\ref{thm:main-splitting} to the negative-Einstein factors, Myers's theorem \cite{myers1941riemannian} to the positive-Einstein factors, and the Ricci-flat splitting result of B\"ohm and Lafuente \cite[Theorem~9.2]{bohm2025non} to the remaining factors gives the following extension.

\begin{maincorollary}\label{cor:parallel-ricci-splitting}
Let $M$ be a connected simply connected Riemannian manifold with parallel Ricci tensor and cocompact symmetries.
Then $M$ splits isometrically as $H\times K$, where $H$ is a homogeneous manifold with non-positive Ricci curvature and $K$ is compact.
Furthermore,
$$
H=\mathbb{R}^r\times H_1\times\cdots\times H_m,
\qquad
K=K_1\times\cdots\times K_\ell,
$$
where each $H_i$ is a homogeneous negative-Einstein manifold and each $K_i$ is a compact Einstein manifold.
\end{maincorollary}

\subsection{Outline of the proof of Theorem A}\label{subsec:proof-outline}

\textcolor{white}{aaa}

Rescaling, we may assume that the Einstein constant is $-1$.
We take $\mathsf{G}$ to be a closed, connected subgroup of $\mathsf{Isom}(M)$ acting cocompactly; see Subsection~\ref{subsec:properties-isometric-actions}.
The proof has two parts:
\begin{enumerate}
    \item Prove the theorem when $\mathsf{G}$ is solvable and acts freely on $M$ (Theorem~\ref{thm:solvable-splitting}).
    \item Reduce to this case using a generalized Iwasawa decomposition $\mathsf{G}=\mathsf{K}\mathsf{S}$, where $\mathsf{K}$ is maximal compact and $\mathsf{S}$ is a closed, connected, simply connected solvable subgroup (Appendix~\ref{app:generalized-iwasawa}).
\end{enumerate}
For (2), the $\mathsf{S}$-action on $M$ is free, proper and cocompact, so part (1) applies.
To meet the finite-center hypothesis needed for the group decomposition, we first use a Bochner argument to show that the center of $\mathsf{G}$ is discrete, and then quotient out a torsion-free subgroup of finite index in said center.
This reduction is carried out in Section~\ref{sec:general-splitting}.

The solvable case is treated in Sections~\ref{sec:shape-commutativity}--\ref{sec:solvable-splitting} where we work under the following assumption.

\begin{trivlist}
\item[\hskip\labelsep\bfseries Assumption F.]
\phantomsection\makeatletter\def\@currentlabel{F}\makeatother\label{ass:main}
\itshape
$(M,\mathrm{g})$ is a connected, complete Riemannian manifold with $\mathrm{ric}_{\mathrm{g}}=-\mathrm{g}$, and $\mathsf{G}$ is a connected solvable Lie group acting freely, properly, isometrically and cocompactly on $M$.
\end{trivlist}

The projection $M\to M/\mathsf{G}$ is a principal bundle with compact base.
We equip this bundle with the principal connection $\mathcal{H}:=\mathcal{V}^{\perp}$ where $\mathcal{V}_x:=T_x(G\cdot x)$ is the vertical distribution.
Our strategy is based on the following criterion; see Subsection~\ref{subsec:holonomy} for the proof.

\begin{prop}\label{prop:splitting-criterion-outline}
Let $M$ be a complete, connected Riemannian manifold, and let $\mathsf{G}$ be a connected Lie group acting freely, properly and isometrically on $M$.
Then the following conditions are equivalent:
\begin{enumerate}
\renewcommand{\labelenumi}{(\roman{enumi})}
    \item The principal bundle $M\to M/\mathsf{G}$ admits a global $\mathsf{G}$-equivariant isometric trivialization.
    \item The $\mathsf{G}$-orbits are totally geodesic, and the holonomy of $M\to M/\mathsf{G}$ is trivial.
\end{enumerate}
\end{prop}

We briefly explain how the two conditions in~(ii) give a global isometric trivialization.
For $x\in M$, the holonomy group $\mathsf{Hol}_x\subseteq\mathsf{G}$ consists of those $g\in\mathsf{G}$ for which $x$ and $g\cdot x$ can be joined by a horizontal curve.
Trivial holonomy gives a horizontal section $\Sigma$ meeting every orbit exactly once.
Fixing $x_0\in\Sigma$, the map $(g\cdot x_0,y)\mapsto g\cdot y$ identifies $(\mathsf{G}\cdot x_0)\times\Sigma$ with $M$.
Total geodesicity of the orbits makes this identification isometric.
Our task is therefore to prove $L=0$ and $\mathsf{Hol}_x=\{e\}$, where $L$ denotes the shape operator of the $\mathsf{G}$-orbits.

Fundamental vector fields give a global trivialization of the vertical bundle.
In this trivialization, a vertical vector field $U$ is represented by a $\mathfrak{g}$-valued function ${}^{\mathfrak{g}}U$, and a vertical endomorphism $E$ by an $\mathrm{End}(\mathfrak{g})$-valued function ${}^{\mathfrak{g}}E$.
We call these representations their \emph{Killing coordinates}.
The orbit metrics are likewise represented by a map $h:M\to\mathrm{Sym}^2_+(\mathfrak{g}^*)$.
For a horizontal vector $X$, the well known variation formula (see e.g.\cite[Lemma~5.4]{BL23}) reads, in Killing coordinates,
$$dh(X)=2h\bigl({}^{\mathfrak{g}}L_X\,\cdot,\cdot\bigr).$$
In particular, $L=0$ is equivalent to $h$ being constant along horizontal curves.

B\"ohm and Lafuente \cite{BL23} study the action of the connected nilradical $\mathsf{N}\subseteq\mathsf{G}$ under hypotheses more general than Assumption~\ref{ass:main}.
They prove that the $\mathsf{N}$-horizontal distribution is integrable \cite[Theorem~8.1]{BL23}.
In our setting, \cite[Corollary~8.7]{BL23} gives an $h_x$-orthogonal decomposition $\mathfrak{g}=\mathfrak{a}_x\oplus\mathfrak{n}$ with $\mathfrak{a}_x$ an abelian subalgebra and $x\mapsto \mathfrak{a}_x$ constant along $\mathsf{G}$-horizontal curves.
This extends Lauret's standardness theorem for Einstein solvmanifolds \cite{standard} to the cocompact setting.

Let $\widehat{\mathcal{V}}$ and $\widehat{\mathcal{H}}$ denote the distributions associated to the $\mathsf{N}$-action.
Set $\mathcal{A}:=\mathcal{V}\cap\widehat{\mathcal{H}}$, so that $\mathcal{A}_x$ corresponds to $\mathfrak{a}_x$ in Killing coordinates.
Then, we have the orthogonal decomposition
$$
TM=\mathcal H\oplus
\underbrace{\mathcal A\oplus\widehat{\mathcal V}}_{\mathcal V}.
$$
For each $X\in \mathcal{H}$, the operator $L_X$ preserves each of the distributions $\mathcal{A}$ and $\widehat{\mathcal{V}}$.
Moreover, $L_X|_{\widehat{\mathcal{V}}}=\hat{L}_X$, where $\hat{L}$ is the shape operator of the $\mathsf{N}$-orbits (Lemma~\ref{lem:vertical-tensor-blocks}).
Thus, proving $L=0$ amounts to showing that $L_X|_{\mathcal{A}}=0$ and $\hat{L}_X=0$ for every such $X$.

We prove the vanishing of these two restrictions in Sections~\ref{sec:affine-conjugation-orbits}--\ref{sec:solvable-splitting}, building on the results of Sections~\ref{sec:shape-commutativity}--\ref{sec:trivial-n-holonomy}. 
We first introduce some preliminary tools that will be used extensively throughout the paper.
\begin{prop}[see Proposition~\ref{prop:cocompact-maximum-principle}]\label{prop:compact-divergence-maximum}
Let $B$ be a compact, connected Riemannian manifold without boundary, let $X$ be a smooth vector field, and let $u$ be a smooth function on $B$.
\begin{enumerate}
\renewcommand{\labelenumi}{(\roman{enumi})}
    \item If $\mathrm{div}\,X-\langle X,\nabla u\rangle\geq0$, then equality holds everywhere.
    \item If $\Delta u-\langle X,\nabla u\rangle\geq0$, then $u$ is constant.
\end{enumerate}
\end{prop}

A vertical field is \emph{horizontally Killing} if its Killing coordinates are constant along horizontal curves.
Every vector in $\mathcal{A}_x$ extends uniquely to an invariant, horizontally Killing, $\mathcal{A}$-valued field (Proposition~\ref{prop:invariant-a-extension}), by the equivariance and horizontal constancy of $x\mapsto\mathfrak{a}_x$ and the abelianness of $\mathfrak{a}_x$.

In Section~\ref{sec:shape-commutativity}, we prove $\mathrm{ad}_A$ is $h_x$-normal for each $A\in\mathfrak{a}_x$, generalizing \cite[Theorem~4.10(1)]{Heber98} to the cocompact setting. Furthermore, for every $x\in M$, every $A,\tilde{A}\in \mathfrak{a}_x$ and every $X_x\in \mathcal{H}_x$ we prove
\begin{equation}\label{cososo}
     [\mathrm{ad}_{\tilde{A}},\mathrm{ad}_A^{t_x}]=0,
    \qquad
    [{}^{\mathfrak{g}}L_{X_x},\mathrm{ad}_A]=0.
\end{equation}
To prove this, we associate to each $\mathcal{A}$-valued field $A$ a horizontal field $\mathrm{W}_A$ defined by
$$\ip{\mathrm{W}_A}{X}=\mathrm{Tr}([\mathrm{ad}_{{}^{\raisebox{0.15em}{$\scriptstyle\mathfrak{g}$}}\! A},\mathrm{ad}_{{}^{\raisebox{0.15em}{$\scriptstyle\mathfrak{g}$}}\! A}^t]\cdot{^\mathfrak{g}L_X}) \qquad \text{for all} \ X\in \mathcal{H}$$
If $A$ is invariant, so is $\mathrm{W}_A$, which descends to $M/\mathsf{G}$.
If $A$ is also horizontally Killing, the vertical O'Neill equations and Heber's Ricci estimate \cite[proof preceding Lemma~4.1]{Heber98} yield an inequality as in Proposition~\ref{prop:compact-divergence-maximum}(i) with $X=\mathrm{W}_A$. Equality then holds, and the extension result yields \eqref{cososo}.

Section~\ref{sec:algebraic-rigidity} uses the commutativity relations established in Section~\ref{sec:shape-commutativity} to extend further structural properties of Einstein solvmanifolds to our setting.
In particular, $\mathsf{G}$ is simply connected and centerless.

Denote $\widehat{\mathsf{Hol}}_x$ the holonomy of $M\to M/\mathsf{N}$. In Section~\ref{sec:trivial-n-holonomy}, we prove
\[
    \widehat{\mathsf{Hol}}_x=\{e\},
    \qquad
    \mathsf{Hol}_x\subseteq\exp(\mathfrak{a}_x).
\]
We use the centerlessness of $\mathsf{G}$ and the fact that the GIT map ${}^{\mathfrak{n}}\beta^+:M\to\mathrm{End}(\mathfrak{n})$ is $\mathsf{G}$-equivariant, constant along $\mathsf{N}$-horizontal curves and pointwise a derivation \cite[Theorem~8.1 and Corollary~8.4]{BL23}.

We now return to proving $\hat{L}_X=0$, or equivalently, that $\hat{h}:=h|_{\mathfrak{n}}$ is constant along $\mathsf{G}$-horizontal curves.
Fix an $\mathsf{N}$-horizontal leaf $\widehat{\Sigma}$ and let $\mathfrak{a}_0$ be the common value of $\mathfrak{a}_x$ along it.
Since the $\mathsf{N}$-holonomy is trivial, $\mathsf{A}_0:=\exp(\mathfrak{a}_0)$ acts on $\widehat{\Sigma}$ with compact quotient $\widehat{\Sigma}/\mathsf{A}_0\equiv M/\mathsf{G}$.
The map
\[
    \hat{h}|_{\widehat{\Sigma}}:
    \widehat{\Sigma}\longrightarrow\mathrm{Sym}^2_+(\mathfrak{n}^*)
\]
is $\mathsf{A}_0$-equivariant, where the action on the target is induced by $\mathrm{Ad}(\mathsf{A}_0)|_{\mathfrak{n}}\subseteq\mathrm{GL}(\mathfrak{n})$.
We will show that, for every $x\in\widehat{\Sigma}$ and $X_x\in\mathcal{H}_x$, the vector $d\hat{h}_x(X_x)$ is both tangent and orthogonal to $\mathsf{A}_0\cdot\hat{h}_x$, and hence vanishes.

The tangency statement is proved in Section~\ref{sec:affine-conjugation-orbits}.
Fixing $x_0\in\widehat{\Sigma}$, we show that the orbit $O:=\mathsf{A}_0\cdot\hat{h}_{x_0}$ is a closed, totally geodesic submanifold of $\mathrm{Sym}^2_+(\mathfrak{n}^*)$.
Since this symmetric space is Hadamard, the squared distance $\varphi(h):=d(h,O)^2$ is smooth and convex.
By the harmonicity and equivariance of $\hat{h}|_{\widehat{\Sigma}}$ \cite[Theorem~2.2(iii)]{bohm2025non}, the composition $\varphi\circ\hat{h}$ is subharmonic and $\mathsf{A}_0$-invariant.
Compactness of $\widehat{\Sigma}/\mathsf{A}_0$ ensures it attains a maximum, so the strong maximum principle makes it constant.
It vanishes at $x_0$, hence everywhere.
Thus the whole image of $\hat{h}|_{\widehat{\Sigma}}$ lies in $O$, proving the required tangency (Theorem~\ref{thm:a0-orbit-image}).

Section~\ref{sec:solvable-splitting} proves that $M/\mathsf{N}$ has non-positive Ricci curvature (Theorem~\ref{thm:orbit-ricci-lower-bound}).
By the identity
\[
    \mathrm{ric}^{M/\mathsf{N}}(\widehat{X},\widehat{X})
    =-\|\widehat{X}\|^2+\|\hat{L}_{\widehat{X}}\|^2
\]
from \cite[Theorem~2.2(iv)]{bohm2025non}, this is equivalent to $\|\hat{L}_{\widehat{X}}\|\leq\|\widehat{X}\|$ for every $\mathsf{N}$-horizontal $\widehat{X}$.
The operator norm of the pointwise map $\widehat{X}\mapsto\hat{L}_{\widehat{X}}$ is $\mathsf{G}$-invariant and, since $M/\mathsf{G}$ is compact, attains a maximum at some $x_0\in M$.
Our goal is to show that this maximum is at most $1$. 
To do so, take the $\mathsf{N}$-horizontal section $\widehat{\Sigma}$ through $x_0$.
Each invariant, horizontally Killing, $\mathcal{A}$-valued field $A$ defines a one-form $\omega_A$ on $\widehat{\Sigma}$ by $\omega_A(\widehat{X}):=\mathrm{Tr}(\hat{L}_A\hat{L}_{\widehat{X}})$.
The commutativity relations from Section~\ref{sec:shape-commutativity} make these forms harmonic.
The tangency statement proved in Section~\ref{sec:affine-conjugation-orbits} allows us to choose $A$ so that $\|\omega_A\|$ attains a maximum at the same point $x_0$.
The Bochner formula at $x_0$ then gives the desired bound.

We now use this estimate to prove $L=0$.
By Lauret's Ricci formula \cite[equation~(25)]{Lauret11} and the preceding bound, we have $\mathrm{ric}^{\mathsf{G}}(A,A)=-\|\hat{L}_A\|^2\geq-\|A\|^2$ for every $\mathcal{A}$-valued field $A$.
Choose $A$ to be invariant and horizontally Killing, as permitted by the extension result, and set $f:=\|A\|^2/2$.
Invariance makes $f$ a function on the compact quotient $M/\mathsf{G}$.
The following identity is an extension of \cite[Lemma~8.1]{bohm2025non} to invariant horizontally Killing fields (see Lemma~\ref{lem:horizontally-killing-norm}):
\[
    \Delta_Bf-\langle\nabla f,H\rangle
    =\mathrm{ric}^{\mathsf{G}}(A,A)+\|A\|^2
    +\sum_i\bigl(2\|L_{X_i}A\|^2+\|\mathrm{A}_{X_i}A\|^2\bigr).
\]
Here $\Delta_B$ denotes the Laplace--Beltrami operator of $M/\mathsf{G}$, $H$ is the mean curvature of the $\mathsf{G}$-orbits, $\mathrm{A}$ is O'Neill's integrability tensor, and $\{X_i\}$ is an orthonormal horizontal frame.
The Ricci estimate makes the right-hand side a sum of nonnegative quantities, so the maximum principle makes $f$ constant and forces equality.
Thus, we obtain $L_X(A)=0$ and $\|\hat{L}_A\|=\|A\|$.
Since the value of $A$ can be prescribed arbitrarily in $\mathcal{A}_x$, this gives $L_X|_{\mathcal{A}}=0$.
Combining the shape-operator bound with $\|\hat{L}_A\|=\|A\|$, we show that $\omega_A(\widehat{X})=\langle A,\widehat{X}\rangle$ for every $\widehat{X}\in \widehat{\mathcal{H}}$.
In particular, for every $\mathsf{G}$-horizontal vector $X$,
\[
    \mathrm{Tr}(\hat{L}_A\hat{L}_X)
    =\omega_A(X)
    =\langle A,X\rangle
    =0.
\]
Equivalently, $d\hat{h}_x(X_x)$ is perpendicular to $\mathsf{A}_0\cdot \hat{h}_x$.
This proves $\hat{L}_X=0$, and hence $L=0$.

Finally, let $a\in\mathsf{Hol}_x\subseteq\exp(\mathfrak{a}_x)$.
$\mathsf{G}$-horizontal constancy and $\mathsf{G}$-equivariance of $\hat{h}$ give
\[
    \mathrm{Ad}(a)|_{\mathfrak{n}}\cdot\hat{h}_x
    =\hat{h}_{a\cdot x}
    =\hat{h}_x.
\]
Since $\mathsf{A}_0:=\exp(\mathfrak{a}_x)$ acts freely on the space of inner products on $\mathfrak{n}$, we get $a=e$, proving $\mathsf{Hol}_x=\{e\}$.

\subsection{Organization of the article}\label{subsec:organization}

Sections~\ref{sec:isometric-group-actions}--\ref{sec:adjoint-connection} contain preliminary material.
In Section~\ref{sec:isometric-group-actions}, we review free proper isometric actions, O'Neill tensors, holonomy and Killing coordinates.
Section~\ref{sec:einstein-cocompact-symmetries} recalls results from \cite{BL23} on Einstein manifolds with cocompact symmetries, focusing on the nilradical action.
In Section~\ref{sec:adjoint-connection}, we express the adjoint bundle connection in terms of ordinary derivatives of Killing coordinates and introduce horizontally Killing vector fields.

Sections~\ref{sec:shape-commutativity}--\ref{sec:solvable-splitting} prove Theorem~\ref{thm:main-splitting} for free actions of solvable groups, following the strategy outlined in Subsection~1.1.
Section~\ref{sec:general-splitting} establishes its general equivariant version.

Appendix~\ref{app:inner-products} describes the geometry of the symmetric space of inner products.
Appendix~\ref{app:generalized-iwasawa} proves the generalized Iwasawa decomposition used in Section~\ref{sec:general-splitting}.
A notation index is included at the end.

\subsubsection*{Acknowledgements}\phantomsection\makeatletter\def\@currentlabel{Acknowledgements}\makeatother\label{sssec:acknowledgements} 
I am deeply grateful to my supervisor, Ramiro Lafuente, for suggesting this problem to me, for his guidance, and for his feedback.
This project has been supported by an Australian Government Research Training Program (RTP)
scholarship

\subsubsection*{AI disclosure}\phantomsection\makeatletter\def\@currentlabel{AI disclosure}\makeatother\label{sssec:ai-disclosure}
During the preparation of most of this paper, AI software such as ChatGPT, Claude, and Gemini served as pseudo-reliable sources of known facts.
With the recent advances in AI capabilities in mathematics, ChatGPT 5.6 and ChatGPT 6.1 have made the following two main contributions to this paper:

\begin{itemize}
    \item Theorem~\ref{thm:orbit-ricci-lower-bound}: the orbit space $M/\mathsf{N}$ has non-positive Ricci curvature.
    \item Theorem~\ref{thm:generalized-iwasawa}: a generalized Iwasawa decomposition for connected Lie groups with finite center.
\end{itemize}

For the first result, I supplied ChatGPT 5.6 (Sol) with the results I had obtained up to that point, including some not included in this paper.
It combined the harmonicity of $\omega_A$ from Proposition~\ref{prop:harmonic-shape-one-form} and the identity
$\hat{L}_{\widehat{X}}=\hat{L}_{\alpha(\widehat{X})}$
from Proposition~\ref{prop:shape-alpha} with a maximum principle strategy suggested by Ramiro Lafuente for bounding the eigenvalues of the shape operator.

For the second result, while working on an unrelated task, ChatGPT 6.1 (Astra) produced the solvable \emph{Iwasawa-like} decomposition described in Proposition~\ref{prop:solvable-radical-decomposition}.
I then asked it to combine this with a Levi decomposition to obtain the general statement.
The resulting statement and proof was developed through several iterations.
The theorem can also be deduced from the proof of \cite[Lemma~6.7]{Cornulier08}, as explained in Remark~\ref{rem:cornulier-iwasawa}.

ChatGPT 6.1 (Sol) was also used to assist with LaTeX formatting and to improve the clarity and English of the exposition.

\section{Isometric group actions}\label{sec:isometric-group-actions}
In this section, we establish the framework and notation used throughout the paper, including O'Neill tensors, curvature formulas, Killing coordinates and holonomy.

Throughout this section, $M$ is a complete connected Riemannian manifold and $\mathsf{G}$ is a connected Lie group acting on $M$ freely, properly and isometrically.

\subsection{Isometric actions and Riemannian submersions}\label{subsec:actions-submersions}
\textcolor{white}{aaa}

Let $m_g:M\rightarrow M$ denote the action of $g\in \mathsf{G}$.
A vector field $\xi$ on $M$ is called \emph{invariant} if $\xi_{g\cdot x}=d(m_g)_x(\xi_x)$ for every $x\in M$ and every $g\in \mathsf{G}$.
A section $E\in\Gamma(\mathrm{End}\,TM)$ is called \emph{equivariant} if $d(m_g)_x\circ E_x=E_{g\cdot x}\circ d(m_g)_x$ for every $x\in M$ and every $g\in \mathsf{G}$.
A function $f:M\rightarrow\mathbb{R}$ is called \emph{invariant} if $f(g\cdot x)=f(x)$ for every $g\in \mathsf{G}$ and every $x\in M$.

For every $x\in M$, we denote by $\mathcal{V}_x:=T_x(\mathsf{G}\cdot x)$ the vertical distribution and by $\mathcal{H}_x:=\mathcal{V}_x^\perp\subseteq T_xM$ the horizontal distribution.
For a vector field $\xi$, we denote its vertical and horizontal components by $\mathcal{V}(\xi)$ and $\mathcal{H}(\xi)$, respectively.
Vertical vector fields will be denoted by $U,V,W$, while horizontal vector fields will be denoted by $X,Y,Z$.
Invariant horizontal vector fields are called \emph{basic} vector fields.

Since the action is free, proper and isometric, there is a unique Riemannian metric on $B:=M/\mathsf{G}$ such that the projection $\pi:M\rightarrow B$ is a Riemannian submersion.
That is, $d\pi_x:\mathcal{H}_x\rightarrow T_{\pi(x)}B$ is a linear isometry for every $x\in M$; see \cite[Chapter~9]{Besse87}.
The completeness of $M$ allows us to apply the following proposition.

\begin{prop}[R.
Hermann, {\cite[Theorem~9.42]{Besse87}}]\label{prop:horizontal-lifting}
Let $\pi:M\rightarrow B$ be a Riemannian submersion.
If $M$ is complete, then $B$ is complete and $\pi:M\rightarrow B$ is a locally trivial fibration.
Moreover:
\begin{enumerate}
\renewcommand{\labelenumi}{(\roman{enumi})}
    \item For every piecewise smooth curve $\bar{\sigma}:[0,1]\rightarrow B$ starting at $b\in B$ and every $x\in\pi^{-1}(b)$, there exists a unique piecewise smooth horizontal curve $\sigma:[0,1]\rightarrow M$ starting at $x$ such that $\pi\circ\sigma=\bar{\sigma}$.
    \item For every vector field $\bar{X}$ on $B$, there exists a unique basic vector field $X$ on $M$ such that $d\pi_x(X_x)=\bar{X}_{\pi(x)}$ for every $x\in M$.
\end{enumerate}
\end{prop}

The following lemma will be useful, for example, when we introduce horizontally Killing vector fields in Subsection~\ref{subsec:horizontally-killing}.

\begin{lemma}[{\cite[Lemma~3.5]{BL23}}]\label{lem:basic-vertical-bracket}
Let $X$ be a basic vector field and $U$ a vertical vector field.
Then $[X,U]$ is vertical.
Furthermore, if $U$ is Killing, then $[X,U]=0$.
\end{lemma}

\subsection{O'Neill tensors}\label{subsec:oneill-tensors}
\textcolor{white}{aaa}

Proofs of the standard identities recalled in this subsection can be found in \cite[Chapter~9]{Besse87}.
For every Riemannian submersion, such as the projection $\pi:M\rightarrow M/\mathsf{G}$, we define the associated O'Neill $(1,2)$-tensors $\mathrm{A}$ and $L$ by
$$
    \mathrm{A}_\xi\tilde{\xi}:=\mathcal{H}(\nabla_{\mathcal{H}(\xi)}\mathcal{V}(\tilde{\xi}))
    +\mathcal{V}(\nabla_{\mathcal{H}(\xi)}\mathcal{H}(\tilde{\xi}))
    \quad ; \quad
    L_\xi\tilde{\xi}:=\mathcal{V}(\nabla_{\mathcal{V}(\tilde{\xi})}\mathcal{H}(\xi)),
$$
for vector fields $\xi,\tilde{\xi}$ on $M$.
Both expressions are tensorial in their two variables.

\begin{remark}\label{rem:oneill-tensor-conventions}
In the literature, it is more common to consider the O'Neill tensor $\mathrm{T}$ instead of $L$; see \cite[Chapter~9]{Besse87}.
They are related by $L_XU=\mathrm{T}_UX$, so $L=0$ if and only if $\mathrm{T}=0$.
We use $L$ to follow the notation of \cite{BL23,bohm2025non}.
\end{remark}

Regarding the O'Neill tensor $\mathrm{A}$, notice that $\mathrm{A}_U=0$, $\mathrm{A}_XU=\mathcal{H}(\nabla_XU)$ and $\mathrm{A}_XY=\mathcal{V}(\nabla_XY)$.
The last two expressions satisfy
$$
    \ip{\mathrm{A}_XY}{U}=-\ip{Y}{\mathrm{A}_XU},
$$
so $\mathrm{A}|_{\mathcal{H}\times\mathcal{H}}=0$ if and only if $\mathrm{A}|_{\mathcal{H}\times\mathcal{V}}=0$, and either condition is equivalent to $\mathrm{A}=0$.
We also have
$$
    \mathrm{A}_XY=\frac{1}{2}\mathcal{V}[X,Y],
$$
so $\mathrm{A}=0$ if and only if the horizontal distribution is integrable.
Thus, for free proper isometric actions, $\mathrm{A}=0$ if and only if the action is \emph{polar}.

Regarding the O'Neill tensor $L$, notice that $L_U=0$, $L_XY=0$ and $L_XU=\mathcal{V}(\nabla_UX)$.
We regard $L_X$ as a section of $\mathrm{End}\,\mathcal{V}$; it is the shape operator of the orbits in the direction $X$, with the convention used in \cite{BL23}.

The following lemma is a consequence of Lemma~\ref{lem:orbit-metric-differential} below; see also \cite[Theorem~9.56]{Besse87}.

\begin{lemma}\label{lem:equivariant-orbits-totally-geodesic}
If the orbits are totally geodesic submanifolds, that is, if $L=0$, then all orbits are pairwise $\mathsf{G}$-equivariantly isometric.
\end{lemma}

The importance of the O'Neill tensors lies in the fact that they control the local geometry of the Riemannian submersion $M\rightarrow M/\mathsf{G}$ in the following sense.

\begin{prop}[{\cite[9.26]{Besse87}}]\label{prop:local-isometric-trivialization}
Let $M$ be a complete connected Riemannian manifold and let $\mathsf{G}$ be a connected Lie group acting on $M$ freely, properly and isometrically.
Then the following conditions are equivalent:
\begin{enumerate}
\renewcommand{\labelenumi}{(\roman{enumi})}
    \item Both O'Neill tensors $\mathrm{A}$ and $L$ of the Riemannian submersion $\pi:M\rightarrow M/\mathsf{G}$ vanish.
    \item The fiber bundle $\pi:M\rightarrow M/\mathsf{G}$ has local isometric trivializations.
\end{enumerate}
\end{prop}

In Subsection~\ref{subsec:holonomy}, we will see how to obtain a global isometric trivialization.

\subsection{O'Neill curvature formulas}\label{subsec:oneill-curvature}
\textcolor{white}{aaa}

To introduce O'Neill's vertical curvature formula, relating the Ricci curvature of $M$ to the Ricci curvature of its orbits, we first define some vertical tensors associated with $\mathrm{A}$ and $L$.
The Ricci curvature of the orbits induces a vertical $(1,1)$-tensor $\mathrm{Ric}^{\mathsf{G}}\in\Gamma(\mathrm{End}\,\mathcal{V})$ on $M$.
If $\mathrm{Ric}$ denotes the Ricci endomorphism of $M$, we define $\mathrm{Ric}^{\mathcal{V}}\in\Gamma(\mathrm{End}\,\mathcal{V})$ by $\mathrm{Ric}^{\mathcal{V}}U:=\mathcal{V}(\mathrm{Ric}\,U)$.
We seek a formula relating $\mathrm{Ric}^{\mathsf{G}}$ to $\mathrm{Ric}^{\mathcal{V}}$.

Throughout this subsection, $\{X_1,\ldots,X_n\}$ is a local orthonormal horizontal frame and $\{U_1,\ldots,U_m\}$ is a local orthonormal vertical frame.
All the definitions below are independent of the choice of these frames.

First, we define $\mathrm{A}^* \mathrm{A}\in\Gamma(\mathrm{End}\,\mathcal{V})$ by
$$
    \ip{\mathrm{A}^* \mathrm{A}U}{V}:=\sum_i\ip{\mathrm{A}_{X_i}U}{\mathrm{A}_{X_i}V}.
$$

Secondly, for horizontal vector fields $X,Y$, the endomorphism $(\nabla_XL)_Y$ need not preserve the vertical distribution.
In fact, the horizontal component of $(\nabla_XL)_Y(U)$ is $\mathrm{A}_X(L_YU)$.
We therefore define $\tau\in\Gamma(\mathrm{End}\,\mathcal{V})$ by
$$
    \tau(U):=\mathcal{V}\left(\sum_i(\nabla_{X_i}L)_{X_i}(U)\right).
$$
Both $\mathrm{A}^* \mathrm{A}$ and $\tau$ are pointwise self-adjoint.

Lastly, we define the mean curvature vector field $H$ of the orbits by
$$
    H:=\sum_i\mathcal{H}(\nabla_{U_i}U_i).
$$
It is a basic vector field on $M$; see \cite[Section~3]{BL23}.
It is related to $L$ and $\tau$ by the following formulas (see \cite[Lemma~2.3]{BL23}):
$$
    \mathrm{Tr}\,L_X=-\ip{H}{X}
    , \quad
    \mathrm{Tr}\,\tau=-\mathrm{div}_B H,
    \quad X\in \mathcal{H}.
$$

We now have all the ingredients to state O'Neill's vertical Ricci formula.

\begin{prop}[{\cite[Theorem~2.2]{BL23}}]\label{prop:oneill-vertical-ricci}
Let $M$ be a complete Riemannian manifold and let $\mathsf{G}$ be a Lie group acting on $M$ freely, properly and isometrically.
Then
$$
    \mathrm{Ric}^{\mathcal{V}}=\mathrm{Ric}^{\mathsf{G}}+\mathrm{A}^* \mathrm{A}+L_H-\tau.
$$
\end{prop}

\subsection{The global trivialization of the vertical bundle}\label{subsec:trivialization}
\textcolor{white}{aaa}

Every $U\in\mathfrak{g}$ induces a Killing vector field $U^*$ on $M$, given by
$$
    U_x^*:=\left.\frac{d}{dt}\right|_{t=0}\bigl(\exp(tU)\cdot x\bigr),
    \qquad x\in M.
$$
We call $U^*$ a \emph{fundamental vector field}.
Since the action is on the left, we have $[U,V]^*=-[U^*,V^*]$ for every $U,V\in\mathfrak{g}$.

The map $\theta:M\times\mathfrak{g}\rightarrow\mathcal{V}$ defined by $\theta(x,U):=U_x^*$ gives a global trivialization of the vertical bundle.
We denote by $\theta_x:\mathfrak{g}\rightarrow\mathcal{V}_x$ the isomorphism $U\mapsto U_x^*$.
Using this trivialization, we can take coordinates of vertical tensors.
We call them \emph{Killing coordinates} and denote them by a left superscript $\mathfrak{g}$.
The Killing coordinates of a vertical vector field $U$ form a smooth map ${}^{\mathfrak{g}}U:M\rightarrow\mathfrak{g}$, while those of a vertical $(1,1)$-tensor $E\in\Gamma(\mathrm{End}\,\mathcal{V})$ form a smooth map ${}^{\mathfrak{g}}E:M\rightarrow\mathrm{End}(\mathfrak{g})$.
Explicitly,
$$
    {}^{\mathfrak{g}}U_x:=\theta_x^{-1}(U_x)
    \quad ; \quad
    {}^{\mathfrak{g}}E_x:=\theta_x^{-1}\circ E_x\circ\theta_x.
$$
Thus, vertical vector fields correspond bijectively to smooth maps $M\rightarrow\mathfrak{g}$, and vertical $(1,1)$-tensors correspond bijectively to smooth maps $M\rightarrow\mathrm{End}(\mathfrak{g})$.

\begin{lemma}[{see e.g. \cite[Appendix~A]{bohm2022euclidean}}]\label{lem:fundamental-field-equivariance}
For every $g\in \mathsf{G}$ and every $x\in M$, we have
$$
    \theta_{g\cdot x}^{-1}\cdot d(m_g)_x\cdot\theta_x=\mathrm{Ad}(g).
$$
Thus, for fixed $g$, the Killing coordinates of $dm_g$ are constant on $M$ and equal to $\mathrm{Ad}(g)$.
\end{lemma}

\begin{remark}\label{rem:killing-coordinate-equivariance}
A vertical vector field $U$ on $M$ is invariant if and only if
$$
    {}^{\mathfrak{g}}U_{g\cdot x}=\mathrm{Ad}(g)({}^{\mathfrak{g}}U_x)
$$
for every $g\in \mathsf{G}$ and every $x\in M$.
Similarly, a vertical $(1,1)$-tensor $E$ on $M$ is equivariant if and only if
$$
    {}^{\mathfrak{g}}E_{g\cdot x}=\mathrm{Ad}(g)\cdot{}^{\mathfrak{g}}E_x\cdot\mathrm{Ad}(g)^{-1}
$$
for every $g\in \mathsf{G}$ and every $x\in M$.
\end{remark}

We can also take Killing coordinates of the metric of each orbit.
Let $\mathrm{Sym}^2_+(\mathfrak{g}^*)$ be the symmetric space of inner products on $\mathfrak{g}$; see Appendix~\ref{app:inner-products}.
We define $h:M\rightarrow\mathrm{Sym}^2_+(\mathfrak{g}^*)$ by
$$
    h_x(U,V):=\langle U^*,V^*\rangle_x,
$$
so that $\theta_x:(\mathfrak{g},h_x)\rightarrow(\mathcal{V}_x,\langle\cdot,\cdot\rangle_x)$ is a linear isometry.
The orbit map $\mathsf{G}\rightarrow \mathsf{G}\cdot x$, $g\mapsto g\cdot x$, is therefore an isometry when $\mathsf{G}$ is equipped with the left-invariant metric determined by $h_x$.
For the action of $\mathsf{GL}(\mathfrak{g})$ on $\mathrm{Sym}^2_+(\mathfrak{g}^*)$ defined in Appendix~\ref{app:inner-products}, we have the following result (see {\cite[Lemma~5.3]{BL23}}).

\begin{lemma}\label{lem:orbit-metric-equivariance}
The map $h:M\rightarrow\mathrm{Sym}^2_+(\mathfrak{g}^*)$ is equivariant:
$$
    h_{g\cdot x}=\mathrm{Ad}(g)\cdot h_x, \quad \forall x\in M,\ \forall g\in \mathsf{G}.
$$
\end{lemma}

The following shows that $L_X$ describes the derivative of $h$ in the horizontal direction $X$ (see {\cite[Lemma~5.4]{BL23}}).

\begin{lemma}\label{lem:orbit-metric-differential}
Let $X$ be a horizontal vector field on $M$.
Then,
$$
    dh_x(X_x)=2h_x({}^{\mathfrak{g}}L_{X_x}\cdot,\cdot), \quad \forall x\in M.
$$
Using the notation of Appendix~\ref{app:inner-products}, this is equivalent to $dh_x(X_x)=-\rho({}^{\mathfrak{g}}L_{X_x})h_x.$
\end{lemma}

\subsection{Holonomy}\label{subsec:holonomy}
\textcolor{white}{aaa}

Consider the principal bundle $\pi:M\rightarrow M/\mathsf{G}$ with the connection given by the horizontal distribution $\mathcal{H}=\mathcal{V}^{\perp}$.
For each $x\in M$, we define the \emph{holonomy group} $\mathsf{Hol}_x\subseteq \mathsf{G}$ as the set of elements $g\in \mathsf{G}$ such that $x$ and $g\cdot x$ can be joined by a piecewise smooth horizontal curve in $M$.
It is an immersed Lie subgroup of $\mathsf{G}$; see \cite[Chapter~II, Theorem~4.2]{KobayashiNomizu63}.

The following observations show that changing the point $x$ changes the holonomy group only by conjugation:
\begin{enumerate}
\renewcommand{\labelenumi}{(\roman{enumi})}
    \item $\mathsf{Hol}_{g\cdot x}=g\,\mathsf{Hol}_x\,g^{-1}$ for every $g\in \mathsf{G}$.
    \item $\mathsf{Hol}_{\sigma(t)}=\mathsf{Hol}_{\sigma(0)}$ for every piecewise smooth horizontal curve $\sigma:[0,1]\rightarrow M$ and every $t\in[0,1]$.
\end{enumerate}
The first identity follows by translating horizontal curves by $g$, while the second follows by concatenating horizontal curves and their translates.
Moreover, for every two points $x,y\in M$, there exist $g\in \mathsf{G}$ and a piecewise smooth horizontal curve $\sigma:[0,1]\rightarrow M$ such that $\sigma(0)=x$ and $g\cdot\sigma(1)=y$.
Indeed, one can horizontally lift a curve joining $\pi(x)$ to $\pi(y)$, starting at $x$; completeness of $M$ ensures that the lift exists on the whole interval; see Proposition~\ref{prop:horizontal-lifting}.
Thus, the groups $\mathsf{Hol}_x$, $x\in M$, are pairwise conjugate.

The Ambrose--Singer theorem \cite{MR63739} describes the Lie algebra of $\mathsf{Hol}_x$ in terms of the O'Neill tensor $\mathrm{A}$; see also \cite[Chapter~II, Theorem~8.1]{KobayashiNomizu63}.
Let $P_x\subseteq M$ be the set of points that can be joined to $x$ by a piecewise smooth horizontal curve.
Then
\begin{equation}\label{eq:holonomy-algebra}
    \mathrm{Lie}(\mathsf{Hol}_x)
    =
    \operatorname{span}_{\mathbb{R}}\left\{
    {}^{\mathfrak{g}}\bigl(\mathrm{A}_{X_p}Y_p\bigr)
    :
    p\in P_x,\quad X_p,Y_p\in\mathcal{H}_p
    \right\},
\end{equation}
This is the usual curvature formulation of the theorem, expressed in Killing coordinates using $\mathcal{V}[X,Y]=2\mathrm{A}_XY$.
In particular, we obtain the following consequence.

\begin{prop}\label{prop:discrete-holonomy}
The O'Neill tensor $\mathrm{A}$ vanishes if and only if
$\mathrm{Lie}(\mathsf{Hol}_x)=\{0\}$ for one, and hence every, $x\in M$.
\end{prop}

We now prove Proposition \ref{prop:splitting-criterion-outline} discussed in the Introduction, which is the global version of the local isometric trivialization criterion given in Proposition~\ref{prop:local-isometric-trivialization}.

\begin{proof}[Proof of Proposition~\ref{prop:splitting-criterion-outline}]
Suppose first that the bundle admits a global isometric trivialization $\mathsf{G}\times(M/\mathsf{G})\to M$.
The orbits correspond to $\mathsf{G}\times\{b\}$ for $b\in M/\mathsf{G}$, which are totally geodesic submanifolds.
The submanifolds $\{g\}\times(M/\mathsf{G})$ are horizontal sections meeting every orbit exactly once, so the holonomy is trivial.

Conversely, suppose that the holonomy is trivial and $L=0$.
Then $\mathrm{A}=0$ by Proposition~\ref{prop:discrete-holonomy}, so let $\Sigma$ be the horizontal leaf through $x_0\in M$.
Horizontal lifting shows that $\Sigma$ meets every orbit, and $x,g\cdot x\in\Sigma$ implies $g\in\mathsf{Hol}_x=\{e\}$.
Thus, $\pi|_\Sigma:\Sigma\rightarrow M/\mathsf{G}$ is a bijective local isometry, and hence an isometry.

Freeness and the unique orbit intersections give a smooth equivariant bijection
$$
    \Phi:(\mathsf{G}\cdot x_0)\times\Sigma\rightarrow M,
    \qquad
    \Phi(g\cdot x_0,x):=g\cdot x.
$$
By Lemma~\ref{lem:orbit-metric-differential}, $L=0$ gives $h_x=h_{x_0}$ on $\Sigma$.
Hence, $\Phi$ preserves the metric on the first factor at $x_0$, and therefore everywhere by equivariance.
On the second factor, it is the restriction of the isometry $m_g$.
The two factors have vertical and horizontal images, which are orthogonal.
Thus, $\Phi$ is an isometry, giving the required trivialization after identifying $\Sigma$ with $M/\mathsf{G}$.
\end{proof}

\subsection{Divergence theorem and maximum principle for cocompact symmetries}\label{subsec:divergence-maximum-principle}
\textcolor{white}{aaa}

Invariant smooth functions $f:M\rightarrow\mathbb{R}$ correspond bijectively to smooth functions $\bar{f}:B\rightarrow\mathbb{R}$.
By Proposition~\ref{prop:horizontal-lifting}, basic vector fields $X$ on $M$ correspond bijectively to vector fields $\bar{X}$ on $B=M/\mathsf{G}$.
Thus, $\mathrm{div}\bar{X}:B\rightarrow\mathbb{R}$ determines an invariant function on $M$, which we denote by $\mathrm{div}_BX$.

For every basic vector field $X$, we have (see \cite[Lemma~B.3]{BL23}):
$$
    \mathrm{div}X=\mathrm{div}_BX-\ip{X}{H}.
$$
We can describe $\mathrm{div}_BX$ locally as follows.
If $\{\bar{X}_1,\ldots,\bar{X}_n\}$ is a local orthonormal frame on $B$, its horizontal lifts $\{X_1,\ldots,X_n\}$ form a local basic orthonormal frame for the horizontal distribution.
For this frame, we have
$$
    \mathrm{div}_BX=\sum_i\ip{\nabla_{X_i}X}{X_i}.
$$

If $f:M\rightarrow\mathbb{R}$ is an invariant smooth function, then $\nabla f$ is basic, and we define
$$
    \Delta_Bf:=\mathrm{div}_B(\nabla f).
$$
This is the invariant function on $M$ corresponding to $\Delta\bar{f}$ on $B$.

Proposition~\ref{prop:compact-divergence-maximum} on the compact quotient $M/\mathsf{G}$ then reads as follows.
\begin{prop}\label{prop:cocompact-maximum-principle}
Let $M$ be a complete connected Riemannian manifold and let $\mathsf{G}$ be a connected Lie group acting on $M$ freely, properly, isometrically and cocompactly.
Let $X$ be a basic vector field on $M$ and let $f:M\rightarrow\mathbb{R}$ be an invariant smooth function.
Then:
\begin{enumerate}
\renewcommand{\labelenumi}{(\roman{enumi})}
    \item If $\mathrm{div}_BX-\ip{X}{\nabla f}\geq0$, then equality holds.
    \item If $\Delta_Bf-\ip{\nabla f}{X}\geq0$, then $f$ is constant.
\end{enumerate}
\end{prop}
Part~(i) follows by applying the divergence theorem on $M/\mathsf{G}$ to the vector field induced by $e^{-f}X$.
Part~(ii) follows from the strong maximum principle \cite[Theorem~3.5]{GilbargTrudinger01}.

\section{Einstein manifolds with cocompact symmetries}\label{sec:einstein-cocompact-symmetries}

In this section, we review some key concepts and results from \cite{BL23}.
The authors work in the following setting: $M$ is a complete connected negative-Einstein manifold, and $\mathsf{G}$ is a connected Lie group acting on $M$ effectively, properly, isometrically and cocompactly, with a single orbit type.
We first discuss these assumptions in Subsection~\ref{subsec:properties-isometric-actions}.
This will be useful later in Section~\ref{sec:general-splitting}.

For the remainder of this section, starting with Subsection~\ref{subsec:nilradical-action}, we work in the setting above and, to simplify the exposition, additionally assume that the $\mathsf{G}$-action is free.
Under this assumption, both effectiveness and the single-orbit-type condition are automatic.
In Sections~\ref{sec:shape-commutativity}--\ref{sec:solvable-splitting}, we will work with free, proper, isometric and cocompact actions, so this additional assumption entails no significant loss for our purposes.

We focus on two main aspects of the analysis carried out in \cite{BL23}.
The first concerns the action of the connected nilradical $\mathsf{N}\subseteq \mathsf{G}$, that is, the maximal connected normal nilpotent subgroup of $\mathsf{G}$.
One of the crucial results established in \cite[Theorem~8.1]{BL23} is that this action is polar.
The second concerns the construction of an $\mathsf{N}$-vertical $(1,1)$-tensor $\beta^+$ using real geometric invariant theory (GIT), following Lauret's work \cite{standard}; see \cite[Sections~4--5 and Appendix~D]{BL23}.
Its relation with the geometry of the $\mathsf{N}$-orbits is described in \cite[Theorem~8.1 and Proposition~9.7]{BL23} and will be discussed below; see Theorem~\ref{thm:beta-properties}.
\vskip2mm

{\bf{\noindent Convention.}}
Throughout this paper, whenever it is ambiguous whether a statement or piece of terminology refers to the $\mathsf{G}$-action or the $\mathsf{N}$-action, it will be understood to refer to the $\mathsf{G}$-action unless explicitly stated otherwise.

\subsection{Properties of isometric actions}\label{subsec:properties-isometric-actions}
\textcolor{white}{aaa}

In this subsection, we consider a Lie group $\mathsf{G}$ acting smoothly and isometrically on a complete connected Riemannian manifold $M$ and discuss the assumptions on the action made in \cite{BL23}.

\subsubsection*{Effectiveness.}\phantomsection\makeatletter\def\@currentlabel{Effectiveness}\makeatother\label{sssec:effectiveness}
The action is \emph{effective} if the only element of $\mathsf{G}$ fixing every point of $M$ is the identity.
Equivalently, the Lie group homomorphism $\mathsf{G}\rightarrow\mathsf{Isom}(M)$, $g\mapsto m_g$, is injective.
We may then identify $\mathsf{G}$ with an immersed Lie subgroup of $\mathsf{Isom}(M)$.

\subsubsection*{Properness.}\phantomsection\makeatletter\def\@currentlabel{Properness}\makeatother\label{sssec:properness}
The action is \emph{proper} if the map $\mathsf{G}\times M\rightarrow M\times M$, $(g,x)\mapsto(x,g\cdot x)$, is proper.
For an effective isometric action, after identifying $\mathsf{G}$ with an immersed Lie subgroup of $\mathsf{Isom}(M)$, properness is equivalent to $\mathsf{G}$ being a closed subgroup; see \cite[Lemma~A.1]{bohm2025non}.
Properness guarantees that the orbits $\mathsf{G}\cdot x\subseteq M$ are closed submanifolds and that the isotropy groups $\mathsf{G}_x\subseteq \mathsf{G}$ are compact.
The following lemma gives a useful converse.

\begin{lemma}\label{lem:closed-orbit-compact-isotropy-proper-action}
Let $M$ be a complete connected Riemannian manifold and let $\mathsf{G}$ be a Lie group acting isometrically on $M$.
Suppose there exists $x_0\in M$ such that $\mathsf{G}\cdot x_0$ is closed and $\mathsf{G}_{x_0}$ is compact.
Then the action is proper.
\end{lemma}

\begin{proof}
Since $\mathsf{G}\cdot x_0$ is closed, it is an embedded submanifold by \cite[Theorem~2.1.12]{PiccioneTausk2000}.
Hence the natural map $\bar{\phi}:\mathsf{G}/\mathsf{G}_{x_0}\rightarrow \mathsf{G}\cdot x_0$ is a diffeomorphism.
Since $\mathsf{G}_{x_0}$ is compact, the quotient projection $q:\mathsf{G}\rightarrow \mathsf{G}/\mathsf{G}_{x_0}$ is proper by \cite[Lemma~3.2]{KramerMuratGarcia2026}.
Thus the orbit map $\phi=\bar{\phi}\circ q:\mathsf{G}\rightarrow \mathsf{G}\cdot x_0$ is proper.

We verify the sequential criterion for properness.
Take sequences $g_i\in \mathsf{G}$ and $x_i\in M$ such that $x_i\rightarrow x$ and $g_i\cdot x_i\rightarrow y$.
Since the action is isometric,
$$
    d(g_i\cdot x_0,x_0)
    \leq d(g_i\cdot x_0,g_i\cdot x_i)+d(g_i\cdot x_i,x_0)
    =d(x_0,x_i)+d(g_i\cdot x_i,x_0).
$$
Thus $g_i\cdot x_0$ is bounded.
By completeness and the Hopf--Rinow theorem, it has a convergent subsequence.
Passing to this subsequence, write $g_i\cdot x_0\rightarrow z$.
Since $\mathsf{G}\cdot x_0$ is closed, $z\in \mathsf{G}\cdot x_0$.
The set
$$
    K:=\{g_i\cdot x_0:i\in\mathbb{N}\}\cup\{z\}
$$
is therefore a compact subset of $\mathsf{G}\cdot x_0$.
Properness of $\phi$ implies that $\phi^{-1}(K)$ is compact.
Since $g_i\in\phi^{-1}(K)$, the sequence $g_i$ has a convergent subsequence in $\mathsf{G}$.
Therefore the action is proper.
\end{proof}

The following result will be used many times throughout the paper.

\begin{lemma}\label{lem:simply-connected-solvable-free-action}
A connected, simply connected, solvable Lie group $\mathsf{S}$ is diffeomorphic to a Euclidean space and has no nontrivial compact subgroups.
Consequently, every proper smooth action of such a group on a manifold is free.
\end{lemma}
\begin{proof}
By \cite[Corollary~1.126]{Knapp02},
$\mathsf{S}$ is diffeomorphic to $\mathbb{R}^n$ and is therefore contractible.
By the Cartan--Iwasawa--Malcev theorem, every compact subgroup is contained
in a maximal compact subgroup which is homotopy equivalent to $\mathsf{S}$.
By Poincaré duality, a compact connected manifold is contractible
if and only if it is a point.
Hence, $\mathsf{S}$ has no nontrivial compact subgroups.
Finally, the stabilizers of a proper action are compact
and therefore trivial, so every proper smooth action
of $\mathsf{S}$ is free.
\end{proof}

\subsubsection*{Principal orbits.}\phantomsection\makeatletter\def\@currentlabel{Principal orbits}\makeatother\label{sssec:principal-orbits}
An action has a \emph{single orbit type} if the isotropy groups $\mathsf{G}_x$, $x\in M$, are pairwise conjugate in $\mathsf{G}$.
For a proper action on a connected manifold, this is equivalent to every orbit being principal; see \cite[Theorem~3.82]{AB15}.

If the action is proper, isometric and has a single orbit type, then $M/\mathsf{G}$ admits a natural smooth structure and a unique Riemannian metric for which the quotient map $M\rightarrow M/\mathsf{G}$ is a Riemannian submersion; see \cite[Theorem~3.95 and Remark~3.99]{AB15}.

\subsubsection*{Cocompactness.}\phantomsection\makeatletter\def\@currentlabel{Cocompactness}\makeatother\label{sssec:cocompactness}
An action of a Lie group $\mathsf{G}$ on a smooth manifold $M$ is \emph{cocompact} if the orbit space $M/\mathsf{G}$ is compact.
Equivalently, there exists a compact subset $K\subseteq M$ such that $M=\mathsf{G}\cdot K$.

A Riemannian manifold $M$ has \emph{cocompact symmetries} if it admits a cocompact isometric action by a connected Lie group.
Equivalently, $M/\mathsf{Isom}(M)^o$ is compact.
Every Riemannian manifold with cocompact symmetries is complete.

\begin{lemma}\label{lem:cocompact-subgroup}
Let $\mathsf{G}$ be a Lie group acting smoothly, properly and cocompactly on a smooth manifold $M$, and let $\mathsf{H}\subseteq \mathsf{G}$ be a closed subgroup.
Then the $\mathsf{H}$-action on $M$ is cocompact if and only if $\mathsf{G}/\mathsf{H}$ is compact.
\end{lemma}

\begin{proof}
Choose a compact subset $K\subseteq M$ with $M=\mathsf{G}\cdot K$.
If $\mathsf{G}/\mathsf{H}$ is compact, then so is $\mathsf{H}\backslash \mathsf{G}$.
Consequently, there exists a compact subset $C\subseteq \mathsf{G}$ such that $\mathsf{G}=\mathsf{H}C$.
It follows that $M=\mathsf{H}\cdot(C\cdot K)$.
The set $C\cdot K$ is compact, being the image of $C\times K$ under the action map.
Hence the $\mathsf{H}$-action is cocompact.

Conversely, suppose that $M=\mathsf{H}\cdot K'$ for a compact subset $K'\subseteq M$.
Fix $x_0\in M$ and set
$$
    C:=\{g\in \mathsf{G}:g\cdot x_0\in K'\}.
$$
Properness of the $\mathsf{G}$-action implies that $C$ is compact.
For every $g\in \mathsf{G}$, there exist $h\in \mathsf{H}$ and $k\in K'$ such that $g\cdot x_0=h\cdot k$.
Thus $h^{-1}g\in C$, and therefore $\mathsf{G}=\mathsf{H}C$.
This proves that $\mathsf{H}\backslash \mathsf{G}$, and hence $\mathsf{G}/\mathsf{H}$, is compact.
\end{proof}

\subsection{The action of the nilradical}\label{subsec:nilradical-action}
\textcolor{white}{aaa}

For the remainder of this section, let $M$ be a complete connected Riemannian manifold with $\mathrm{ric}_{\mathrm{g}}=-\mathrm{g}$, and let $\mathsf{G}$ be a connected Lie group acting freely, properly, isometrically and cocompactly on $M$.
We review results from \cite{BL23} on the action of the connected nilradical $\mathsf{N}$, which is closed in $\mathsf{G}$, and on the associated Riemannian submersion $M\rightarrow M/\mathsf{N}$.

\begin{notation}\phantomsection\makeatletter\def\@currentlabel{Notation}\makeatother\label{not:nilradical-tensors}
We use a hat to distinguish the tensors and distributions associated with the Riemannian submersion $M\rightarrow M/\mathsf{N}$ from those associated with $M\rightarrow M/\mathsf{G}$.
Thus, its O'Neill tensors are denoted by $\widehat{\mathrm{A}}$ and $\hat{L}$, and the mean curvature vector field of the $\mathsf{N}$-orbits by $\widehat{H}$.
We denote the $\mathsf{N}$-vertical and $\mathsf{N}$-horizontal distributions by $\widehat{\mathcal{V}}$ and $\widehat{\mathcal{H}}$, respectively.
Accordingly, $\mathsf{N}$-vertical vector fields are denoted by $\widehat{U},\widehat{V},\widehat{W}$, and $\mathsf{N}$-horizontal vector fields by $\widehat{X},\widehat{Y},\widehat{Z}$.
The map $\hat{h}:M\rightarrow\mathrm{Sym}^2_+(\mathfrak{n}^*)$ denotes the Killing coordinates of the metric of the $\mathsf{N}$-orbits.
\end{notation}
\begin{remark}\label{rem:preservation-n-horizontal}
The $\mathsf{G}$-action preserves the $\mathsf{N}$-vertical and $\mathsf{N}$-horizontal distributions.
More precisely, recalling that $m_g(x)=g\cdot x$, we have for every $g\in\mathsf{G}$ and $x\in M$:
$$
    d(m_g)_x\bigl(\widehat{\mathcal{V}}_x\bigr)=\widehat{\mathcal{V}}_{g\cdot x},
    \qquad
    d(m_g)_x\bigl(\widehat{\mathcal{H}}_x\bigr)=\widehat{\mathcal{H}}_{g\cdot x}.
$$
Indeed, normality of $\mathsf{N}$ implies that $m_g$ maps $\mathsf{N}\cdot x$ diffeomorphically onto $\mathsf{N}\cdot(g\cdot x)$, proving the first equality.
Moreover, $m_g$ is an isometry, so its differential also maps the orthogonal complement $\widehat{\mathcal{H}}_x$ of $\widehat{\mathcal{V}}_x$ onto the orthogonal complement $\widehat{\mathcal{H}}_{g\cdot x}$ of $\widehat{\mathcal{V}}_{g\cdot x}$.
\end{remark}

A central result of \cite{BL23} states that the $\mathsf{N}$-horizontal distribution is integrable, or equivalently, that the corresponding O'Neill tensor $\widehat{\mathrm{A}}$ vanishes.

\begin{theorem}[{\cite[Theorem~8.1(i)]{BL23}}]\label{thm:nilradical-polarity}
Let $M$ be a complete connected negative-Einstein manifold admitting a free, proper, isometric and cocompact action of a connected Lie group $\mathsf{G}$, and let $\mathsf{N}$ denote the nilradical of $\mathsf{G}$.
Then the horizontal distribution $\widehat{\mathcal{H}}$ associated with the Riemannian submersion $M\rightarrow M/\mathsf{N}$ is integrable.
Equivalently, $\widehat{\mathrm{A}}=0$.
\end{theorem}

The following consequence gives an algebraic interpretation of Theorem~\ref{thm:nilradical-polarity}.
In the solvable case, it can be regarded as a cocompact analogue of Lauret's standardness theorem \cite[Theorem~1]{standard}.
Recall that, for each $x\in M$, $h_x$ denotes the inner product on $\mathfrak{g}$ induced by the metric on the $\mathsf{G}$-orbit through $x$; see Subsection~\ref{subsec:trivialization}.

\begin{corollary}[{\cite[Corollary~8.7]{BL23}}]\label{cor:horizontal-complement-subalgebra}
Under the assumptions of Theorem~\ref{thm:nilradical-polarity}, let $\mathfrak{a}_x\subseteq\mathfrak{g}$ be the $h_x$-orthogonal complement of $\mathfrak{n}$ for each $x\in M$.
Then $\mathfrak{a}_x$ is a Lie subalgebra of $\mathfrak{g}$.
Moreover, the map $x\mapsto\mathfrak{a}_x$ is constant along $\mathsf{N}$-horizontal curves.
\end{corollary}

\begin{remark}\label{rem:horizontal-complement-equivariance}
We record two additional observations concerning Corollary~\ref{cor:horizontal-complement-subalgebra}.
\begin{enumerate}
\renewcommand{\labelenumi}{(\roman{enumi})}
    \item The map $x\mapsto\mathfrak{a}_x$ is $\mathsf{G}$-equivariant: $\mathfrak{a}_{g\cdot x}=\mathrm{Ad}(g)\cdot\mathfrak{a}_x$ for every $g\in \mathsf{G}$ and $x\in M$.
    \item If $\mathfrak{g}$ is solvable, then each $\mathfrak{a}_x$ is an abelian subalgebra.
    Thus, for $A\in\mathfrak{a}_x$, the endomorphism $\mathrm{ad}_A$ has the following matrix form with respect to $\mathfrak{g}=\mathfrak{a}_x\oplus\mathfrak{n}$:
$$
    \mathrm{ad}_A=
    \begin{pmatrix}
        0&0\\
        0&\mathrm{ad}_A|_{\mathfrak{n}}
    \end{pmatrix}.
$$
\end{enumerate}
\end{remark}

Corollary~\ref{cor:horizontal-complement-subalgebra} naturally leads to a refinement of the vertical-horizontal decomposition associated with the $\mathsf{G}$-action.
Define the distribution $\mathcal{A}\subseteq TM$ by
$$
    \mathcal{A}:=\mathcal{V}\cap\widehat{\mathcal{H}}.
$$
This is the orthogonal complement of $\widehat{\mathcal{V}}$ inside $\mathcal{V}$, so $\mathcal{V}=\mathcal{A}\oplus\widehat{\mathcal{V}}$.
Under the global trivialization of the $\mathsf{G}$-vertical bundle from Subsection~\ref{subsec:trivialization}, a vector $U_x\in\mathcal{V}_x$ belongs to $\mathcal{A}_x$ if and only if ${}^{\mathfrak{g}}U_x\in\mathfrak{a}_x$.
Combining this with $TM=\mathcal{H}\oplus\mathcal{V}$, we obtain the orthogonal decomposition
$$
    TM=\mathcal{H}\oplus\mathcal{A}\oplus\widehat{\mathcal{V}}.
$$
By Remark~\ref{rem:preservation-n-horizontal}, the distribution $\mathcal{A}$ is $\mathsf{G}$-invariant:
\begin{equation}\label{eq:a-invariance}
    d(m_g)_x(\mathcal{A}_x)=\mathcal{A}_{g\cdot x}
    \qquad\text{for every }g\in\mathsf{G}\text{ and }x\in M.
\end{equation}

\begin{remark}\label{rem:intermediate-quotient}
Since both $\mathcal{V}$ and $\widehat{\mathcal{H}}$ are integrable, the distribution $\mathcal{A}$ is integrable as well.
In general, however, this does not produce a Riemannian submersion onto the corresponding leaf space: the quotient $M/\mathcal{A}$ need not carry a suitable smooth structure, and the transverse metric need not descend to the quotient.

There is nevertheless a natural substitute.
Since $\mathsf{N}$ is normal in $\mathsf{G}$, the quotient group $\mathsf{G}/\mathsf{N}$ acts isometrically on $M/\mathsf{N}$, and the corresponding orbit projection gives the Riemannian submersion $M/\mathsf{N}\rightarrow M/\mathsf{G}$.
The differential of the projection $M\rightarrow M/\mathsf{N}$ maps $\mathcal{A}$ onto the vertical distribution of this submersion, and $\mathcal{H}$ onto its horizontal distribution.
\end{remark}

Another useful consequence of Theorem~\ref{thm:nilradical-polarity} is the following block decomposition of several tensors associated with the $\mathsf{G}$-submersion.

\begin{lemma}\label{lem:vertical-tensor-blocks}
Let $X$ be a $\mathsf{G}$-horizontal vector field.
With respect to the orthogonal decomposition $\mathcal{V}=\mathcal{A}\oplus\widehat{\mathcal{V}}$, we have
$$
    L_X=\begin{pmatrix} *&0\\0&\hat{L}_X\end{pmatrix},
    \quad
    \mathrm{A}^*\mathrm{A}=\begin{pmatrix} *&0\\0&0\end{pmatrix}.
$$
\end{lemma}

\begin{proof}
Since $L_X$ and $\mathrm{A}^*\mathrm{A}$ are pointwise self-adjoint endomorphisms of $\mathcal{V}$, it suffices to prove that
$$
    L_X|_{\widehat{\mathcal{V}}}=\hat{L}_X
    \qquad\text{and}\qquad
    \mathrm{A}^*\mathrm{A}|_{\widehat{\mathcal{V}}}=0.
$$
We first prove the second identity.
By Theorem~\ref{thm:nilradical-polarity}, we have $\mathcal{H}(\nabla_X\widehat{U})=0$, and hence $\mathrm{A}_X\widehat{U}=0$.
Therefore, for every vertical vector field $U$ and a local orthonormal horizontal frame $\{X_i\}$,
$$
    \ip{\mathrm{A}^*\mathrm{A}(\widehat{U})}{U}
    =\sum_i\ip{\mathrm{A}_{X_i}\widehat{U}}{\mathrm{A}_{X_i}U}
    =0.
$$
It follows that $\mathrm{A}^*\mathrm{A}(\widehat{U})=0$.

We now prove the first identity.
Since both $L$ and $\hat{L}$ are tensorial, it is enough to argue at a fixed point, where we may extend $\widehat{U}$ to an $\mathsf{N}$-vertical Killing field and $X$ to a $\mathsf{G}$-basic vector field.
By Lemma~\ref{lem:basic-vertical-bracket}, $[\widehat{U},X]=0$.
Therefore,
$$
    L_X\widehat{U}
    =\mathcal{V}(\nabla_{\widehat{U}}X)
    =\mathcal{V}(\nabla_X\widehat{U})
    =\widehat{\mathcal{V}}(\nabla_X\widehat{U})
    =\widehat{\mathcal{V}}(\nabla_{\widehat{U}}X)
    =\hat{L}_X\widehat{U}.
$$
In the third equality, we have used Theorem~\ref{thm:nilradical-polarity}, which implies that $\nabla_X\widehat{U}$ is $\mathsf{N}$-vertical.
\end{proof}

We now present some important results stated in \cite{bohm2025non}, whose proofs are mostly contained in \cite{BL23}.
We denote by $\hat{\tau}\in\Gamma(\mathrm{End}\,\widehat{\mathcal{V}})$ the tensor $\tau$ from Subsection~\ref{subsec:oneill-curvature} associated with the $\mathsf{N}$-action.
Let $\widehat{\Sigma}\subseteq M$ be an $\mathsf{N}$-horizontal leaf, which exists by Theorem~\ref{thm:nilradical-polarity}.
The projection $\widehat{\Sigma}\rightarrow M/\mathsf{N}$ is a local isometry, so the Ricci curvature of $\widehat{\Sigma}$ agrees with that of $M/\mathsf{N}$ under the identification $\widehat{\mathcal{H}}\equiv T(M/\mathsf{N})$.

\begin{theorem}[{\cite[Theorem~2.2(iii)--(iv) and proof of (iii)]{bohm2025non}}]\label{thm:nilradical-metric-identities}
Under the assumptions of this section, the following hold:
\begin{enumerate}
\renewcommand{\labelenumi}{(\roman{enumi})}
    \item $\hat{\tau}=0$.
    \item The map
$$
    \hat{h}|_{\widehat{\Sigma}}:\widehat{\Sigma}\longrightarrow\mathrm{Sym}^2_+(\mathfrak{n}^*)
$$
    is harmonic, where $\widehat{\Sigma}$ carries the induced metric and $\mathrm{Sym}^2_+(\mathfrak{n}^*)$ carries the symmetric-space metric described in Appendix~\ref{app:inner-products}.
    \item The Ricci endomorphism of $M/\mathsf{N}$ satisfies
$$
    \mathrm{Ric}^{M/\mathsf{N}}=-\mathrm{Id}_{\widehat{\mathcal{H}}}+\hat{L}^*\hat{L}
$$
    where the section $\hat{L}^*\hat{L}$ of $\mathrm{End}\,\widehat{\mathcal{H}}$ is defined by
$$
    \langle\hat{L}^*\hat{L}(\widehat{X}),\widehat{Y}\rangle
    :=\mathrm{Tr}(\hat{L}_{\widehat{X}}\cdot\hat{L}_{\widehat{Y}}) \qquad \text{for all} \ \widehat{X},\widehat{Y}\in\widehat{\mathcal{H}}.
$$
\end{enumerate}
\end{theorem}

Recall that $\widehat{H}$ denotes the mean curvature vector field of the $\mathsf{N}$-orbits in $M$; see Subsection~\ref{subsec:oneill-curvature}.
If an $\mathsf{N}$-orbit is regarded instead as a submanifold of the $\mathsf{G}$-orbit containing it, then its mean curvature vector field is precisely the $\mathsf{G}$-vertical component $\mathcal{V}(\widehat{H})$.
Its expression in Killing coordinates is characterized by the identity
$$
    h_x\bigl({}^{\mathfrak{g}}\mathcal{V}(\widehat{H})_x,U\bigr)
    =\operatorname{Tr}(\operatorname{ad}_U),
    \qquad x\in M,\quad U\in\mathfrak{g};
$$
see \cite[equations~(2.1)--(2.2)]{Heber98}.
In particular, ${}^{\mathfrak{g}}\mathcal{V}(\widehat{H})_x\in\mathfrak{a}_x$ for every $x\in M$, so $\mathcal{V}(\widehat{H})$ is an $\mathcal{A}$-valued vector field.

We conclude the subsection with a consequence of \cite[Lemma~C.2]{BL23}.
For $E\in\mathrm{End}(\mathfrak{g})$, we denote by
$$
    \mathrm{S}_x(E):=\frac12(E+E^{t_x})
$$
its $h_x$-symmetric part, where $E^{t_x}$ denotes the adjoint with respect to $h_x$.
We use the same notation for endomorphisms of $\mathfrak{n}$, with respect to $\hat{h}_x$, and omit the subscript when the point is clear.

\begin{lemma}\label{lem:shape-symmetric-adjoint}
Let $A$ be an $\mathcal{A}$-valued vector field.
Then
$$
    {}^{\mathfrak{n}}\hat{L}_A
    =-\mathrm{S}\bigl(\mathrm{ad}_{{}^{\raisebox{0.15em}{$\scriptstyle\mathfrak{g}$}}\! A}|_{\mathfrak{n}}\bigr).
$$
\end{lemma}

\subsection{Geometric Invariant Theory}\label{subsec:git}

\textcolor{white}{aaa}

The real geometric invariant theory of the change-of-basis action of $\mathsf{GL}(\mathfrak{n})$ on $\Lambda^2(\mathfrak{n}^*)\otimes\mathfrak{n}$ assigns a stratum to each nonzero nilpotent Lie bracket; see \cite[Section~2]{standard}.
For a non-abelian metric nilpotent Lie algebra $(\mathfrak{n},\bar{h})$, the corresponding stratum determines a $\bar{h}$-self-adjoint endomorphism $\bar{\beta}$, with $\mathrm{Tr}\,\bar{\beta}=-1$, once the representative is chosen as in \cite[Lemma~D.2]{BL23} (the \emph{well-gauged} choice).
Lauret proved that the endomorphism
$$
    \bar{\beta}^+:=\frac{\bar{\beta}}{\mathrm{Tr}(\bar{\beta}^2)}+\mathrm{Id}_{\mathfrak{n}}
$$
is positive definite; see \cite[Lemma~2.17]{standard}.

We briefly recall how B\"ohm and Lafuente extend this construction to the $\mathsf{N}$-vertical bundle in their setting; see \cite[Section~4 and Definition~5.5]{BL23}.
Fix a background inner product $\bar{h}$ on $\mathfrak{n}$ and a $\bar{h}$-orthonormal eigenbasis of $\bar{\beta}$, ordered by nondecreasing eigenvalues.
Let $\mathsf{B}_{\bar{\beta}}\subseteq\mathsf{GL}(\mathfrak{n})$ be the group of lower triangular matrices with positive diagonal entries in this basis.
For each $x\in M$, there is a unique $q_x\in\mathsf{B}_{\bar{\beta}}$ such that
$$
    \hat{h}_x=q_x\cdot\bar{h},\qquad
    (q_x\cdot\bar{h})(U,V):=\bar{h}(q_x^{-1}U,q_x^{-1}V).
$$
The map $x\mapsto q_x$ is smooth, and one defines
$$
    {}^{\mathfrak{n}}\beta_x^+:=q_x\bar{\beta}^+q_x^{-1},\qquad
    \beta_x^+:=\hat{\theta}_x\circ{}^{\mathfrak{n}}\beta_x^+\circ\hat{\theta}_x^{-1},
$$
where $\hat{\theta}_x:\mathfrak{n}\rightarrow\widehat{\mathcal{V}}_x$ is the trivialization of the $\mathsf{N}$-vertical bundle from Subsection~\ref{subsec:trivialization}.
Thus, $\beta^+$ is a smooth section of $\mathrm{End}\,\widehat{\mathcal{V}}$.
When $\mathfrak{n}$ is abelian, one simply sets $\bar{\beta}^+=\mathrm{Id}_{\mathfrak{n}}$, ${}^{\mathfrak{n}}\beta_x^+=\mathrm{Id}_{\mathfrak{n}}$ and $\beta_x^+=\mathrm{Id}_{\widehat{\mathcal{V}}_x}$.

The following theorem collects the properties that will be used in this paper; see \cite[Theorem~8.1, Corollary~8.4, Proposition~9.7 and Corollary~9.9]{BL23}.
Here, $\mathrm{Ric}^{\mathsf{N}}$ denotes the Ricci endomorphism of the $\mathsf{N}$-orbits.

\begin{theorem}\label{thm:beta-properties}
Under the assumptions of this section, the section $\beta^+$ of $\mathrm{End}\,\widehat{\mathcal{V}}$ has the following properties:
\begin{enumerate}
\renewcommand{\labelenumi}{(\roman{enumi})}
    \item $\beta^+$ is $\mathsf{G}$-equivariant, pointwise self-adjoint and positive definite.

    \item The mean curvature vector field and the Ricci endomorphism of the $\mathsf{N}$-orbits satisfy
$$
    \hat{L}_{\widehat{H}}=-\beta^+,\qquad
    \mathrm{Ric}^{\mathsf{N}}=-\mathrm{Id}_{\widehat{\mathcal{V}}}+\beta^+.
$$

    \item For every $x\in M$ and every $\mathsf{N}$-horizontal vector $\widehat{X}_x\in\widehat{\mathcal{H}}_x$,
$$
    {}^{\mathfrak{n}}\beta_x^+\in\mathrm{Der}(\mathfrak{n}),\qquad
    {}^{\mathfrak{n}}\hat{L}_{\widehat{X}_x}\in\mathrm{Der}(\mathfrak{n}).
$$

    \item The map ${}^{\mathfrak{n}}\beta^+:M\rightarrow\mathrm{End}(\mathfrak{n})$ is constant along $\mathsf{N}$-horizontal curves.

    \item If $\mathfrak{g}$ is solvable, then
$$
    [{}^{\mathfrak{n}}\beta_x^+,\mathrm{ad}_A|_{\mathfrak{n}}]=0
    \qquad\text{for every }x\in M\text{ and }A\in\mathfrak{a}_x.
$$
    See \cite[equation~(35), in the proof of Corollary~11.7]{BL23}.
\end{enumerate}
\end{theorem}

\begin{remark}\label{rem:beta-equivariance-corrected-parallelism}
Since $\mathsf{N}$ is normal in $\mathsf{G}$, each $m_g$ restricts to an isometry from $\mathsf{N}\cdot x$ to $\mathsf{N}\cdot(g\cdot x)$.
These isometries preserve the Ricci endomorphisms of the $\mathsf{N}$-orbits.
Thus, the identity $\beta^+=\mathrm{Ric}^{\mathsf{N}}+\mathrm{Id}_{\widehat{\mathcal{V}}}$ implies that $\beta^+$ is $\mathsf{G}$-equivariant.

The constancy in~(iv) is established in the proof of \cite[Corollary~8.4]{BL23}.
Although the general differentiation formula in \cite[Lemma~5.7]{BL23} requires correction (see Corollary~\ref{cor:self-adjoint-adjoint-derivative}), the argument establishing this constancy remains valid.
Since $\widehat{\mathrm{A}}=0$ and $[\hat{L}_{\widehat{X}},\beta^+]=0$ by \cite[Theorem~8.1]{BL23}, Proposition~\ref{Adam}, applied to the $\mathsf{N}$-action, then gives $\nabla_{\widehat{X}}\beta^+=0$.
\end{remark}

\section{The adjoint bundle connection}\label{sec:adjoint-connection}

Throughout this section, $M$ is a complete connected Riemannian manifold and $\mathsf{G}$ is a connected Lie group acting freely, properly and isometrically on $M$.
The orbit projection $\pi:M\rightarrow M/\mathsf{G}$ is therefore a principal bundle, equipped with the principal connection given by the horizontal distribution $\mathcal{H}=\mathcal{V}^{\perp}$.

The adjoint representation of $\mathsf{G}$ on $\mathfrak{g}$ determines the associated vector bundle $\mathrm{ad}(M):=(M\times\mathfrak{g})/\mathsf{G}$, where $(x,U)$ is identified with $(g\cdot x,\mathrm{Ad}(g)U)$.
This is called the \emph{adjoint bundle}, and its sections can be identified with invariant vertical vector fields on $M$.
The principal connection $\mathcal{H}$ induces a connection on $\mathrm{ad}(M)$, which we denote by $\nabla^{\mathrm{ad}}$.
Under this identification, covariant differentiation of a section in a direction on $M/\mathsf{G}$ corresponds to differentiating an invariant vertical field along the basic horizontal lift of that direction.

We begin by describing $\nabla^{\mathrm{ad}}$ in Killing coordinates.
This description also defines horizontal differentiation of vertical fields that are not invariant.
We then compare it with the Levi-Civita connection, obtaining formulas involving the shape operator $L$ and the O'Neill tensor $\mathrm{A}$.
These formulas allow us to study differentiation of metric adjoints and, for polar actions, relate the curvature of $M$ to the shape operators of the orbits.
Finally, we introduce \emph{horizontally Killing vector fields} and characterize the existence of invariant such fields in terms of holonomy.

\subsection{The connection in Killing coordinates}\label{subsec:adjoint-connection-coordinates}

\textcolor{white}{aaa}

Recall from Subsection~\ref{subsec:trivialization} that fundamental vector fields define a global trivialization $M\times\mathfrak{g}\rightarrow\mathcal{V}$, $(x,U)\mapsto\theta_x(U)=U_x^*$, where $\theta_x:\mathfrak{g}\rightarrow\mathcal{V}_x$ is an isomorphism for every $x\in M$.
In this trivialization, a vertical vector field $U$ is represented by the smooth map ${}^{\mathfrak{g}}U:M\rightarrow\mathfrak{g}$ given by ${}^{\mathfrak{g}}U_x=\theta_x^{-1}(U_x)$, which we call its \emph{Killing coordinates}.
Likewise, a section $E$ of $\mathrm{End}\,\mathcal{V}\rightarrow M$, called a \emph{vertical $(1,1)$-tensor}, is represented by the smooth map ${}^{\mathfrak{g}}E:M\rightarrow\mathrm{End}(\mathfrak{g})$ given by ${}^{\mathfrak{g}}E_x=\theta_x^{-1}\cdot E_x\cdot\theta_x$.

\begin{definition}\label{def:adjoint-connection}
Let $X$ be a horizontal vector field, $U$ a vertical vector field and $E$ a vertical $(1,1)$-tensor.
We define the vertical vector field $\nabla^{\mathrm{ad}}_XU$ and the vertical $(1,1)$-tensor $\nabla^{\mathrm{ad}}_XE$ by
$$
    (\nabla^{\mathrm{ad}}_XU)_x:=\theta_x\bigl(X_x({}^{\mathfrak{g}}U)\bigr),\qquad
    (\nabla^{\mathrm{ad}}_XE)_x:=\theta_x\cdot X_x({}^{\mathfrak{g}}E)\cdot\theta_x^{-1}.
$$
\end{definition}

Thus, ${}^{\mathfrak{g}}(\nabla^{\mathrm{ad}}_XU)=X({}^{\mathfrak{g}}U)$ and ${}^{\mathfrak{g}}(\nabla^{\mathrm{ad}}_XE)=X({}^{\mathfrak{g}}E)$.
This defines covariant differentiation of arbitrary vertical fields in horizontal directions; neither $X$ nor $U$ is required to be invariant.
We use the notation $\nabla^{\mathrm{ad}}$ because, when $X$ is basic and $U$ is invariant, it agrees with the adjoint bundle connection introduced above; see Remark~\ref{rem:induced-adjoint-connection}.

For vertical $(1,1)$-tensors, the usual compatibility identity follows from the product rule in Killing coordinates.

\begin{lemma}\label{lem:adjoint-connection-product-rule}
Let $E$ be a vertical $(1,1)$-tensor and $X$ a horizontal vector field.
Then, for every vertical vector field $U$,
$$
    (\nabla^{\mathrm{ad}}_XE)(U)=\nabla^{\mathrm{ad}}_X(EU)-E(\nabla^{\mathrm{ad}}_XU).
$$
\end{lemma}
\begin{proof}
Since ${}^{\mathfrak{g}}(EU)={}^{\mathfrak{g}}E({}^{\mathfrak{g}}U)$, the product rule gives
\begin{equation}\label{eq:adjoint-product-rule}
    X_x\bigl({}^{\mathfrak{g}}E({}^{\mathfrak{g}}U)\bigr)
    =X_x({}^{\mathfrak{g}}E)({}^{\mathfrak{g}}U_x)+{}^{\mathfrak{g}}E_x\bigl(X_x({}^{\mathfrak{g}}U)\bigr).
\end{equation}
Applying $\theta_x$ and using Definition~\ref{def:adjoint-connection} proves the claim.
\end{proof}

A vertical $(1,1)$-tensor $E$ extends to a $(1,1)$-tensor on $M$ by setting $E(X)=0$ for every horizontal vector field $X$.
We use the same convention for $L_X$ and $\nabla^{\mathrm{ad}}_XE$.
The following identities will be important (see e.g. \cite[Lemma~2.18]{LT25}).

\begin{prop}\label{Adam}
Let $U$ be a vertical vector field, $E$ a vertical $(1,1)$-tensor and $X$ a horizontal vector field.
Then
$$
    \nabla_XU=\nabla^{\mathrm{ad}}_XU+L_XU+\mathrm{A}_XU,\qquad
    \nabla_XE=\nabla^{\mathrm{ad}}_XE+[L_X,E]+[\mathrm{A}_X,E].
$$
\end{prop}
\begin{proof}
Since both identities are tensorial in $X$, we may assume that $X$ is basic.
Let $\{U_1,\ldots,U_n\}$ be a basis of $\mathfrak{g}$, and let $f_1,\ldots,f_n:M\rightarrow\mathbb{R}$ be the components of ${}^{\mathfrak{g}}U$ in this basis.
By definition,
$$
    \nabla^{\mathrm{ad}}_XU=X(f_1)\cdot U_1^*+\cdots+X(f_n)\cdot U_n^*.
$$
On the other hand,
$$
    \nabla_XU=X(f_1)\cdot U_1^*+\cdots+X(f_n)\cdot U_n^*+f_1\cdot\nabla_XU_1^*+\cdots+f_n\cdot\nabla_XU_n^*.
$$
The first $n$ terms sum to the vertical vector field $\nabla^{\mathrm{ad}}_XU$.
Since $X$ is basic and each $U_i^*$ is vertical and Killing, Lemma~\ref{lem:basic-vertical-bracket} gives $[U_i^*,X]=0$.
Hence,
$$
    \mathcal{V}(\nabla_XU_i^*)=\mathcal{V}(\nabla_{U_i^*}X)=L_X(U_i^*).
$$
It follows that
\begin{align*}
    \mathcal{V}(\nabla_XU)
    &=\nabla^{\mathrm{ad}}_XU+f_1\cdot\mathcal{V}(\nabla_XU_1^*)+\cdots+f_n\cdot\mathcal{V}(\nabla_XU_n^*)\\
    &=\nabla^{\mathrm{ad}}_XU+L_X(f_1\cdot U_1^*+\cdots+f_n\cdot U_n^*)=\nabla^{\mathrm{ad}}_XU+L_XU.
\end{align*}
The horizontal component of $\nabla_XU$ is $\mathrm{A}_XU$ by definition.
This proves the first identity.

We now use this identity and Lemma~\ref{lem:adjoint-connection-product-rule} to prove the formula for vertical $(1,1)$-tensors.
For a vertical vector field $U$,
\begin{align*}
    (\nabla_XE)U
    &=\nabla_X(EU)-E(\mathcal{V}(\nabla_XU))
    =\nabla^{\mathrm{ad}}_X(EU)+L_X(EU)+\mathrm{A}_X(EU)-E(\nabla^{\mathrm{ad}}_XU)-E(L_XU)\\
    &=(\nabla^{\mathrm{ad}}_XE)U+[L_X,E]U+\mathrm{A}_X(EU).
\end{align*}
Since $\mathrm{A}_XU$ is horizontal, $E(\mathrm{A}_XU)=0$, which proves the identity on vertical vectors.
For a horizontal vector field $Y$,
$$
    (\nabla_XE)(Y)=\nabla_X(EY)-E(\mathcal{V}(\nabla_XY))=-E(\mathrm{A}_XY)=[\mathrm{A}_X,E]Y.
$$
Here we have used $EY=0$, and hence $\mathrm{A}_X(EY)=0$.
The remaining terms in the asserted identity vanish on horizontal vectors, completing the proof.
\end{proof}

\begin{remark}\label{rem:adjoint-connection-block-form}
For a vertical vector field $U$, the two components of the first identity in Proposition~\ref{Adam} are
$$
    \mathcal{V}(\nabla_XU)=\nabla^{\mathrm{ad}}_XU+L_XU,\qquad
    \mathcal{H}(\nabla_XU)=\mathrm{A}_XU.
$$
For a vertical $(1,1)$-tensor $E$, the second identity has the following matrix representation with respect to the decomposition $TM=\mathcal{V}\oplus\mathcal{H}$:
$$
    \nabla_XE=
    \begin{pmatrix}
        \nabla_X^{\mathrm{ad}}E+[L_X,E] & -E\cdot\mathrm{A}_X\\
        \mathrm{A}_X\cdot E & 0
    \end{pmatrix}.
$$
\end{remark}

For vertical vector fields $U,V$ and a horizontal vector field $X$, Proposition~\ref{Adam} and the self-adjointness of $L_X$ give
$$
    X\langle U,V\rangle
    =\langle\nabla_X^{\mathrm{ad}}U,V\rangle
    +\langle U,\nabla_X^{\mathrm{ad}}V\rangle
    +2\langle L_XU,V\rangle.
$$
Thus, $L_X$ measures the failure of $\nabla^{\mathrm{ad}}$ to preserve the vertical metric.
This explains the correction terms that appear when differentiating metric adjoints in Subsection~\ref{subsec:metric-adjoints}.

The connection also admits the following description in terms of Lie brackets and Lie derivatives.

\begin{prop}\label{prop:adjoint-connection-bracket}
Let $U$ be a vertical vector field and $X$ a basic vector field.
Then
$$
    \nabla^{\mathrm{ad}}_XU=[X,U].
$$
Moreover, if $\mathrm{g}$ denotes the metric on $M$ and $\mathcal{L}$ the Lie derivative, then $(\mathcal{L}_U\mathrm{g})(X,\cdot)$ is the $1$-form metrically dual to $\nabla^{\mathrm{ad}}_XU$.
\end{prop}
\begin{proof}
Since $X$ is basic, $[X,U]$ is vertical by Lemma~\ref{lem:basic-vertical-bracket}.
The Levi-Civita connection is torsion-free, so taking vertical components in $\nabla_XU-\nabla_UX=[X,U]$ and applying Proposition~\ref{Adam} gives
$$
    [X,U]=\mathcal{V}(\nabla_XU)-L_XU=\nabla_X^{\mathrm{ad}}U.
$$
For the second assertion, let $V$ be a vertical vector field.
By Proposition~\ref{Adam},
$$
    \langle\nabla_XU,V\rangle
    =\langle\nabla_X^{\mathrm{ad}}U,V\rangle+\langle L_XU,V\rangle.
$$
On the other hand, orthogonality and the self-adjointness of $L_X$ give
$$
    \langle X,\nabla_VU\rangle
    =-\langle\nabla_VX,U\rangle
    =-\langle L_XV,U\rangle
    =-\langle V,L_XU\rangle.
$$
Consequently, the shape-operator terms cancel in the Lie derivative:
$$
    (\mathcal{L}_U\mathrm{g})(X,V)
    =\langle\nabla_XU,V\rangle+\langle X,\nabla_VU\rangle
    =\langle\nabla_X^{\mathrm{ad}}U,V\rangle.
$$
Finally, for a horizontal vector field $Y$, the identities for O'Neill's tensor $\mathrm{A}$ (see \cite[Section~9.C]{Besse87}) yield
\begin{align*}
    (\mathcal{L}_U\mathrm{g})(X,Y)
    &=\langle\nabla_XU,Y\rangle+\langle X,\nabla_YU\rangle
    =\langle\mathrm{A}_XU,Y\rangle+\langle X,\mathrm{A}_YU\rangle\\
    &=-\langle U,\mathrm{A}_XY\rangle-\langle\mathrm{A}_YX,U\rangle
    =-\langle U,\mathrm{A}_XY+\mathrm{A}_YX\rangle=0.
\end{align*}
This proves the asserted identity on both vertical and horizontal vectors.
\end{proof}

\begin{remark}\label{rem:induced-adjoint-connection}
If $U$ is an invariant vertical vector field, then ${}^{\mathfrak{g}}U_{g\cdot x}=\mathrm{Ad}(g)({}^{\mathfrak{g}}U_x)$.
If $X$ is basic, then $X({}^{\mathfrak{g}}U)$ satisfies the same equivariance identity.
Differentiating equivariant maps in horizontal directions is precisely the construction of the connection on an associated vector bundle induced by the principal connection $\mathcal{H}=\mathcal{V}^{\perp}$.
Thus, when $X$ is basic and $U$ is invariant, Definition~\ref{def:adjoint-connection} agrees with this induced connection on $\mathrm{ad}(M)$.
\end{remark}

\subsection{Metric adjoints}\label{subsec:metric-adjoints}

\textcolor{white}{aaa}

For a $(1,1)$-tensor $E$ on $M$, let $E^t$ denote its pointwise metric adjoint, and let $\mathcal{V}(E)$ denote the vertical $(1,1)$-tensor defined by $\mathcal{V}(E)(U):=\mathcal{V}(EU)$ for vertical $U$.
Then $\mathcal{V}(E^t)=\mathcal{V}(E)^t$ and $\nabla_X(E^t)=(\nabla_XE)^t$.
If $E$ is vertical, Remark~\ref{rem:adjoint-connection-block-form} gives
$$
    \mathcal{V}(\nabla_XE)=\nabla^{\mathrm{ad}}_XE+[L_X,E].
$$
These observations imply the following compatibility formula for $\nabla^{\mathrm{ad}}$ and the metric adjoint.

\begin{prop}\label{prop:adjoint-metric-transpose}
Let $E$ be a vertical $(1,1)$-tensor.
Then, for every horizontal vector field $X$,
$$
    \nabla^{\mathrm{ad}}_X(E^t)=(\nabla^{\mathrm{ad}}_XE)^t-2[L_X,E^t].
$$
In particular, writing $\mathrm{S}(E)=(E+E^t)/2$ for the symmetric part of $E$, we have
$$
    \nabla^{\mathrm{ad}}_X(\mathrm{S}(E))=\mathrm{S}(\nabla^{\mathrm{ad}}_XE)-[L_X,E^t].
$$
\end{prop}
\begin{proof}
By Proposition~\ref{Adam},
$$
    \mathcal{V}(\nabla_XE)=\nabla^{\mathrm{ad}}_XE+[L_X,E],\qquad
    \mathcal{V}(\nabla_X(E^t))=\nabla^{\mathrm{ad}}_X(E^t)+[L_X,E^t].
$$
Since $\mathcal{V}(\nabla_X(E^t))=\mathcal{V}(\nabla_XE)^t$ and $L_X$ is self-adjoint, taking the adjoint of the first identity and comparing it with the second yields
$$
    (\nabla^{\mathrm{ad}}_XE)^t-[L_X,E^t]=\nabla^{\mathrm{ad}}_X(E^t)+[L_X,E^t].
$$
Rearranging proves the first formula.
The second follows by applying the first to the definition of $\mathrm{S}(E)$.
\end{proof}

\begin{remark}
Alternatively, the first identity in Proposition~\ref{prop:adjoint-metric-transpose} can be deduced from Proposition~\ref{lem:inner-product-adjoint-differential} and Lemma~\ref{lem:orbit-metric-differential}.
Indeed, ${}^{\mathfrak{g}}(E^t)_x=({}^{\mathfrak{g}}E_x)^{t_{h_x}}$, so ${}^{\mathfrak{g}}\nabla^{\mathrm{ad}}_{X_x}(E^t)$ is the derivative in the direction $X_x$ of the composition of $x\mapsto({}^{\mathfrak{g}}E_x,h_x)$ with $(F,h)\mapsto F^{t_h}$.
\end{remark}

Write $\mathrm{sk}(E)=(E-E^t)/2$ for the skew-symmetric part of $E$.
Together with Proposition~\ref{Adam}, Proposition~\ref{prop:adjoint-metric-transpose} gives the following corrected version of \cite[Lemma~5.7]{BL23}.

\begin{corollary}\label{cor:self-adjoint-adjoint-derivative}
If $E$ is a pointwise self-adjoint vertical $(1,1)$-tensor and $X$ is a horizontal vector field, then
$$
    \mathrm{S}(\nabla^{\mathrm{ad}}_XE)=\mathcal{V}(\nabla_XE),\qquad
    \mathrm{sk}(\nabla^{\mathrm{ad}}_XE)=-[L_X,E].
$$
\end{corollary}

\subsection{Curvature for polar actions}\label{subsec:curvature-metric-adjoints}

\textcolor{white}{aaa}

Throughout this subsection, we assume that the action is polar, that is, $\mathrm{A}=0$.
We first record the flatness identity for $\nabla^{\mathrm{ad}}$ and then compare it with the curvature of the Levi-Civita connection.

\begin{lemma}\label{lem:polar-adjoint-flatness}
Suppose that $\mathrm{A}=0$, so that the Lie bracket of horizontal vector fields is horizontal.
Then, on both vertical vector fields and vertical $(1,1)$-tensors,
$$
    \nabla^{\mathrm{ad}}_X\nabla^{\mathrm{ad}}_Y-\nabla^{\mathrm{ad}}_Y\nabla^{\mathrm{ad}}_X
    =\nabla^{\mathrm{ad}}_{[X,Y]}
    \qquad\text{for all horizontal vector fields }X,Y.
$$
\end{lemma}
\begin{proof}
The asserted identity is tensorial in $X$ and $Y$, so we may assume that both are basic.
Since $\mathrm{A}=0$, their bracket $[X,Y]$ is also basic.
For a vertical vector field $U$, Proposition~\ref{prop:adjoint-connection-bracket} and the Jacobi identity give
$$
    \nabla^{\mathrm{ad}}_X\nabla^{\mathrm{ad}}_YU-\nabla^{\mathrm{ad}}_Y\nabla^{\mathrm{ad}}_XU
    =[X,[Y,U]]-[Y,[X,U]]=[[X,Y],U]=\nabla^{\mathrm{ad}}_{[X,Y]}U.
$$
The corresponding identity for vertical $(1,1)$-tensors follows by applying Lemma~\ref{lem:adjoint-connection-product-rule} twice.
\end{proof}

Let $\mathrm{R}$ denote the $(1,3)$-curvature tensor of the Levi-Civita connection, with convention $\mathrm{R}(X,Y)=[\nabla_X,\nabla_Y]-\nabla_{[X,Y]}$.
If $\mathrm{A}=0$, covariant differentiation in horizontal directions preserves $\mathcal{V}$ and $\mathcal{H}$.
Consequently, $\mathrm{R}(X,Y)(\mathcal{V})\subseteq\mathcal{V}$ for horizontal $X,Y$, and its restriction to $\mathcal{V}$ is a vertical $(1,1)$-tensor.

\begin{prop}\label{prop:polar-curvature-shape}
If $\mathrm{A}=0$, then, for every pair of horizontal vector fields $X,Y$, the following identities hold on vertical vectors:
$$
    \mathrm{R}(X,Y)=-[L_X,L_Y],\qquad
    (\nabla_XL)_Y=(\nabla_YL)_X.
$$
\end{prop}
\begin{proof}
We first establish the identity
\begin{equation}\label{eq:polar-curvature-symmetry}
    \mathrm{R}(X,Y)+[L_X,L_Y]=(\nabla_XL)_Y-(\nabla_YL)_X.
\end{equation}
The left-hand side is pointwise skew-adjoint on $\mathcal{V}$.
The right-hand side is pointwise self-adjoint, since each $L_Z$ is self-adjoint and covariant differentiation in horizontal directions preserves the vertical metric.
Thus, both sides must vanish, giving the claimed identities.

To prove~\eqref{eq:polar-curvature-symmetry}, Proposition~\ref{Adam} gives, as operators on vertical vector fields,
$$
    [\nabla_X,\nabla_Y]
    =[\nabla^{\mathrm{ad}}_X,\nabla^{\mathrm{ad}}_Y]
    +[L_X,\nabla^{\mathrm{ad}}_Y]+[\nabla^{\mathrm{ad}}_X,L_Y]
    +[L_X,L_Y].
$$
By Lemma~\ref{lem:polar-adjoint-flatness}, $[\nabla^{\mathrm{ad}}_X,\nabla^{\mathrm{ad}}_Y]=\nabla^{\mathrm{ad}}_{[X,Y]}$.
Applying Proposition~\ref{Adam} to $\nabla_{[X,Y]}$ and subtracting, we obtain
\begin{equation}\label{eq:polar-curvature-expansion}
    \mathrm{R}(X,Y)=[L_X,\nabla^{\mathrm{ad}}_Y]+[\nabla^{\mathrm{ad}}_X,L_Y]+[L_X,L_Y]-L_{[X,Y]}.
\end{equation}
For the right-hand side of~\eqref{eq:polar-curvature-symmetry}, the torsion-free property gives
$$
    (\nabla_XL)_Y-(\nabla_YL)_X=\nabla_X(L_Y)-\nabla_Y(L_X)-L_{[X,Y]}.
$$
Proposition~\ref{Adam} and Lemma~\ref{lem:adjoint-connection-product-rule} yield $\nabla_X(L_Y)=[\nabla^{\mathrm{ad}}_X,L_Y]+[L_X,L_Y]$.
Hence,
\begin{equation}\label{eq:shape-derivative-antisymmetrization}
    (\nabla_XL)_Y-(\nabla_YL)_X
    =[\nabla^{\mathrm{ad}}_X,L_Y]+[L_X,\nabla^{\mathrm{ad}}_Y]+2[L_X,L_Y]-L_{[X,Y]}.
\end{equation}
Combining~\eqref{eq:polar-curvature-expansion} and~\eqref{eq:shape-derivative-antisymmetrization} proves~\eqref{eq:polar-curvature-symmetry}.
\end{proof}

\subsection{Horizontally Killing vector fields}\label{subsec:horizontally-killing}

\textcolor{white}{aaa}

In this subsection, we return to our original assumptions, without assuming that the $\mathsf{G}$-action is polar.
For every vertical vector field $U$, the proof of Proposition~\ref{prop:adjoint-connection-bracket} shows that $(\mathcal{L}_U\mathrm{g})(X,Y)=0$ whenever $X,Y$ are horizontal.
We are interested in fields for which the mixed components also vanish.

\begin{definition}\label{def:horizontally-killing}
A vertical vector field $U$ on $M$ is called \emph{horizontally Killing} if it satisfies any of the following equivalent conditions:
\begin{enumerate}
\renewcommand{\labelenumi}{(\roman{enumi})}
\item $\nabla^{\mathrm{ad}}_XU=0$ for every horizontal vector field $X$.
\item ${}^{\mathfrak{g}}U:M\rightarrow\mathfrak{g}$ is constant along horizontal curves.
\item $[X,U]=0$ for every basic vector field $X$.
\item $\mathcal{V}(\nabla_XU)=L_XU$ for every horizontal vector field $X$.
\item $(\mathcal{L}_U\mathrm{g})(X,V)=0$ for every horizontal vector field $X$ and every vertical vector field $V$.
\end{enumerate}
\end{definition}

The equivalences follow from Definition~\ref{def:adjoint-connection} and Propositions~\ref{prop:adjoint-connection-bracket} and~\ref{Adam}.
Invariant horizontally Killing vector fields correspond precisely to parallel sections of the adjoint bundle.
Every vertical Killing field is horizontally Killing, by Lemma~\ref{lem:basic-vertical-bracket}, although it need not be invariant.

The vanishing of the horizontal and mixed components of $\mathcal{L}_U\mathrm{g}$ gives the following characterization of vertical Killing fields.

\begin{remark}\label{rem:vertical-killing-characterization}
A vertical vector field $U$ on $M$ is Killing if and only if it is horizontally Killing and $U|_{\mathsf{G}\cdot x}\in\mathrm{Kill}(\mathsf{G}\cdot x)$ for every $x\in M$.
\end{remark}

Any prescribed vector $U_{x_0}\in\mathcal{V}_{x_0}$ extends to a horizontally Killing vector field: for example, take $U=\theta_{x_0}^{-1}(U_{x_0})\in\mathfrak{g}$ and consider the fundamental field $U^*$.
Such an extension need not be unique among horizontally Killing vector fields.
The following proposition describes the holonomy obstruction to choosing an invariant extension, which is unique whenever it exists.

\begin{prop}\label{prop:holonomy-horizontally-killing}
Let $x_0\in M$ and $U_{x_0}\in\mathcal{V}_{x_0}$.
Then the following conditions are equivalent:
\begin{enumerate}
\renewcommand{\labelenumi}{(\roman{enumi})}
\item There exists an invariant horizontally Killing vertical vector field $U$ on $M$ whose value at $x_0$ is $U_{x_0}$.
\item $\mathrm{Ad}(g_{\mathrm{h}})({}^{\mathfrak{g}}U_{x_0})={}^{\mathfrak{g}}U_{x_0}$ for every $g_{\mathrm{h}}\in\mathsf{Hol}_{x_0}$.
\end{enumerate}
If such a vector field exists, it is unique.
\end{prop}
\begin{proof}

By the holonomy principle for vector bundles, a vector in a fiber over a connected base extends to a unique smooth parallel section if and only if it is fixed by the full holonomy group; see \cite[Theorem~2.7]{Agricola06}.
We apply this principle to the adjoint bundle $\mathrm{ad}(M)\rightarrow M/\mathsf{G}$ with its induced connection $\nabla^{\mathrm{ad}}$.

The prescribed vector $U_{x_0}$ corresponds to $[x_0,{}^{\mathfrak{g}}U_{x_0}]\in\mathrm{ad}(M)_{\pi(x_0)}$.
If a loop based at $\pi(x_0)$ has horizontal lift starting at $x_0$ and ending at $g_{\mathrm{h}}\cdot x_0$, then parallel transport in the adjoint bundle sends this element to
$$
    [g_{\mathrm{h}}\cdot x_0,{}^{\mathfrak{g}}U_{x_0}]
    =[x_0,\mathrm{Ad}(g_{\mathrm{h}}^{-1})({}^{\mathfrak{g}}U_{x_0})].
$$
Since $\mathsf{Hol}_{x_0}$ is closed under inversion, $[x_0,{}^{\mathfrak{g}}U_{x_0}]$ is fixed by the full holonomy group of $\nabla^{\mathrm{ad}}$ if and only if condition~(ii) holds.

Finally, parallel sections of the adjoint bundle correspond precisely to invariant horizontally Killing vertical vector fields on $M$.
Thus, the holonomy principle gives the equivalence of~(i) and~(ii), as well as uniqueness.
\end{proof}

\section{Commutativity relations for the shape operator}\label{sec:shape-commutativity}

As stated in the Introduction, throughout Sections~\ref{sec:shape-commutativity}--\ref{sec:solvable-splitting}, we work under Assumption \ref{ass:main}. In this section, we prove the following commutativity relations between the shape operators and the adjoint maps determined by $\mathfrak{a}_x$.

\begin{theorem}\label{thm:shape-adjoint-commutativity}
Let $M$ be a complete, connected negative-Einstein manifold and let $\mathsf{G}$ be a connected solvable Lie group acting freely, properly, isometrically and cocompactly on $M$.
Then, for every $x\in M$, the following properties hold:
\begin{enumerate}
\renewcommand{\labelenumi}{(\roman{enumi})}
    \item $[\mathrm{ad}_A^{t_x},\mathrm{ad}_{\tilde{A}}]=0$ for every $A,\tilde{A}\in\mathfrak{a}_x$;
    \item $[{}^{\mathfrak{g}}L_{X_x},\mathrm{ad}_A]=0$ for every $A\in\mathfrak{a}_x$ and $X_x\in\mathcal{H}_x$.
\end{enumerate}
Equivalently,
$$
    [{}^{\mathfrak{n}}\hat{L}_{\widehat{X}_x},\mathrm{ad}_A|_{\mathfrak{n}}]=0
    \qquad\text{for every }A\in\mathfrak{a}_x
    \text{ and }\widehat{X}_x\in\widehat{\mathcal{H}}_x.
$$
\end{theorem}

The equivalence follows from Remark~\ref{rem:horizontal-complement-equivariance}\,(ii), Lemma~\ref{lem:vertical-tensor-blocks} and Lemma~\ref{lem:shape-symmetric-adjoint}.
Item~(i) extends the normality statement in \cite[Theorem~4.10(1)]{Heber98} from Einstein solvmanifolds to our cocompact setting.

Combining \cite[Proposition~9.7]{BL23} with Lemma~\ref{lem:shape-symmetric-adjoint}, we have $\mathrm{ad}_A|_{\mathfrak{n}}^{t_x}\in\mathrm{Der}(\mathfrak{n})$ for every $x\in M$ and every $A\in\mathfrak{a}_x$.
By Remark~\ref{rem:horizontal-complement-equivariance}\,(ii), $\mathrm{ad}_A^{t_x}$ vanishes on $\mathfrak{a}_x$ and preserves $\mathfrak{n}$.
Thus, it is a derivation of $\mathfrak{g}$ if and only if
$$
    [\mathrm{ad}_{\tilde{A}}|_{\mathfrak{n}},\mathrm{ad}_A|_{\mathfrak{n}}^{t_x}]=0
    \qquad\text{for every }\tilde{A}\in\mathfrak{a}_x.
$$
Indeed, the derivation identity is already satisfied for pairs of vectors in $\mathfrak{n}$ and for pairs in $\mathfrak{a}_x$, while the displayed condition is precisely the derivation identity for mixed pairs.
Consequently, item~(i) in Theorem \ref{thm:shape-adjoint-commutativity} is equivalent to $\mathrm{ad}_A^{t_x}\in\mathrm{Der}(\mathfrak{g})$ for every $A\in\mathfrak{a}_x$.
We will prove item~(i) in this equivalent form.

The proof has three main ingredients.
First, vectors in $\mathcal{A}$ extend to invariant horizontally Killing vector fields.
Secondly, an algebraic Ricci estimate, obtained from Heber's argument, gives a pointwise nonnegative quantity.
Finally, we express this quantity, together with the squared norms of the commutators in item~(ii), as a weighted divergence on the compact quotient $M/\mathsf{G}$.
We develop these ingredients in the first three subsections and collect the consequences of the theorem in the last one.

\subsection{Invariant horizontally Killing extensions}\label{subsec:commutativity-preliminaries}

\textcolor{white}{aaa}

We begin with the extension result that allows us to apply a global divergence argument to an arbitrary vector in $\mathcal{A}_x$.

\begin{prop}\label{prop:invariant-a-extension}
Under Assumption~\ref{ass:main}, the holonomy group fixes $\mathfrak{a}_x$ pointwise:
$$
    \mathrm{Ad}(g_{\mathrm{h}})A=A
    \qquad\text{for every }A\in\mathfrak{a}_x
    \text{ and }g_{\mathrm{h}}\in\mathsf{Hol}_x.
$$
Equivalently, for every $x_0\in M$ and $A_{x_0}\in\mathcal{A}_{x_0}$, there exists a unique invariant horizontally Killing vertical vector field $A$ on $M$ with the prescribed value at $x_0$.
This extension is $\mathcal{A}$-valued.
\end{prop}

\begin{proof}
Fix $x_0\in M$.
Let $r:\widetilde{\mathsf{G}}\rightarrow\mathsf{G}$ be the universal covering homomorphism, identifying the Lie algebras through $dr_e$.
Let $\widetilde{\mathsf{N}}$ and $\widetilde{\mathsf{A}}_0$ be the connected subgroups of $\widetilde{\mathsf{G}}$ with Lie algebras $\mathfrak{n}$ and $\mathfrak{a}_{x_0}$, respectively.
By Corollary~\ref{cor:horizontal-complement-subalgebra} and Remark~\ref{rem:horizontal-complement-equivariance}\,(ii), $\mathfrak{g}=\mathfrak{a}_{x_0}\ltimes\mathfrak{n}$, with $\mathfrak{a}_{x_0}$ abelian.
Thus, \cite[Lemma~3.18.4]{Varadarajan84} implies that these subgroups are closed and simply connected, and that
$$
    \widetilde{\mathsf{G}}=\widetilde{\mathsf{A}}_0\ltimes\widetilde{\mathsf{N}}.
$$
The exponential maps of the simply connected abelian group $\widetilde{\mathsf{A}}_0$ and the simply connected nilpotent group $\widetilde{\mathsf{N}}$ are diffeomorphisms; see \cite[Theorem~1.127]{Knapp02}.
Consequently,
$$
    \mathfrak{a}_{x_0}\times\mathfrak{n}\longrightarrow\widetilde{\mathsf{G}},
    \qquad
    (A,\widehat{U})\longmapsto\exp(A)\exp(\widehat{U})
$$
is a diffeomorphism.

Let $g_{\mathrm{h}}\in\mathsf{Hol}_{x_0}$.
Lift it to $\tilde{g}_{\mathrm{h}}\in\widetilde{\mathsf{G}}$ and write $\tilde{g}_{\mathrm{h}}=\exp(A)\exp(\widehat{U})$.
Since $\mathfrak{a}_{x_0}$ is abelian and $\mathfrak{n}$ is an ideal, the adjoint maps have the following block forms relative to $\mathfrak{g}=\mathfrak{a}_{x_0}\oplus\mathfrak{n}$:
$$
    \mathrm{ad}_A=
    \begin{pmatrix}0&0\\0&\mathrm{ad}_A|_{\mathfrak{n}}\end{pmatrix},
    \qquad
    \mathrm{ad}_{\widehat{U}}=
    \begin{pmatrix}0&0\\\mathrm{ad}_{\widehat{U}}|_{\mathfrak{a}_{x_0}}&\mathrm{ad}_{\widehat{U}}|_{\mathfrak{n}}\end{pmatrix}.
$$
Exponentiating these matrices, we obtain
$$
    e^{\mathrm{ad}_A}=
    \begin{pmatrix}\mathrm{Id}_{\mathfrak{a}_{x_0}}&0\\0&e^{\mathrm{ad}_A|_{\mathfrak{n}}}\end{pmatrix},
    \qquad
    e^{\mathrm{ad}_{\widehat{U}}}=
    \begin{pmatrix}\mathrm{Id}_{\mathfrak{a}_{x_0}}&0\\ *&e^{\mathrm{ad}_{\widehat{U}}|_{\mathfrak{n}}}\end{pmatrix}.
$$
Since $\mathrm{Ad}(r(\tilde{g}))=\mathrm{Ad}(\tilde{g})$, multiplying these matrices gives
\begin{equation}\label{eq:holonomy-adjoint-blocks}
    \mathrm{Ad}(g_{\mathrm{h}})
    =e^{\mathrm{ad}_A}e^{\mathrm{ad}_{\widehat{U}}}
    =\begin{pmatrix}
        \mathrm{Id}_{\mathfrak{a}_{x_0}}&0\\
        *&e^{\mathrm{ad}_A|_{\mathfrak{n}}}e^{\mathrm{ad}_{\widehat{U}}|_{\mathfrak{n}}}
    \end{pmatrix}.
\end{equation}
On the other hand, $x\mapsto\mathfrak{a}_x$ is equivariant by Remark~\ref{rem:horizontal-complement-equivariance}\,(i) and constant along horizontal curves by Corollary~\ref{cor:horizontal-complement-subalgebra}.
Since $x_0$ and $g_{\mathrm{h}}\cdot x_0$ can be joined by such a curve,
\begin{equation}\label{eq:holonomy-preserves-complement}
    \mathrm{Ad}(g_{\mathrm{h}})(\mathfrak{a}_{x_0})
    =\mathfrak{a}_{g_{\mathrm{h}}\cdot x_0}
    =\mathfrak{a}_{x_0}.
\end{equation}
Combining~\eqref{eq:holonomy-adjoint-blocks} and~\eqref{eq:holonomy-preserves-complement}, the lower-left block must vanish.
Thus $\mathrm{Ad}(g_{\mathrm{h}})|_{\mathfrak{a}_{x_0}}=\mathrm{Id}_{\mathfrak{a}_{x_0}}$.

Proposition~\ref{prop:holonomy-horizontally-killing} now gives the unique invariant horizontally Killing extension.
Corollary~\ref{cor:horizontal-complement-subalgebra} shows that $\mathcal{A}$ is preserved by horizontal parallel transport for $\nabla^{\mathrm{ad}}$, and equation~\eqref{eq:a-invariance} shows that it is preserved by the $\mathsf{G}$-action.
Since the extension is parallel in horizontal directions and invariant, its prescribed value in $\mathcal{A}_{x_0}$ remains $\mathcal{A}$-valued throughout $M$.
\end{proof}

\begin{remark}\label{rem:a-extension-quotient-lift}
Alternatively, we can construct the invariant horizontally Killing extension by taking the $\mathsf{N}$-horizontal lift of a fundamental vector field on $M/\mathsf{N}$.
 Identify ${}^{\mathfrak{g}}A_{x_0}\in\mathfrak{a}_{x_0}$ with its class in $\mathfrak{g}/\mathfrak{n}$ and let $\bar{A}$ be the corresponding fundamental vector field on $M/\mathsf{N}$.
 Since $\mathsf{G}/\mathsf{N}$ is abelian, $\bar{A}$ is $\mathsf{G}/\mathsf{N}$-invariant.
 Its $\mathsf{N}$-horizontal lift to $M$ is $\mathsf{G}$-invariant and has the prescribed value at $x_0$.
 It is $\mathcal{A}$-valued because the differential of $M\rightarrow M/\mathsf{N}$ maps $\mathcal{A}$ isomorphically onto the tangent spaces of the $\mathsf{G}/\mathsf{N}$-orbits. 
 Using Corollary~\ref{cor:horizontal-complement-subalgebra}, one can show that this vector field is horizontally Killing.
\end{remark}

\begin{remark}\label{rem:a-fields-n-horizontal-constancy}
If $A$ is invariant, horizontally Killing and $\mathcal{A}$-valued, then ${}^{\mathfrak{g}}A$ is constant along $\mathsf{N}$-horizontal curves.
It is already constant along $\mathsf{G}$-horizontal curves by definition.
To treat the remaining directions, let $\tilde{A}$ be $\mathcal{A}$-valued and fix $x\in M$.
Along $\sigma(t):=\exp(t\,{}^{\mathfrak{g}}\tilde{A}_x)\cdot x$, invariance gives
$$
    {}^{\mathfrak{g}}A_{\sigma(t)}
    =\mathrm{Ad}\bigl(\exp(t\,{}^{\mathfrak{g}}\tilde{A}_x)\bigr)({}^{\mathfrak{g}}A_x).
$$
Differentiating at $t=0$, we obtain
$$
    \tilde{A}_x({}^{\mathfrak{g}}A)
    =[{}^{\mathfrak{g}}\tilde{A}_x,{}^{\mathfrak{g}}A_x]=0,
$$
since $\mathfrak{a}_x$ is abelian by Remark~\ref{rem:horizontal-complement-equivariance}\,(ii).
The conclusion follows from $\widehat{\mathcal{H}}=\mathcal{H}\oplus\mathcal{A}$.
\end{remark}

We also record a decomposition of the mean curvature vector field.
It will provide the weight in the divergence argument.

\begin{prop}\label{prop:mean-curvature-decomposition}
There exists a $\mathsf{G}$-invariant smooth function $\log\check{v}:M\rightarrow\mathbb{R}$ such that
$$
    H=\nabla\log\check{v}+\mathcal{H}(\widehat{H}).
$$
\end{prop}

\begin{proof}
Let $\{A_i\}$ be a local invariant orthonormal frame for $\mathcal{A}$.
Completing it with a local orthonormal frame for $\widehat{\mathcal{V}}$, the definitions of $H$ and $\widehat{H}$ give
$$
    H=\sum_i\mathcal{H}(\nabla_{A_i}A_i)+\mathcal{H}(\widehat{H}).
$$
The connected abelian group $\mathsf{G}/\mathsf{N}$ acts freely, properly and isometrically on $M/\mathsf{N}$.
Since it is abelian, it is unimodular, so \cite[Lemma~3.2]{BL23} implies that the mean curvature vector field of its orbits is the gradient of a $\mathsf{G}/\mathsf{N}$-invariant function $\log\check{v}$ on $M/\mathsf{N}$.
Under the Riemannian submersion $p:M\rightarrow M/\mathsf{N}$, the fields $A_i$ project to an orthonormal frame tangent to these orbits.
The horizontal lift of their mean curvature vector field is therefore
$$
    \sum_i\mathcal{H}(\nabla_{A_i}A_i)=\nabla(\log\check{v}\circ p).
$$
Using the same notation for the function composed with $p$ proves the formula.
\end{proof}

\subsection{An algebraic Ricci estimate}\label{subsec:commutativity-algebraic-estimate}

\textcolor{white}{aaa}

For an inner product $h$ on a Lie algebra $\mathfrak{g}$, denote by $\mathrm{Ric}_h$ the Ricci endomorphism of the corresponding left-invariant metric.
The following estimate is obtained from the argument preceding \cite[Lemma~4.1]{Heber98}.
We retain the Ricci term, which vanishes in Heber's Einstein setting.

\begin{prop}\label{prop:heber-derivation-estimate}
Let $\mathfrak{g}$ be a solvable Lie algebra with Killing form $B$ and an inner product $h$.
Define $\mathring{H}_h\in\mathfrak{g}$ by
$$
    h(\mathring{H}_h,U)=\mathrm{Tr}(\mathrm{ad}_U)
    \qquad\text{for every }U\in\mathfrak{g}.
$$
For every $D\in\mathrm{Der}(\mathfrak{g})$ and every $h$-orthonormal basis $\{e_i\}$,
$$
    -\mathrm{Tr}(\mathrm{Ric}_h\cdot[D,D^t])
    +\mathrm{Tr}(\mathrm{ad}_{D(\mathring{H}_h)}\cdot D^t)
    +\frac12\sum_iB(D^tDe_i,e_i)\leq0.
$$
Equality holds if and only if $D^t\in\mathrm{Der}(\mathfrak{g})$.
\end{prop}

\begin{corollary}\label{cor:ricci-commutator-estimate}
Let $\mathfrak{g}$ be a solvable Lie algebra with nilradical $\mathfrak{n}$ and an inner product $h$ such that $\mathfrak{n}^{\perp}$ is an abelian subalgebra.
Then, for every $A\in\mathfrak{n}^{\perp}$,
$$
    \mathrm{Tr}(\mathrm{Ric}_h\cdot[\mathrm{ad}_A,\mathrm{ad}_A^t])\geq0,
$$
with equality if and only if $\mathrm{ad}_A^t\in\mathrm{Der}(\mathfrak{g})$.
\end{corollary}

\begin{proof}
Apply Proposition~\ref{prop:heber-derivation-estimate} to $D=\mathrm{ad}_A$.
For $U\in\mathfrak{n}$, the endomorphism $\mathrm{ad}_U$ is nilpotent by \cite[Proposition~1.40]{Knapp02}, so $h(\mathring{H}_h,U)=\mathrm{Tr}(\mathrm{ad}_U)=0$.
Thus $\mathring{H}_h\in\mathfrak{n}^{\perp}$, and the abelianness of this subalgebra gives $D(\mathring{H}_h)=0$.
Choose an $h$-orthonormal basis $\{e_i\}$ adapted to $\mathfrak{g}=\mathfrak{n}^{\perp}\oplus\mathfrak{n}$.
If $e_i\in\mathfrak{n}^{\perp}$, then $De_i=0$.
If $e_i\in\mathfrak{n}$, then $B(\cdot,e_i)=0$, since the nilradical of a solvable Lie algebra is contained in the kernel of its Killing form; see \cite[Corollary~1.29 and Proposition~1.40]{Knapp02}.
Thus the last two terms vanish, giving the claimed inequality and its equality condition.
\end{proof}

We now return to Assumption~\ref{ass:main}.
For each $x\in M$, Remark~\ref{rem:horizontal-complement-equivariance}\,(ii) ensures that Corollary~\ref{cor:ricci-commutator-estimate} applies to the metric Lie algebra $(\mathfrak{g},h_x)$.
Recall from Subsection~\ref{subsec:trivialization} that the orbit map $\mathsf{G}\rightarrow\mathsf{G}\cdot x$ is an isometry when $\mathsf{G}$ carries the left-invariant metric induced by $h_x$.
Its differential at the identity is $\theta_x$, so ${}^{\mathfrak{g}}\mathrm{Ric}^{\mathsf{G}}_x=\mathrm{Ric}_{h_x}$.

For an $\mathcal{A}$-valued vector field $A$, we define $\mathrm{ad}(A)$ as the vertical $(1,1)$-tensor corresponding, via Killing coordinates, to the smooth map $M\rightarrow\mathrm{End}(\mathfrak{g})$ given by
$$
    {}^{\mathfrak{g}}\mathrm{ad}(A)_x:=\mathrm{ad}_{{}^{\raisebox{0.15em}{$\scriptstyle\mathfrak{g}$}}\! A_x}.
$$

We use the inner product $\ip{E}{F}=\mathrm{Tr}(E\cdot F^t)$ on vertical endomorphisms, induced by the metric on $\mathcal{V}$.
In particular, $[\mathrm{ad}(A),\mathrm{ad}(A)^t]$ is self-adjoint and has zero trace.

\begin{lemma}\label{lem:tau-shape-commutator-estimate}
Let $A$ be an $\mathcal{A}$-valued vector field.
Then
$$
    \ip{\tau-L_H}{[\mathrm{ad}(A),\mathrm{ad}(A)^t]}\geq0.
$$
Equality holds at $x\in M$ if and only if $\mathrm{ad}_{{}^{\raisebox{0.15em}{$\scriptstyle\mathfrak{g}$}}\! A_x}^{\,t_x}\in\mathrm{Der}(\mathfrak{g})$.
\end{lemma}

\begin{proof}
By Corollary~\ref{cor:ricci-commutator-estimate} and the preceding identification of the Ricci endomorphisms,
$$
    \ip{\mathrm{Ric}^{\mathsf{G}}}{[\mathrm{ad}(A),\mathrm{ad}(A)^t]}\geq0,
$$
with the stated equality condition.
Since $\mathrm{ric}_{\mathrm{g}}=-\mathrm{g}$, O'Neill's formula in Proposition~\ref{prop:oneill-vertical-ricci} gives
$$
    \mathrm{Ric}^{\mathsf{G}}
    =-\mathrm{Id}_{\mathcal{V}}+\tau-L_H-\mathrm{A}^*\mathrm{A}.
$$
The identity term has zero pairing with $[\mathrm{ad}(A),\mathrm{ad}(A)^t]$ because the commutator has zero trace.
By Remark~\ref{rem:horizontal-complement-equivariance}\,(ii), $\mathrm{ad}_{{}^{\raisebox{0.15em}{$\scriptstyle\mathfrak{g}$}}\! A_x}$ and its transpose vanish on $\mathfrak{a}_x$ and preserve $\mathfrak{n}$.
Therefore $\mathrm{ad}(A)$, $\mathrm{ad}(A)^t$ and their commutator vanish on $\mathcal{A}$ and preserve $\widehat{\mathcal{V}}$.
By Lemma~\ref{lem:vertical-tensor-blocks}, $\mathrm{A}^*\mathrm{A}$ vanishes on $\widehat{\mathcal{V}}$, so its pairing with the commutator also vanishes.
This proves the lemma.
\end{proof}

\subsection{The divergence argument}\label{subsec:proof-shape-commutativity}

\textcolor{white}{aaa}

For an $\mathcal{A}$-valued vector field $A$, define the horizontal vector field $\mathrm{W}_A$ by
$$
    \ip{\mathrm{W}_A}{X}=\ip{L_X}{[\mathrm{ad}(A),\mathrm{ad}(A)^t]}
    \qquad\text{for every horizontal vector field }X.
$$

This definition is inspired in part by the formula for $\nabla\log v_{\beta}$ in \cite[Lemma~5.9]{BL23}.
Its purpose is to fit the estimate in Lemma~\ref{lem:tau-shape-commutator-estimate} into the divergence theorem of Proposition~\ref{prop:cocompact-maximum-principle}\,(i).
We expect $\mathrm{div}_B\mathrm{W}_A-\ip{\mathrm{W}_A}{H}$ to be the sum of the term in that estimate and an additional term, which fortunately turns out to be nonnegative.
The next two lemmas establish this.

If $A$ is invariant, the equivariance of the tensors in this formula shows that $\mathrm{W}_A$ is invariant and hence basic.
As before, $\mathcal{V}(\nabla_XL)_Y$ denotes the restriction to vertical vectors followed by vertical projection.

\begin{lemma}\label{lem:commutator-field-derivative}
Let $A$ be an $\mathcal{A}$-valued horizontally Killing vector field.
Then, for all horizontal vector fields $X,Y$,
$$
    \ip{\nabla_X\mathrm{W}_A}{Y}
    =\ip{\mathcal{V}(\nabla_XL)_Y}{[\mathrm{ad}(A),\mathrm{ad}(A)^t]}
    +2\ip{[L_X,\mathrm{ad}(A)]}{[L_Y,\mathrm{ad}(A)]}.
$$
\end{lemma}

\begin{proof}
Since $A$ is horizontally Killing, $\nabla_X^{\mathrm{ad}}\mathrm{ad}(A)=0$.
The product rule and Proposition~\ref{prop:adjoint-metric-transpose} give
\begin{equation}\label{eq:commutator-adjoint-derivative}
    \nabla_X^{\mathrm{ad}}[\mathrm{ad}(A),\mathrm{ad}(A)^t]
    =[\mathrm{ad}(A),\nabla_X^{\mathrm{ad}}(\mathrm{ad}(A)^t)]
    =-2[\mathrm{ad}(A),[L_X,\mathrm{ad}(A)^t]].
\end{equation}
Differentiating $\ip{\mathrm{W}_A}{Y}=\mathrm{Tr}(L_Y\cdot[\mathrm{ad}(A),\mathrm{ad}(A)^t])$, we obtain
\begin{multline}\label{eq:commutator-field-product-rule}
    \ip{\nabla_X\mathrm{W}_A}{Y}
    +\mathrm{Tr}(L_{\nabla_XY}\cdot[\mathrm{ad}(A),\mathrm{ad}(A)^t])
    \\
    {}=\mathrm{Tr}\bigl((\nabla_X^{\mathrm{ad}}L_Y)\cdot[\mathrm{ad}(A),\mathrm{ad}(A)^t]\bigr)
    +\mathrm{Tr}\bigl(L_Y\cdot\nabla_X^{\mathrm{ad}}[\mathrm{ad}(A),\mathrm{ad}(A)^t]\bigr).
\end{multline}
By Proposition~\ref{Adam},
\begin{equation}\label{eq:shape-adjoint-derivative-comparison}
    \nabla_X^{\mathrm{ad}}L_Y
    =\mathcal{V}(\nabla_XL)_Y+L_{\nabla_XY}-[L_X,L_Y].
\end{equation}
The commutator $[L_X,L_Y]$ is skew-adjoint, while $[\mathrm{ad}(A),\mathrm{ad}(A)^t]$ is self-adjoint, so their trace pairing vanishes.
Substituting~\eqref{eq:shape-adjoint-derivative-comparison} into~\eqref{eq:commutator-field-product-rule} and cancelling the terms containing $L_{\nabla_XY}$ yields
\begin{equation}\label{eq:commutator-field-derivative-reduction}
    \ip{\nabla_X\mathrm{W}_A}{Y}
    =\ip{\mathcal{V}(\nabla_XL)_Y}{[\mathrm{ad}(A),\mathrm{ad}(A)^t]}
    +\mathrm{Tr}\bigl(L_Y\cdot\nabla_X^{\mathrm{ad}}[\mathrm{ad}(A),\mathrm{ad}(A)^t]\bigr).
\end{equation}
For the last term, using~\eqref{eq:commutator-adjoint-derivative} and cyclicity of the trace gives
$$
    \mathrm{Tr}\bigl(L_Y\cdot\nabla_X^{\mathrm{ad}}[\mathrm{ad}(A),\mathrm{ad}(A)^t]\bigr)
    =-2\mathrm{Tr}\bigl([L_Y,\mathrm{ad}(A)]\cdot[L_X,\mathrm{ad}(A)^t]\bigr).
$$
Since $[L_X,\mathrm{ad}(A)^t]=-[L_X,\mathrm{ad}(A)]^t$, this is $2\ip{[L_X,\mathrm{ad}(A)]}{[L_Y,\mathrm{ad}(A)]}$.
Substituting into~\eqref{eq:commutator-field-derivative-reduction} proves the identity.
\end{proof}

\begin{lemma}\label{lem:commutator-divergence-inequality}
Let $A$ be an invariant horizontally Killing $\mathcal{A}$-valued vector field.
Then
$$
    \mathrm{div}_B\mathrm{W}_A-\ip{\mathrm{W}_A}{H}\geq0.
$$
Equality holds at $x\in M$ if and only if $\mathrm{ad}_{{}^{\raisebox{0.15em}{$\scriptstyle\mathfrak{g}$}}\! A_x}^{\,t_x}\in\mathrm{Der}(\mathfrak{g})$ and $[L_{X_x},\mathrm{ad}(A)_x]=0$ for every $X_x\in\mathcal{H}_x$.
\end{lemma}

\begin{proof}
Taking the horizontal trace in Lemma~\ref{lem:commutator-field-derivative}, for a local orthonormal horizontal frame $\{X_i\}$ we obtain
$$
    \mathrm{div}_B\mathrm{W}_A-\ip{\mathrm{W}_A}{H}
    =\ip{\tau-L_H}{[\mathrm{ad}(A),\mathrm{ad}(A)^t]}
    +2\sum_i\|[L_{X_i},\mathrm{ad}(A)]\|^2.
$$
The conclusion follows from Lemma~\ref{lem:tau-shape-commutator-estimate}.
\end{proof}

The decomposition of $H$ in Proposition~\ref{prop:mean-curvature-decomposition} almost puts this inequality in the form required by Proposition~\ref{prop:cocompact-maximum-principle}\,(i).
The remaining term vanishes because of the relation between the $\mathsf{N}$-orbit geometry and $\beta^+$.

\begin{proof}[Proof of Theorem~\ref{thm:shape-adjoint-commutativity}]
Fix $x\in M$ and $A_0\in\mathfrak{a}_x$.
By Proposition~\ref{prop:invariant-a-extension}, extend $\theta_x(A_0)$ to an invariant horizontally Killing $\mathcal{A}$-valued vector field $A$.
We first show that
$$
    \ip{\mathrm{W}_A}{\mathcal{H}(\widehat{H})}=0.
$$
Since $[\mathrm{ad}(A),\mathrm{ad}(A)^t]$ vanishes on $\mathcal{A}$ and preserves $\widehat{\mathcal{V}}$, Lemma~\ref{lem:vertical-tensor-blocks} gives
$$
    \ip{\mathrm{W}_A}{\mathcal{H}(\widehat{H})}
    =\mathrm{Tr}\bigl(\hat{L}_{\mathcal{H}(\widehat{H})}\cdot[\mathrm{ad}(A),\mathrm{ad}(A)^t]\bigr).
$$
By Theorem~\ref{thm:beta-properties}\,(ii),
$$
    \hat{L}_{\mathcal{H}(\widehat{H})}
    =-\beta^+-\hat{L}_{\mathcal{V}(\widehat{H})}.
$$
Theorem~\ref{thm:beta-properties}\,(v) gives $[\beta^+,\mathrm{ad}(A)]=0$, and cyclicity of the trace implies
$$
    \mathrm{Tr}\bigl(\beta^+\cdot[\mathrm{ad}(A),\mathrm{ad}(A)^t]\bigr)
    =\mathrm{Tr}\bigl([\beta^+,\mathrm{ad}(A)]\cdot \mathrm{ad}(A)^t\bigr)=0.
$$
For the other summand, we work in Killing coordinates to make the commutativity of the adjoint maps explicit.
Since $\mathcal{V}(\widehat{H})$ is $\mathcal{A}$-valued, Lemma~\ref{lem:shape-symmetric-adjoint} gives
$$
    {}^{\mathfrak{n}}\hat{L}_{\mathcal{V}(\widehat{H})}
    =-\mathrm{S}\bigl(\mathrm{ad}_{{}^{\raisebox{0.15em}{$\scriptstyle\mathfrak{g}$}}\! \mathcal{V}(\widehat{H})}|_{\mathfrak{n}}\bigr).
$$
Since $[\mathrm{ad}_{{}^{\raisebox{0.15em}{$\scriptstyle\mathfrak{g}$}}\! A},\mathrm{ad}_{{}^{\raisebox{0.15em}{$\scriptstyle\mathfrak{g}$}}\! A}^{t}]$ is self-adjoint and vanishes on $\mathfrak{a}_*$, it suffices to observe that
$$
    \mathrm{Tr}\bigl(\mathrm{ad}_{{}^{\raisebox{0.15em}{$\scriptstyle\mathfrak{g}$}}\! \mathcal{V}(\widehat{H})}\cdot[\mathrm{ad}_{{}^{\raisebox{0.15em}{$\scriptstyle\mathfrak{g}$}}\! A},\mathrm{ad}_{{}^{\raisebox{0.15em}{$\scriptstyle\mathfrak{g}$}}\! A}^{t}]\bigr)
    =\mathrm{Tr}\bigl([\mathrm{ad}_{{}^{\raisebox{0.15em}{$\scriptstyle\mathfrak{g}$}}\! \mathcal{V}(\widehat{H})},\mathrm{ad}_{{}^{\raisebox{0.15em}{$\scriptstyle\mathfrak{g}$}}\! A}]\cdot\mathrm{ad}_{{}^{\raisebox{0.15em}{$\scriptstyle\mathfrak{g}$}}\! A}^{t}\bigr)
    =0,
$$
by cyclicity of the trace and the abelianness of $\mathfrak{a}_*$.
This proves the desired orthogonality.

Combining Proposition~\ref{prop:mean-curvature-decomposition} with Lemma~\ref{lem:commutator-divergence-inequality}, we now have
$$
    \mathrm{div}_B\mathrm{W}_A-\ip{\mathrm{W}_A}{\nabla\log\check{v}}
    =\mathrm{div}_B\mathrm{W}_A-\ip{\mathrm{W}_A}{H}\geq0.
$$
Since $\mathrm{W}_A$ is basic and $\log\check{v}$ is invariant, Proposition~\ref{prop:cocompact-maximum-principle}\,(i) implies equality everywhere.
The equality condition in Lemma~\ref{lem:commutator-divergence-inequality} gives
$$
    \mathrm{ad}_{A_0}^{t_x}\in\mathrm{Der}(\mathfrak{g}),
    \qquad
    [{}^{\mathfrak{g}}L_{X_x},\mathrm{ad}_{A_0}]=0
    \quad\text{for every }X_x\in\mathcal{H}_x.
$$
The second identity proves item~(ii).
Item~(i) follows from the equivalence discussed at the beginning of this section.
Since $x$ and $A_0$ were arbitrary, the theorem follows.
\end{proof}

\subsection{Algebraic and geometric consequences}\label{subsec:commutativity-consequences}

\textcolor{white}{aaa}

We first record the consequences for the adjoint maps.
These are the cocompact counterparts of the algebraic properties used in Heber's study of Einstein solvmanifolds.

\begin{corollary}\label{cor:normal-adjoint-derivations}
For every $x\in M$, the following properties hold:
\begin{enumerate}
\renewcommand{\labelenumi}{(\roman{enumi})}
    \item For every $A\in\mathfrak{a}_x$, $\mathrm{ad}_A$ is $h_x$-normal.
    Equivalently, $\mathrm{ad}_A|_{\mathfrak{n}}$ is $\hat{h}_x$-normal.
    \item For every $A\in\mathfrak{a}_x$, $\mathrm{ad}_A^{t_x}\in\mathrm{Der}(\mathfrak{g})$.
    \item The family
$$
    \{\mathrm{S}_x(\mathrm{ad}_A),\mathrm{sk}_x(\mathrm{ad}_A):A\in\mathfrak{a}_x\}
$$
    consists of commuting derivations of $\mathfrak{g}$ that vanish on $\mathfrak{a}_x$.
\end{enumerate}
\end{corollary}

\begin{proof}
Item~(i) follows from Theorem~\ref{thm:shape-adjoint-commutativity}\,(i) by taking $\tilde{A}=A$.
Item~(ii) was obtained in the proof of the theorem.
For item~(iii), the symmetric and skew-symmetric parts are linear combinations of $\mathrm{ad}_A$ and $\mathrm{ad}_A^{t_x}$, so they are derivations vanishing on $\mathfrak{a}_x$.
The maps $\mathrm{ad}_A$ commute because $\mathfrak{a}_x$ is abelian by Remark~\ref{rem:horizontal-complement-equivariance}\,(ii); their transposes commute as well.
The mixed commutators vanish by Theorem~\ref{thm:shape-adjoint-commutativity}\,(i), proving the claim.
\end{proof}

The commutativity relations also give parallelism of the shape operators associated with invariant horizontally Killing $\mathcal{A}$-valued fields.

\begin{corollary}\label{cor:parallel-a-shape}
Let $A$ be an $\mathcal{A}$-valued vector field.
Then, for every $\mathsf{N}$-horizontal vector field $\widehat{X}$,
$$
    [\hat{L}_{\widehat{X}},\hat{L}_A]=0.
$$
Furthermore, if $A$ is invariant and horizontally Killing, then
$$
    \widehat{\nabla}^{\mathrm{ad}}_{\widehat{X}}(\hat{L}_A)=0,
    \qquad
    \nabla_{\widehat{X}}(\hat{L}_A)=0.
$$
Here $\widehat{\nabla}^{\mathrm{ad}}$ is the adjoint connection associated with the $\mathsf{N}$-action.
\end{corollary}

\begin{proof}
By Lemma~\ref{lem:shape-symmetric-adjoint}, ${}^{\mathfrak{n}}\hat{L}_A=-\mathrm{S}(\mathrm{ad}_{{}^{\raisebox{0.15em}{$\scriptstyle\mathfrak{g}$}}\! A}|_{\mathfrak{n}})$.
The equivalent formulation of Theorem~\ref{thm:shape-adjoint-commutativity} says that ${}^{\mathfrak{n}}\hat{L}_{\widehat{X}}$ commutes with $\mathrm{ad}_{{}^{\raisebox{0.15em}{$\scriptstyle\mathfrak{g}$}}\! A}|_{\mathfrak{n}}$.
Since ${}^{\mathfrak{n}}\hat{L}_{\widehat{X}}$ is self-adjoint, it also commutes with its transpose, proving the first assertion.

Suppose now that $A$ is invariant and horizontally Killing.
Remark~\ref{rem:a-fields-n-horizontal-constancy} gives $\widehat{X}(\mathrm{ad}_{{}^{\raisebox{0.15em}{$\scriptstyle\mathfrak{g}$}}\! A}|_{\mathfrak{n}})=0$.
Applying Proposition~\ref{prop:adjoint-metric-transpose} to the $\mathsf{N}$-action gives
$$
    \widehat{X}\bigl(\mathrm{S}(\mathrm{ad}_{{}^{\raisebox{0.15em}{$\scriptstyle\mathfrak{g}$}}\! A}|_{\mathfrak{n}})\bigr)
    =-[{}^{\mathfrak{n}}\hat{L}_{\widehat{X}},(\mathrm{ad}_{{}^{\raisebox{0.15em}{$\scriptstyle\mathfrak{g}$}}\! A}|_{\mathfrak{n}})^t]
    =0.
$$
Thus ${}^{\mathfrak{n}}\hat{L}_A$ is constant along $\mathsf{N}$-horizontal curves.
By the definition of the adjoint connection for the $\mathsf{N}$-action, this is precisely $\widehat{\nabla}^{\mathrm{ad}}_{\widehat{X}}(\hat{L}_A)=0$.
Since $\widehat{\mathrm{A}}=0$, Proposition~\ref{Adam} applied to that action yields
$$
    \nabla_{\widehat{X}}(\hat{L}_A)
    =\widehat{\nabla}^{\mathrm{ad}}_{\widehat{X}}(\hat{L}_A)
    +[\hat{L}_{\widehat{X}},\hat{L}_A]=0.
$$
\end{proof}

We next interpret the theorem in terms of the map $\hat{h}:M\rightarrow\mathrm{Sym}^2_+(\mathfrak{n}^*)$.
By Theorem~\ref{thm:beta-properties}\,(iii), the shape operators ${}^{\mathfrak{n}}\hat{L}_{\widehat{X}}$ are derivations of $\mathfrak{n}$.
Together with Lemma~\ref{lem:orbit-metric-differential} for the $\mathsf{N}$-action, this implies that the image of $\hat{h}$ along a connected $\mathsf{N}$-horizontal leaf lies in a single $\mathsf{Aut}(\mathfrak{n})$-orbit.
Theorem~\ref{thm:shape-adjoint-commutativity} restricts these derivations to a fixed centralizer and gives the following refinement.

\begin{corollary}\label{cor:centralizer-orbit}
Let $\widehat{\Sigma}$ be a connected leaf of the $\mathsf{N}$-horizontal distribution, and set $\mathfrak{a}_0:=\mathfrak{a}_x$ for $x\in\widehat{\Sigma}$.
This subalgebra is independent of $x$ by Corollary~\ref{cor:horizontal-complement-subalgebra}.
Define
$$
    \mathsf{J}:=\{\varphi\in\mathsf{Aut}(\mathfrak{n}):
    [\varphi,\mathrm{ad}_A|_{\mathfrak{n}}]=0
    \text{ for every }A\in\mathfrak{a}_0\}.
$$
Then $\mathsf{J}$ is a closed subgroup of $\mathsf{Aut}(\mathfrak{n})$, and the image of $\hat{h}|_{\widehat{\Sigma}}$ lies in a single $\mathsf{J}$-orbit.
\end{corollary}

\begin{proof}
The automorphism group $\mathsf{Aut}(\mathfrak{n})$ is closed in $\mathsf{GL}(\mathfrak{n})$, and the centralizer of any endomorphism is closed.
Thus $\mathsf{J}$ is a closed subgroup of $\mathsf{GL}(\mathfrak{n})$ and hence is a closed subgroup of $\mathsf{Aut}(\mathfrak{n})$.
Its Lie algebra is
$$
    \mathfrak{j}=\{D\in\mathrm{Der}(\mathfrak{n}):
    [D,\mathrm{ad}_A|_{\mathfrak{n}}]=0
    \text{ for every }A\in\mathfrak{a}_0\}.
$$
By Theorem~\ref{thm:beta-properties}\,(iii) and Theorem~\ref{thm:shape-adjoint-commutativity}, ${}^{\mathfrak{n}}\hat{L}_{\widehat{X}}\in\mathfrak{j}$ for every vector field $\widehat{X}$ tangent to $\widehat{\Sigma}$.
Thus, Lemma~\ref{lem:orbit-metric-differential} gives
$$
    d\hat{h}_x(\widehat{X}_x)
    =-\rho({}^{\mathfrak{n}}\hat{L}_{\widehat{X}_x})\hat{h}_x
    \in T_{\hat{h}_x}(\mathsf{J}\cdot\hat{h}_x).
$$
As recalled in Appendix~\ref{app:inner-products}, the natural $\mathsf{GL}(\mathfrak{n})$-action on $\mathrm{Sym}^2_+(\mathfrak{n}^*)$ is proper.
Its restriction to the closed subgroup $\mathsf{J}$ is therefore proper as well. 
Suppose, for a contradiction, that there exist $x_0,y_0\in\widehat{\Sigma}$ such that $\hat{h}_{x_0}$ and $\hat{h}_{y_0}$ lie in different $\mathsf{J}$-orbits.
Since the action is proper, there exists a smooth $\mathsf{J}$-invariant function $F$ on $\mathrm{Sym}^2_+(\mathfrak{n}^*)$ such that $F(\hat{h}_{x_0})\neq F(\hat{h}_{y_0})$; see \cite[Proposition~7]{CrainicMestre18}.
However, since $F$ is $\mathsf{J}$-invariant and $d\hat{h}(T\widehat{\Sigma})$ is tangent to the $\mathsf{J}$-orbits, we have $d(F\circ\hat{h}|_{\widehat{\Sigma}})=0$.
Connectedness of $\widehat{\Sigma}$ implies that $F\circ\hat{h}|_{\widehat{\Sigma}}$ is constant, contradicting $F(\hat{h}_{x_0})\neq F(\hat{h}_{y_0})$.
Thus the image of $\hat{h}|_{\widehat{\Sigma}}$ lies in a single $\mathsf{J}$-orbit.
\end{proof}

The abelianness of $\mathfrak{a}_0$, established in Remark~\ref{rem:horizontal-complement-equivariance}\,(ii), gives $\mathrm{Ad}(\exp(\mathfrak{a}_0))|_{\mathfrak{n}}\subseteq\mathsf{J}$.
In Section~\ref{sec:affine-conjugation-orbits}, we will show that the image of $\hat{h}|_{\widehat{\Sigma}}$ lies in an orbit of this smaller group.
In Section~\ref{sec:solvable-splitting}, we will prove that $\hat{h}$ is constant along $\mathsf{G}$-horizontal curves.
The following remark describes the corresponding isometries of the $\mathsf{N}$-orbits.

\begin{remark}\label{rem:affine-orbit-isometries}
Since $\mathsf{N}$ is simply connected by Corollary~\ref{cor:solvable-exponential-decomposition}\,(i), proved in the next section, every automorphism of $\mathfrak{n}$ integrates uniquely to an automorphism of $\mathsf{N}$.
The image of $\hat{h}|_{\widehat{\Sigma}}$ lying in a $\mathsf{Aut}(\mathfrak{n})$-orbit therefore means that, for every $x_0,y_0\in\widehat{\Sigma}$, there is an automorphism $\varphi\in\mathsf{Aut}(\mathsf{N})$ and a $\varphi$-equivariant isometry $f:\mathsf{N}\cdot x_0\rightarrow\mathsf{N}\cdot y_0$ with $f(x_0)=y_0$.
Explicitly, $\varphi$-equivariance means
$$
    f(n\cdot x)=\varphi(n)\cdot f(x)
    \qquad\text{for }n\in\mathsf{N},\ x\in\mathsf{N}\cdot x_0.
$$
We call such an isometry \emph{affine}.
The successive restrictions on the image of $\hat{h}$ give the following refinements:
\begin{enumerate}
\renewcommand{\labelenumi}{(\roman{enumi})}
    \item Corollary~\ref{cor:centralizer-orbit} allows $\varphi$ to be chosen so that its differential belongs to $\mathsf{J}$.
    Equivalently,
$$
    \varphi(ana^{-1})=a\varphi(n)a^{-1}
    \qquad\text{for }a\in\exp(\mathfrak{a}_0),\ n\in\mathsf{N}.
$$
    \item Theorem~\ref{thm:a0-orbit-image} allows $\varphi$ to be chosen as conjugation by an element of $\exp(\mathfrak{a}_0)$:
$$
    \varphi(n)=ana^{-1}
    \qquad\text{for some }a\in\exp(\mathfrak{a}_0).
$$
    \item By Proposition~\ref{prop:totally-geodesic-orbits}, if $x_0$ and $y_0$ are joined by a $\mathsf{G}$-horizontal curve, then one may take $\varphi=\mathrm{Id}_{\mathsf{N}}$.
    Thus the corresponding $\mathsf{N}$-orbits are $\mathsf{N}$-equivariantly isometric.
\end{enumerate}
\end{remark}

\section{Extending algebraic results from Einstein solvmanifolds}\label{sec:algebraic-rigidity}

Throughout this section, we work under Assumption~\ref{ass:main}.
We extend several algebraic properties of Einstein solvmanifolds to the solvable cocompact setting.
The commutativity relations established in Section~\ref{sec:shape-commutativity} provide the main algebraic input.

\begin{theorem}\label{thm:solvable-algebraic-rigidity}
Let $M$ be a complete, connected negative-Einstein manifold and let $\mathsf{G}$ be a connected solvable Lie group acting freely, properly, isometrically and cocompactly on $M$.
Then the following properties hold:
\begin{enumerate}
\renewcommand{\labelenumi}{(\roman{enumi})}
    \item For every $x\in M$ and $A\in\mathfrak{a}_x$, the endomorphism $\mathrm{ad}_A$ is $h_x$-skew-symmetric if and only if $A=0$.
    \item $\mathfrak{g}$ is of real type: for every $U\in\mathfrak{g}\setminus\mathfrak{n}$, the endomorphism $\mathrm{ad}_U$ has an eigenvalue with nonzero real part.
    \item $\mathsf{G}$ is simply connected and centerless.
    \item For every $x\in M$ and $g\in\mathsf{G}$, we have $\mathrm{Ad}(g)|_{\mathfrak{n}}\in\mathsf{O}(\mathfrak{n},\hat{h}_x)$ if and only if $g\in Z(\mathsf{N})$.
\end{enumerate}
\end{theorem}

For Einstein solvmanifolds, item~(i) is \cite[Theorem~4.10(3)]{Heber98}, while item~(ii) follows from Heber's structure results and is stated in \cite[Proposition~2.8]{Jab15}.
Heber assumes standardness, which holds for every Einstein solvmanifold by \cite{standard}.
Simple connectedness in item~(iii) was proved in the more general setting of nonflat Ricci soliton solvmanifolds in \cite[Theorem~1.1]{Jbl2015}; centerlessness also follows from Heber's structure results.
The argument for item~(iv) adapts \cite[Section~5]{ThompsonSimplyConnected}, where orthogonality is assumed on the whole Lie algebra.

These properties will be used at different stages of the proof of Theorem~\ref{thm:main-splitting}.
In view of Lemma~\ref{lem:orbit-ricci-a-directions}, item~(i) gives $\mathrm{ric}^{\mathsf{G}}(A_x,A_x)<0$ for every nonzero $A_x\in\mathcal{A}_x$; it will also give the injectivity needed to construct the bundle map $\alpha$ in Proposition~\ref{prop:shape-alpha}.
Item~(ii) will be used in the proofs of items~(iii) and~(iv).
Item~(iii) gives the global exponential decomposition in Corollary~\ref{cor:solvable-exponential-decomposition} and eliminates central elements in the holonomy argument of Section~\ref{sec:trivial-n-holonomy}.
Finally, item~(iv) will be used to prove that the $\mathsf{A}_0$-action on $\mathrm{Sym}^2_+(\mathfrak{n}^*)$ is free and proper in Lemma~\ref{lem:a0-action-free-proper}.
This freeness will then eliminate the $\mathsf{G}$-holonomy in Proposition~\ref{prop:trivial-g-holonomy}.

We first record two differential identities and prove items~(i)--(iii).
We then introduce an orthogonal decomposition of the nilradical and use it to prove item~(iv).

For every $x\in M$, the orbit $\mathsf{G}\cdot x$ is isometric to $\mathsf{G}$ endowed with the left-invariant metric determined by $h_x$; see Subsection~\ref{subsec:trivialization}.
Applying Lauret's Ricci formula \cite[equation~(25)]{Lauret11} and using the abelianness of $\mathfrak{a}_x$, we obtain the following identity by polarization.

\begin{lemma}\label{lem:orbit-ricci-a-directions}
Let $A,\tilde{A}$ be $\mathcal{A}$-valued vector fields on $M$.
Then
$$
    \mathrm{ric}^{\mathsf{G}}(A,\tilde{A})
    =-\mathrm{Tr}\bigl(
    \mathrm{S}(\mathrm{ad}_{{}^{\raisebox{0.15em}{$\scriptstyle\mathfrak{g}$}}\! A})
    \cdot\mathrm{S}(\mathrm{ad}_{{}^{\raisebox{0.15em}{$\scriptstyle\mathfrak{g}$}}\!\tilde{A}})
    \bigr)
    =-\mathrm{Tr}(\hat{L}_A\cdot\hat{L}_{\tilde{A}}).
$$
If $A$ and $\tilde{A}$ are also invariant and horizontally Killing, then $\mathrm{ric}^{\mathsf{G}}(A,\tilde{A})$ is constant on $M$.

\end{lemma}

\begin{proof}
The first equality is Lauret's Ricci formula recalled above.
By Remark~\ref{rem:horizontal-complement-equivariance}\,(ii), the endomorphisms $\mathrm{ad}_{{}^{\raisebox{0.15em}{$\scriptstyle\mathfrak{g}$}}\! A}$ and $\mathrm{ad}_{{}^{\raisebox{0.15em}{$\scriptstyle\mathfrak{g}$}}\!\tilde{A}}$ vanish on $\mathfrak{a}_x$ and preserve $\mathfrak{n}$.
Lemma~\ref{lem:shape-symmetric-adjoint} therefore gives the second equality.

Suppose that $A$ and $\tilde{A}$ are invariant and horizontally Killing.
By Corollary~\ref{cor:parallel-a-shape}, their shape operators are parallel for the adjoint connection of the $\mathsf{N}$-action.
Hence, for every $\mathsf{N}$-horizontal vector field $\widehat{X}$, the product rule and invariance of the trace give
$$
    \widehat{X}\bigl(\mathrm{ric}^{\mathsf{G}}(A,\tilde{A})\bigr)
    =-\mathrm{Tr}\bigl(
    (\widehat{\nabla}^{\mathrm{ad}}_{\widehat{X}}\hat{L}_A)\cdot\hat{L}_{\tilde{A}}
    +\hat{L}_A\cdot(\widehat{\nabla}^{\mathrm{ad}}_{\widehat{X}}\hat{L}_{\tilde{A}})
    \bigr)
    =0.
$$
Moreover, $\mathrm{ric}^{\mathsf{G}}(A,\tilde{A})$ is $\mathsf{G}$-invariant.
Thus, its differential vanishes everywhere.
Connectedness of $M$ proves constancy.
\end{proof}

The next identity extends \cite[Lemma~8.1]{bohm2025non} to invariant horizontally Killing vector fields.
As in Subsection~\ref{subsec:divergence-maximum-principle}, we regard functions on $B:=M/\mathsf{G}$ as invariant functions on $M$.

\begin{lemma}\label{lem:horizontally-killing-norm}
Let $U$ be a vertical, invariant and horizontally Killing vector field on $M$, and set $f:=\|U\|^2/2$.
Then $f$ is invariant and
$$
    \Delta_Bf-\ip{\nabla f}{H}
    =\mathrm{ric}^{\mathsf{G}}(U,U)+\|U\|^2+2\|LU\|^2+\|\mathrm{A}U\|^2,
$$
where, for a local orthonormal horizontal frame $\{X_i\}$,
$$
    \|LU\|^2:=\sum_i\|L_{X_i}U\|^2,
    \qquad
    \|\mathrm{A}U\|^2:=\sum_i\|\mathrm{A}_{X_i}U\|^2.
$$
\end{lemma}

\begin{proof}
The proof follows the same computation as \cite[Lemma~8.1]{bohm2025non}, using $\mathcal{V}(\nabla_XU)=L_XU$ for every horizontal vector field $X$, which characterizes horizontally Killing fields by Definition~\ref{def:horizontally-killing}\,(iv).
The additional term $\|\mathrm{A}U\|^2$ is retained from the general vertical Ricci formula in Proposition~\ref{prop:oneill-vertical-ricci}.
\end{proof}

\begin{proof}[Proof of Theorem~\ref{thm:solvable-algebraic-rigidity}\,(i)]
Fix $x_0\in M$ and $A_0\in\mathfrak{a}_{x_0}$ such that $\mathrm{S}_{x_0}(\mathrm{ad}_{A_0})=0$.
By Proposition~\ref{prop:invariant-a-extension}, there exists a unique invariant, horizontally Killing, $\mathcal{A}$-valued vector field $A$ such that ${}^{\mathfrak{g}}A_{x_0}=A_0$.
Lemma~\ref{lem:orbit-ricci-a-directions} gives $\mathrm{ric}^{\mathsf{G}}(A,A)(x_0)=0$ and shows that $\mathrm{ric}^{\mathsf{G}}(A,A)$ is constant on $M$.
Thus, $\mathrm{ric}^{\mathsf{G}}(A,A)=0$ on $M$.

Applying Lemma~\ref{lem:horizontally-killing-norm} to the invariant function $f:=\|A\|^2/2$, we obtain
$$
    \Delta_Bf-\ip{\nabla f}{H}
    =\|A\|^2+2\|LA\|^2+\|\mathrm{A}A\|^2\geq0.
$$
The maximum principle in Proposition~\ref{prop:cocompact-maximum-principle}\,(ii) implies that $f$ is constant.
The left-hand side therefore vanishes, so $A\equiv0$.
In particular, $A_0=0$.
\end{proof}

We use the terminology \emph{real type} from \cite[Definition~2.1]{BL19}; the equivalent condition is called \emph{almost-complete solvability} in \cite[Definition~1.6]{Jab15}.

\begin{proof}[Proof of Theorem~\ref{thm:solvable-algebraic-rigidity}\,(ii)]
Fix $x\in M$ and write $U=A+\widehat{U}$, with $A\in\mathfrak{a}_x$ and $\widehat{U}\in\mathfrak{n}$.
As observed in \cite[Remark~4.3(c)]{Heber98}, $\mathrm{ad}_U$ and $\mathrm{ad}_A$ have the same eigenvalues, counted with multiplicity.
Indeed, Lie's theorem gives a basis of the complexification $\mathfrak{g}^{\mathbb{C}}$ in which all adjoint endomorphisms are triangular.
Since $\widehat{U}\in\mathfrak{n}$, the endomorphism $\mathrm{ad}_{\widehat{U}}$ is nilpotent and has zero diagonal in this basis.
Thus, $\mathrm{ad}_U$ and $\mathrm{ad}_A$ have the same diagonal entries.

Suppose that every eigenvalue of $\mathrm{ad}_U$ has zero real part.
The same is then true of $\mathrm{ad}_A$.
By Corollary~\ref{cor:normal-adjoint-derivations}\,(i), $\mathrm{ad}_A$ is $h_x$-normal.
A real normal endomorphism whose eigenvalues all have zero real part is skew-symmetric.
Item~(i) therefore gives $A=0$, and hence $U\in\mathfrak{n}$.
This proves item~(ii).
\end{proof}

In Lemma~\ref{lem:discrete-center}, we will show under more general assumptions that $\mathfrak{z}(\mathfrak{g})=\{0\}$.
We use this fact in the following proof.

\begin{proof}[Proof of Theorem~\ref{thm:solvable-algebraic-rigidity}\,(iii)]
Let $r:\widetilde{\mathsf{G}}\rightarrow\mathsf{G}$ be the universal covering homomorphism, and identify its Lie algebra with $\mathfrak{g}$ through $dr_e$.
Fix $x\in M$.
As established in the proof of Proposition~\ref{prop:invariant-a-extension}, the map
$$
    \mathfrak{a}_x\times\mathfrak{n}\longrightarrow\widetilde{\mathsf{G}},
    \qquad
    (A,\widehat{U})\longmapsto\exp(A)\exp(\widehat{U})
$$
is a diffeomorphism.

Take $\tilde{g}\in Z(\widetilde{\mathsf{G}})$ and write $\tilde{g}=\exp(A)\exp(\widehat{U})$.
Since $\mathrm{Ad}(\tilde{g})=\mathrm{Id}_{\mathfrak{g}}$,
$$
    e^{\mathrm{ad}_A}e^{\mathrm{ad}_{\widehat{U}}}=\mathrm{Id}_{\mathfrak{g}},
    \qquad
    e^{\mathrm{ad}_A}=e^{-\mathrm{ad}_{\widehat{U}}}.
$$
The right-hand side is unipotent because $\mathrm{ad}_{\widehat{U}}$ is nilpotent.
The left-hand side is $h_x$-normal by Corollary~\ref{cor:normal-adjoint-derivations}\,(i).
A normal unipotent endomorphism is the identity, since it is diagonalizable over $\mathbb{C}$ and all its eigenvalues are $1$.
Consequently,
$$
    e^{\mathrm{ad}_A}=\mathrm{Id}_{\mathfrak{g}}
    =e^{\mathrm{ad}_{\widehat{U}}}.
$$
Every eigenvalue $\lambda$ of $\mathrm{ad}_A$ therefore satisfies $e^\lambda=1$, so $\mathrm{Re}(\lambda)=0$.
By item~(ii), we conclude that $A=0$.
Since $\mathrm{ad}_{\widehat{U}}$ is nilpotent and its exponential is the identity, it must vanish.
Thus $\widehat{U}\in\mathfrak{z}(\mathfrak{g})=\{0\}$, and $\tilde{g}=e$.

We have proved that $Z(\widetilde{\mathsf{G}})=\{e\}$.
The kernel of $r$ is a discrete central subgroup of $\widetilde{\mathsf{G}}$, so it is trivial.
Hence $r$ is an isomorphism, and $\mathsf{G}$ is simply connected and centerless.
\end{proof}

By Theorem~\ref{thm:solvable-algebraic-rigidity}\,(iii), $\mathsf{G}$ is simply connected, so the Lie algebra splitting $\mathfrak{g}=\mathfrak{a}_x\ltimes\mathfrak{n}$ integrates to a semidirect product of closed, simply connected subgroups; see \cite[Theorem~1.125]{Knapp02}.
The exponential maps of the abelian and nilpotent factors are diffeomorphisms by \cite[Theorem~1.127]{Knapp02}.
Normality of $\mathsf{N}$ allows either order of the factors, giving the following.

\begin{corollary}\label{cor:solvable-exponential-decomposition}
For every $x\in M$, the following hold:
\begin{enumerate}
\renewcommand{\labelenumi}{(\roman{enumi})}
    \item $\mathsf{N}$ is simply connected, and the exponential map $\exp:\mathfrak{n}\to\mathsf{N}$ is a diffeomorphism.
    \item $\exp(\mathfrak{a}_x)$ is a closed, connected, simply connected abelian subgroup of $\mathsf{G}$, and $\exp:\mathfrak{a}_x\to\exp(\mathfrak{a}_x)$ is a diffeomorphism.
    \item We have $\mathsf{G}=\exp(\mathfrak{a}_x)\ltimes\mathsf{N}$.
    In particular, the map
$$
    \mathfrak{a}_x\times\mathfrak{n}\longrightarrow\mathsf{G},
    \qquad
    (A,\widehat{U})\longmapsto\exp(A)\cdot\exp(\widehat{U})
$$
    is a diffeomorphism.
\end{enumerate}
\end{corollary}

\subsection{Orthogonal adjoint actions on the nilradical}\label{subsec:algebraic-rigidity-orthogonal-action}

\textcolor{white}{aaa}

To prove item~(iv), we will use the following lemma, which gives an orthogonal decomposition of $\mathfrak{n}$ obtained from its lower central series.

\begin{lemma}\label{lem:nilradical-orthogonal-filtration}
For every $x\in M$, there is an $h_x$-orthogonal decomposition $\mathfrak{n}=\mathfrak{n}_0\oplus\cdots\oplus\mathfrak{n}_r$ such that
\begin{enumerate}
\renewcommand{\labelenumi}{(\roman{enumi})}
    \item $[\mathfrak{n},\mathfrak{n}_i]\subseteq\mathfrak{n}_{i+1}\oplus\cdots\oplus\mathfrak{n}_r$ for every $i$;
    \item $[\mathfrak{a}_x,\mathfrak{n}_i]\subseteq\mathfrak{n}_i$ for every $i$.
\end{enumerate}
\end{lemma}

\begin{proof}
Fix $x\in M$.
Consider the lower central series of $\mathfrak{n}$, defined by $c_0:=\mathfrak{n}$ and $c_{i+1}:=[\mathfrak{n},c_i]$.
Since $\mathfrak{n}$ is nilpotent, there exists $r$ such that $c_r\neq\{0\}$ and $c_{r+1}=\{0\}$.
Define $\mathfrak{n}_i$ to be the $h_x$-orthogonal complement of $c_{i+1}$ in $c_i$.
Thus,
$$
    c_i=\mathfrak{n}_i\oplus c_{i+1}
    \qquad\text{orthogonally},
$$
and, inductively, $c_i=\mathfrak{n}_i\oplus\cdots\oplus\mathfrak{n}_r$.
In particular, $\mathfrak{n}=\mathfrak{n}_0\oplus\cdots\oplus\mathfrak{n}_r$ is an orthogonal decomposition.
Moreover,
$$
    [\mathfrak{n},\mathfrak{n}_i]
    \subseteq[\mathfrak{n},c_i]
    =c_{i+1}
    =\mathfrak{n}_{i+1}\oplus\cdots\oplus\mathfrak{n}_r,
$$
which proves item~(i).

Every derivation of $\mathfrak{n}$ preserves its lower central series.
For $A\in\mathfrak{a}_x$, the restriction $\mathrm{ad}_A|_{\mathfrak{n}}$ is a derivation of $\mathfrak{n}$.
As explained at the beginning of Section~\ref{sec:shape-commutativity}, its $\hat{h}_x$-transpose is also a derivation by \cite[Proposition~9.7]{BL23} and Lemma~\ref{lem:shape-symmetric-adjoint}; this also follows from Corollary~\ref{cor:normal-adjoint-derivations}\,(ii).
Hence both $\mathrm{ad}_A|_{\mathfrak{n}}$ and $\mathrm{ad}_A^{t_x}|_{\mathfrak{n}}$ preserve every $c_i$.

Let $\widehat{U}\in\mathfrak{n}_i$.
Since $\mathrm{ad}_A$ preserves $c_i$, we have $[A,\widehat{U}]\in c_i$.
On the other hand,
$$
    h_x([A,\widehat{U}],\widehat{V})
    =h_x(\widehat{U},\mathrm{ad}_A^{t_x}\widehat{V})=0
    \qquad\text{for every }\widehat{V}\in c_{i+1},
$$
because $\mathrm{ad}_A^{t_x}$ preserves $c_{i+1}$ and $\widehat{U}$ is orthogonal to $c_{i+1}$.
Thus $[A,\widehat{U}]$ belongs to the orthogonal complement of $c_{i+1}$ in $c_i$, which is $\mathfrak{n}_i$.
This proves item~(ii).
\end{proof}

\begin{proof}[Proof of Theorem~\ref{thm:solvable-algebraic-rigidity}\,(iv)]
Fix $x\in M$ and take the decomposition from Lemma~\ref{lem:nilradical-orthogonal-filtration}.
Let $g\in\mathsf{G}$ be such that $\mathrm{Ad}(g)|_{\mathfrak{n}}$ is $\hat{h}_x$-orthogonal.
By Corollary~\ref{cor:solvable-exponential-decomposition}\,(iii), we may write $g=\exp(A)\exp(\widehat{U})$, with $A\in\mathfrak{a}_x$ and $\widehat{U}\in\mathfrak{n}$.
Hence
\begin{equation}\label{eq:Ad-decomposition}
    \mathrm{Ad}(g)|_{\mathfrak{n}}
    =e^{\mathrm{ad}_A|_{\mathfrak{n}}}
    e^{\mathrm{ad}_{\widehat{U}}|_{\mathfrak{n}}}.
\end{equation}

By Lemma~\ref{lem:nilradical-orthogonal-filtration}\,(ii), $\mathrm{ad}_A$ preserves each $\mathfrak{n}_i$, so $e^{\mathrm{ad}_A|_{\mathfrak{n}}}$ does as well.
By item~(i) of the same lemma, $\mathrm{ad}_{\widehat{U}}$ maps $\mathfrak{n}_i$ into $\mathfrak{n}_{i+1}\oplus\cdots\oplus\mathfrak{n}_r$.
Thus its exponential preserves every subspace $\mathfrak{n}_i\oplus\cdots\oplus\mathfrak{n}_r$. 
Equation~\eqref{eq:Ad-decomposition} now shows that $\mathrm{Ad}(g)|_{\mathfrak{n}}$ preserves $\mathfrak{n}_i\oplus\cdots\oplus\mathfrak{n}_r$ for every $i$.

Since $\mathrm{Ad}(g)|_{\mathfrak{n}}$ is orthogonal, it also preserves the orthogonal complement of $\mathfrak{n}_{i+1}\oplus\cdots\oplus\mathfrak{n}_r$ in $\mathfrak{n}_i\oplus\cdots\oplus\mathfrak{n}_r$, namely $\mathfrak{n}_i$.
Therefore,
$$
    \mathrm{Ad}(g)(\mathfrak{n}_i)=\mathfrak{n}_i
    \qquad\text{for every }i.
$$

Rearranging~\eqref{eq:Ad-decomposition}, we obtain
\begin{equation}\label{eq:Ad-nilpotent}
    e^{-\mathrm{ad}_A|_{\mathfrak{n}}}\mathrm{Ad}(g)|_{\mathfrak{n}}
    -\mathrm{Id}_{\mathfrak{n}}
    =e^{\mathrm{ad}_{\widehat{U}}|_{\mathfrak{n}}}-\mathrm{Id}_{\mathfrak{n}}.
\end{equation}
The left-hand side preserves each $\mathfrak{n}_i$, because both $e^{-\mathrm{ad}_A|_{\mathfrak{n}}}$ and $\mathrm{Ad}(g)|_{\mathfrak{n}}$ do.
By Lemma~\ref{lem:nilradical-orthogonal-filtration}\,(i), every positive power of $\mathrm{ad}_{\widehat{U}}|_{\mathfrak{n}}$ maps $\mathfrak{n}_i$ into $\mathfrak{n}_{i+1}\oplus\cdots\oplus\mathfrak{n}_r$.
The right-hand side of~\eqref{eq:Ad-nilpotent} therefore has the same property.
Since $\mathfrak{n}_i$ is orthogonal to $\mathfrak{n}_{i+1}\oplus\cdots\oplus\mathfrak{n}_r$, both sides vanish on each summand.
Consequently,
\begin{equation}\label{eq:Ad-conclusion}
    e^{\mathrm{ad}_{\widehat{U}}|_{\mathfrak{n}}}
    =\mathrm{Id}_{\mathfrak{n}},
    \qquad
    e^{\mathrm{ad}_A|_{\mathfrak{n}}}
    =\mathrm{Ad}(g)|_{\mathfrak{n}}.
\end{equation}

Since $e^{\mathrm{ad}_A|_{\mathfrak{n}}}$ is orthogonal, every eigenvalue $\lambda$ of $\mathrm{ad}_A|_{\mathfrak{n}}$ satisfies $1=|e^\lambda|=e^{\mathrm{Re}(\lambda)}$.
Thus every eigenvalue of $\mathrm{ad}_A|_{\mathfrak{n}}$ has zero real part.
Since $\mathrm{ad}_A$ vanishes on $\mathfrak{a}_x$, the same holds for all eigenvalues of $\mathrm{ad}_A$.
Theorem~\ref{thm:solvable-algebraic-rigidity}\,(ii) gives $A=0$.

It follows that $g=\exp(\widehat{U})$.
Since $\mathrm{ad}_{\widehat{U}}|_{\mathfrak{n}}$ is nilpotent and its exponential is the identity by~\eqref{eq:Ad-conclusion}, we have $\mathrm{ad}_{\widehat{U}}|_{\mathfrak{n}}=0$.
Thus $\widehat{U}\in\mathfrak{z}(\mathfrak{n})$, and hence $g\in Z(\mathsf{N})$.
\end{proof}

\section{The $\mathsf{N}$-holonomy is trivial}\label{sec:trivial-n-holonomy}
Throughout this section, we work under Assumption~\ref{ass:main}.
Recall that the holonomy group $\mathsf{Hol}_x\subseteq\mathsf{G}$ of the principal bundle $M\to M/\mathsf{G}$ consists of those elements $g\in\mathsf{G}$ for which $x$ and $g\cdot x$ can be joined by a piecewise smooth $\mathsf{G}$-horizontal curve.
We denote the holonomy group of $M\to M/\mathsf{N}$ by $\widehat{\mathsf{Hol}}_x\subseteq\mathsf{N}$.

The polarity of the $\mathsf{N}$-action (Theorem~\ref{thm:nilradical-polarity}) implies that $\mathrm{Lie}(\widehat{\mathsf{Hol}}_x)=\{0\}$, by Proposition~\ref{prop:discrete-holonomy}.
Moreover, Lemma~\ref{lem:vertical-tensor-blocks}, equation~\eqref{eq:holonomy-algebra} and the constancy of $\mathfrak{a}_x$ along horizontal curves give $\mathrm{Lie}(\mathsf{Hol}_x)\subseteq\mathfrak{a}_x$ for every $x\in M$.
The following theorem strengthens these infinitesimal conclusions to statements about the full holonomy groups.

\begin{theorem}\label{thm:trivial-n-holonomy}
Under Assumption~\ref{ass:main}, $\widehat{\mathsf{Hol}}_x=\{e\}$ and $\mathsf{Hol}_x\subseteq\exp(\mathfrak{a}_x)$ for every $x\in M$.
\end{theorem}

The main application of this theorem is the corollary below.
By Theorem~\ref{thm:nilradical-polarity}, the $\mathsf{N}$-horizontal distribution is integrable; we call its maximal connected integral manifolds $\mathsf{N}$-horizontal leaves.
By Corollary~\ref{cor:horizontal-complement-subalgebra}, the map $x\mapsto\mathfrak{a}_x$ is constant along each such leaf.

\begin{corollary}\label{cor:a0-action-on-section}
Let $\widehat{\Sigma}\subseteq M$ be an $\mathsf{N}$-horizontal leaf, and set $\mathfrak{a}_0:=\mathfrak{a}_x$ for any $x\in\widehat{\Sigma}$.
Then the action of $\mathsf{A}_0:=\exp(\mathfrak{a}_0)\subseteq\mathsf{G}$ restricts to an action on $\widehat{\Sigma}$ with compact orbit space $\widehat{\Sigma}/\mathsf{A}_0\equiv M/\mathsf{G}$.
\end{corollary}

\begin{proof}
Let $x\in\widehat{\Sigma}$ and $A\in\mathfrak{a}_0$.
By Lemma~\ref{lem:exponential-horizontal-curve} below, the curve $t\mapsto\exp(tA)\cdot x$ is $\mathsf{N}$-horizontal.
It therefore lies in $\widehat{\Sigma}$, and in particular $\exp(A)\cdot x\in\widehat{\Sigma}$.
Thus $\mathsf{A}_0$ preserves $\widehat{\Sigma}$.

Let $\pi:M\to M/\mathsf{N}$ and $\overline{\pi}:\mathsf{G}\to\mathsf{G}/\mathsf{N}$ denote the respective projections.
Since the $\mathsf{N}$-holonomy is trivial by Theorem~\ref{thm:trivial-n-holonomy}, $\pi|_{\widehat{\Sigma}}:\widehat{\Sigma}\to M/\mathsf{N}$ is an isometry.
By Corollary~\ref{cor:solvable-exponential-decomposition}\,(iii), the restriction $\overline{\pi}|_{\mathsf{A}_0}:\mathsf{A}_0\to\mathsf{G}/\mathsf{N}$ is a Lie group isomorphism.
Moreover,
$$
    \pi|_{\widehat{\Sigma}}(a\cdot x)
    =\overline{\pi}(a)\cdot\pi|_{\widehat{\Sigma}}(x)
    \qquad\text{for every }a\in\mathsf{A}_0,\ x\in\widehat{\Sigma}.
$$
Thus, in the category of Riemannian manifolds equipped with isometric group actions, the pair $(\widehat{\Sigma},\mathsf{A}_0)$ is isomorphic to the pair $(M/\mathsf{N},\mathsf{G}/\mathsf{N})$.
In particular,
$$
    \widehat{\Sigma}/\mathsf{A}_0
    \equiv\frac{M/\mathsf{N}}{\mathsf{G}/\mathsf{N}}
    \equiv M/\mathsf{G},
$$
which is compact by Assumption~\ref{ass:main}.
\end{proof}

Before proving Theorem~\ref{thm:trivial-n-holonomy}, we establish two facts about curves.
The first decomposes a curve into a group component and a $\mathsf{G}$-horizontal curve, while the second establishes the $\mathsf{N}$-horizontality of curves of the form $t\mapsto\exp(A(t))\cdot\sigma(t)$.

\begin{lemma}\label{lem:horizontal-decomposition}
Let $\mathsf{G}$ be a Lie group acting freely, properly and isometrically on a complete Riemannian manifold $M$.
Then, for every piecewise smooth curve $\sigma:[0,1]\to M$, there exist a unique piecewise smooth $\mathsf{G}$-horizontal curve $\overline{\sigma}:[0,1]\to M$ and a unique piecewise smooth curve $g:[0,1]\to\mathsf{G}$ with $g(0)=e$ such that $\sigma(t)=g(t)\cdot\overline{\sigma}(t)$ for all $t$.
\end{lemma}

\begin{proof}
Let $\pi:M\to M/\mathsf{G}$ be the projection.
By Proposition~\ref{prop:horizontal-lifting}, the curve $\pi\circ\sigma$ has a unique piecewise smooth $\mathsf{G}$-horizontal lift $\overline{\sigma}$ with $\overline{\sigma}(0)=\sigma(0)$.
Since $\pi\circ\sigma=\pi\circ\overline{\sigma}$, the points $\sigma(t)$ and $\overline{\sigma}(t)$ lie in the same $\mathsf{G}$-orbit, and freeness gives a unique $g(t)\in\mathsf{G}$ such that $\sigma(t)=g(t)\cdot\overline{\sigma}(t)$.
Local trivializations of the principal bundle show that $g$ is piecewise smooth.
Moreover, $g(0)=e$ because $\sigma(0)=\overline{\sigma}(0)$.
Uniqueness of $\overline{\sigma}$ follows from Proposition~\ref{prop:horizontal-lifting}, and uniqueness of $g$ then follows from freeness.
\end{proof}

We return to Assumption~\ref{ass:main}.
\begin{lemma}\label{lem:exponential-horizontal-curve}
Let $\sigma:[0,1]\to M$ be a piecewise smooth $\mathsf{N}$-horizontal curve with $\sigma(0)=x$, and let $A:[0,1]\to\mathfrak{a}_x$ be a piecewise smooth curve.
Then $t\mapsto\exp(A(t))\cdot\sigma(t)$ is an $\mathsf{N}$-horizontal curve in $M$.
\end{lemma}

\begin{proof}
Let $\overline{m}_g:\mathsf{G}\to \mathsf{G}$ denote left multiplication by $g$, and set $a(t):=\exp(A(t))$.
Since $\mathfrak{a}_x$ is abelian (Remark~\ref{rem:horizontal-complement-equivariance}\,(ii)), we have
$$
    \exp(A(t+s))=\exp(A(t))\cdot\exp(A(t+s)-A(t))\qquad\text{for all }t,t+s\in[0,1].
$$
That is,
$$
    a(t+s)=\overline{m}_{a(t)}\bigl(\exp(A(t+s)-A(t))\bigr).
$$
Differentiating at $s=0$ and using $d\exp_0=\mathrm{Id}_{\mathfrak{g}}$, we obtain
\begin{equation}\label{eq:derivative-exponential-curve}
    a'(t)=d(\overline{m}_{a(t)})_e\bigl(A'(t)\bigr).
\end{equation}
We now use this identity to differentiate the curve $\widetilde{\sigma}(t):=\exp(A(t))\cdot\sigma(t)=a(t)\cdot\sigma(t)$.
Recall that $m_g:M\to M$ denotes the map $x\mapsto g\cdot x$, while $r_x:\mathsf{G}\to M$ denotes the map $g\mapsto g\cdot x$.
We have
\begin{equation}\label{eq:derivative-action-curve}
    \widetilde{\sigma}'(t)=d(m_{a(t)})_{\sigma(t)}\bigl(\sigma'(t)\bigr)
    +d(r_{\sigma(t)})_{a(t)}\bigl(a'(t)\bigr).
\end{equation}
Associativity of the action gives the identity $r_x\circ\overline{m}_g=m_g\circ r_x$.
Note that $d(r_x)_e(U)=U_x^*$ by the definition of the fundamental vector field.
Thus, using equation~\eqref{eq:derivative-exponential-curve},
\begin{align*}
    d(r_{\sigma(t)})_{a(t)}\bigl(a'(t)\bigr)
    &=d(r_{\sigma(t)})_{a(t)}\bigl(d(\overline{m}_{a(t)})_e(A'(t))\bigr)=
    d(r_{\sigma(t)}\circ\overline{m}_{a(t)})_e\bigl(A'(t)\bigr)\\
    &=d(m_{a(t)}\circ r_{\sigma(t)})_e\bigl(A'(t)\bigr)
    =d(m_{a(t)})_{\sigma(t)}\bigl((A'(t))_{\sigma(t)}^*\bigr).
\end{align*}

Substituting this into equation~\eqref{eq:derivative-action-curve}, we obtain
\begin{equation}\label{eq:horizontal-derivative}
    \widetilde{\sigma}'(t)=d(m_{a(t)})_{\sigma(t)}\bigl(\sigma'(t)+(A'(t))_{\sigma(t)}^*\bigr).
\end{equation}
Since $\sigma$ is $\mathsf{N}$-horizontal, Corollary \ref{cor:horizontal-complement-subalgebra} gives $\mathfrak{a}_{\sigma(t)}=\mathfrak{a}_x$.
Hence $A'(t)\in\mathfrak{a}_x=\mathfrak{a}_{\sigma(t)}$, so $(A'(t))_{\sigma(t)}^*\in\mathcal{A}_{\sigma(t)}$.
In particular,
$$
    \sigma'(t)+(A'(t))_{\sigma(t)}^*\in\widehat{\mathcal{H}}_{\sigma(t)}\qquad\text{for all }t.
$$
By Remark~\ref{rem:preservation-n-horizontal} and equation~\eqref{eq:horizontal-derivative}, we conclude that $\widetilde{\sigma}'(t)\in\widehat{\mathcal{H}}_{\widetilde{\sigma}(t)}$, as desired.
\end{proof}

\begin{remark}\label{rem:n-horizontal-translates}
By Remark~\ref{rem:preservation-n-horizontal}, if $\sigma$ is an $\mathsf{N}$-horizontal curve and $g\in\mathsf{G}$, then $t\mapsto g\cdot\sigma(t)$ is also an $\mathsf{N}$-horizontal curve.
\end{remark}

We now turn to the proof of Theorem~\ref{thm:trivial-n-holonomy}.
We first show that the $\mathsf{N}$-holonomy is central, using the $\mathsf{N}$-vertical endomorphism $\beta^+$ introduced in Subsection~\ref{subsec:git}.

\begin{lemma}\label{lem:central-n-holonomy}
$\widehat{\mathsf{Hol}}_x\subseteq Z(\mathsf{N})$ for every $x\in M$.
\end{lemma}

\begin{proof}
Throughout this proof, all adjoint endomorphisms are restricted to $\mathfrak{n}$.
By Theorem~\ref{thm:beta-properties}\,(i),(iv) and Lemma~\ref{lem:fundamental-field-equivariance}, the map ${}^{\mathfrak{n}}\beta^+:M\to\mathrm{End}(\mathfrak{n})$ is $\mathsf{N}$-equivariant and constant along $\mathsf{N}$-horizontal curves.
Let $n_0\in\widehat{\mathsf{Hol}}_x$.
Since $x$ and $n_0\cdot x$ can be joined by an $\mathsf{N}$-horizontal curve, we have
$$
    \mathrm{Ad}(n_0)\cdot{}^{\mathfrak{n}}\beta_x^+
    \cdot\mathrm{Ad}(n_0)^{-1}
    ={}^{\mathfrak{n}}\beta_{n_0\cdot x}^+
    ={}^{\mathfrak{n}}\beta_x^+.
$$
Thus $[\mathrm{Ad}(n_0),{}^{\mathfrak{n}}\beta_x^+]=0$.

By Corollary~\ref{cor:solvable-exponential-decomposition}\,(i), write $n_0=\exp(\widehat{U})$ with $\widehat{U}\in\mathfrak{n}$, so that $\mathrm{Ad}(n_0)=e^{\mathrm{ad}_{\widehat{U}}}$.
It follows that $[e^{k\cdot\mathrm{ad}_{\widehat{U}}},{}^{\mathfrak{n}}\beta_x^+]=0$ for every $k\in\mathbb{N}$.
Since $\mathfrak{n}$ is nilpotent, there exists $r\in\mathbb{N}$ such that $(\mathrm{ad}_{\widehat{U}})^r=0$ on $\mathfrak{n}$.
Hence
$$
    P(t):=[e^{t\cdot\mathrm{ad}_{\widehat{U}}},{}^{\mathfrak{n}}\beta_x^
    +]=t\cdot[\mathrm{ad}_{\widehat{U}},{}^{\mathfrak{n}}\beta_x^
    +]+\frac{1}{2}t^2[(\mathrm{ad}_{\widehat{U}})^2,{}^{\mathfrak{n}}\beta_x^
    +]+\cdots+\frac{1}{(r-1)!}t^{r-1}[(\mathrm{ad}_{\widehat{U}})^{r-1},{}^{\mathfrak{n}}\beta_x^
    +]
$$
is a matrix whose entries are polynomials in $t$.
Since $P(k)=0$ for every $k\in\mathbb{N}$, each entry has infinitely many roots and therefore vanishes identically.
Thus $P\equiv0$, and hence $[\mathrm{ad}_{\widehat{U}},{}^{\mathfrak{n}}\beta_x^+]=0$.

By Theorem~\ref{thm:beta-properties}\,(iii), ${}^{\mathfrak{n}}\beta_x^+\in\mathrm{Der}(\mathfrak{n})$, so $\mathrm{ad}_{{}^{\mathfrak{n}}\beta_x^+(\widehat{U})}=[{}^{\mathfrak{n}}\beta_x^+,\mathrm{ad}_{\widehat{U}}]=0$.
Thus ${}^{\mathfrak{n}}\beta_x^+(\widehat{U})\in\mathfrak{z}(\mathfrak{n})$.
Since derivations preserve the center, ${}^{\mathfrak{n}}\beta_x^+(\mathfrak{z}(\mathfrak{n}))\subseteq\mathfrak{z}(\mathfrak{n})$.
Moreover, ${}^{\mathfrak{n}}\beta_x^+$ is positive definite by Theorem~\ref{thm:beta-properties}\,(i), and hence invertible.
Its restriction to $\mathfrak{z}(\mathfrak{n})$ is therefore an isomorphism of this subspace, so $\widehat{U}\in\mathfrak{z}(\mathfrak{n})$.
Consequently, $n_0=\exp(\widehat{U})\in Z(\mathsf{N})$.
\end{proof}

\begin{proof}[Proof of Theorem~\ref{thm:trivial-n-holonomy}]
We first show that $\mathsf{Hol}_x\subseteq\exp(\mathfrak{a}_x)$.
Normality of $\mathsf{N}$ allows us to reverse the order of the factors in Corollary~\ref{cor:solvable-exponential-decomposition}.
Let $g_0\in\mathsf{Hol}_x$.
By Corollary~\ref{cor:solvable-exponential-decomposition}\,(iii), write $g_0=\exp(\widehat{U})\cdot\exp(A)$ with $\widehat{U}\in\mathfrak{n}$ and $A\in\mathfrak{a}_x$.
Consider the curve
$$
    \widetilde{\sigma}(t):=\exp(\widehat{U})\cdot\exp((1-t)A)\cdot x,
    \qquad t\in[0,1].
$$
By Lemma~\ref{lem:exponential-horizontal-curve} and Remark~\ref{rem:n-horizontal-translates}, this curve is $\mathsf{N}$-horizontal and joins $g_0\cdot x$ to $\exp(\widehat{U})\cdot x$.
Since $g_0\in\mathsf{Hol}_x$, there is a $\mathsf{G}$-horizontal curve $\sigma$ joining $x$ to $g_0\cdot x$.
Concatenating $\sigma$ and $\widetilde{\sigma}$ gives an $\mathsf{N}$-horizontal curve joining $x$ to $\exp(\widehat{U})\cdot x$.
Therefore, $\exp(\widehat{U})\in\widehat{\mathsf{Hol}}_x\subseteq Z(\mathsf{N})$ by Lemma~\ref{lem:central-n-holonomy}.

Since $[\mathfrak{a}_x,\mathfrak{n}]\subseteq\mathfrak{n}$, as $\mathfrak{n}$ is an ideal, we have, with respect to the decomposition $\mathfrak{g}=\mathfrak{a}_x\oplus\mathfrak{n}$,
$$
    \mathrm{ad}_{\widehat{U}}=\begin{pmatrix}0&0\\ \mathrm{ad}_{\widehat{U}}|_{\mathfrak{a}_x}&\mathrm{ad}_{\widehat{U}}|_{\mathfrak{n}}\end{pmatrix}
    \quad\Longrightarrow\quad e^{\mathrm{ad}_{\widehat{U}}}
    =\begin{pmatrix}\mathrm{Id}_{\mathfrak{a}_x}&0\\ *&e^{\mathrm{ad}_{\widehat{U}}|_{\mathfrak{n}}}\end{pmatrix}.
$$
Moreover, since $\exp(\widehat{U})\in Z(\mathsf{N})$, we know that $\mathrm{Ad}(\exp(\widehat{U}))|_{\mathfrak{n}}=\mathrm{Id}_{\mathfrak{n}}$, so
$$
    \mathrm{Ad}(\exp(\widehat{U}))=e^{\mathrm{ad}_{\widehat{U}}}
    =\begin{pmatrix}\mathrm{Id}_{\mathfrak{a}_x}&0\\ *&\mathrm{Id}_{\mathfrak{n}}\end{pmatrix}.
$$
If we show that $*=0$, then $\mathrm{Ad}(\exp(\widehat{U}))=\mathrm{Id}_{\mathfrak{g}}$, and hence $\exp(\widehat{U})\in Z(\mathsf{G})=\{e\}$ because $\mathsf{G}$ is connected and centerless by Theorem~\ref{thm:solvable-algebraic-rigidity}\,(iii).
We would then conclude that $g_0\in\exp(\mathfrak{a}_x)$, as desired.
We will now show that $*=0$.

By Remark~\ref{rem:horizontal-complement-equivariance}\,(ii), we have
$$
    \mathrm{ad}_A=\begin{pmatrix}0&0\\ 0&\mathrm{ad}_A|_{\mathfrak{n}}\end{pmatrix}
    \quad\Longrightarrow\quad e^{\mathrm{ad}_A}
    =\begin{pmatrix}\mathrm{Id}&0\\ 0&e^{\mathrm{ad}_A|_{\mathfrak{n}}}\end{pmatrix}.
$$
Therefore,
$$
    \mathrm{Ad}(g_0)=e^{\mathrm{ad}_{\widehat{U}}}\cdot e^{\mathrm{ad}_A}
    =\begin{pmatrix}\mathrm{Id}&0\\ *&\mathrm{Id}\end{pmatrix}\cdot\begin{pmatrix}\mathrm{Id}&0\\ 0&e^{\mathrm{ad}_A|_{\mathfrak{n}}}\end{pmatrix}
    =\begin{pmatrix}\mathrm{Id}&0\\ *&e^{\mathrm{ad}_A|_{\mathfrak{n}}}\end{pmatrix},
$$
where the lower-left block is unchanged.
Since $g_0\in\mathsf{Hol}_x$, equivariance of $x\mapsto\mathfrak{a}_x$ and its constancy along $\mathsf{G}$-horizontal curves give $\mathrm{Ad}(g_0)(\mathfrak{a}_x)=\mathfrak{a}_{g_0\cdot x}=\mathfrak{a}_x$.
Thus $*=0$, as desired.
This proves that $\mathsf{Hol}_x\subseteq\exp(\mathfrak{a}_x)$.

We now show that $\widehat{\mathsf{Hol}}_x=\{e\}$.
Let $n_0\in\widehat{\mathsf{Hol}}_x$ and choose a piecewise smooth $\mathsf{N}$-horizontal curve $\sigma:[0,1]\to M$ joining $x$ to $n_0\cdot x$.
By Lemma~\ref{lem:horizontal-decomposition}, we may write $\sigma(t)=g(t)\cdot\overline{\sigma}(t)$, where $g(0)=e$ and $\overline{\sigma}$ is $\mathsf{G}$-horizontal.

We claim that $g(t)\in\exp(\mathfrak{a}_x)$ for every $t$.
By Corollary~\ref{cor:solvable-exponential-decomposition}\,(iii), write $g(t)=n(t)\cdot a(t)$, with $n(t)\in\mathsf{N}$, $a(t)\in\exp(\mathfrak{a}_x)$ and $n(0)=a(0)=e$.
The factors are piecewise smooth, and Corollary~\ref{cor:solvable-exponential-decomposition}\,(ii) allows us to write $a(t)=\exp(A(t))$ for a piecewise smooth curve $A(t)\in\mathfrak{a}_x$.
By Lemma~\ref{lem:exponential-horizontal-curve}, $t\mapsto a(t)\cdot\overline{\sigma}(t)$ is $\mathsf{N}$-horizontal.
This curve and $\sigma$ start at $x$ and project to the same curve in $M/\mathsf{N}$.
Uniqueness of horizontal lifts (Proposition~\ref{prop:horizontal-lifting}) therefore gives $\sigma(t)=a(t)\cdot\overline{\sigma}(t)$.
Freeness of the action implies $n(t)=e$, proving the claim.

In particular, $\overline{\sigma}$ is a $\mathsf{G}$-horizontal curve joining $x$ to $a(1)^{-1}\cdot n_0\cdot x$, so $a(1)^{-1}n_0\in\mathsf{Hol}_x\subseteq\exp(\mathfrak{a}_x)$.
Since $a(1)\in\exp(\mathfrak{a}_x)$, it follows that $n_0\in\mathsf{N}\cap\exp(\mathfrak{a}_x)=\{e\}$ by Corollary~\ref{cor:solvable-exponential-decomposition}\,(iii).
\end{proof}

\section{The $\mathsf{N}$-orbits are pairwise isometric by affine conjugations}\label{sec:affine-conjugation-orbits}
Throughout this section, we work under Assumption~\ref{ass:main}.
We use the notation and formulas from Appendix~\ref{app:inner-products}.
The goal of this section is to prove Theorem~\ref{thm:a0-orbit-image} below.

Let $\widehat{\Sigma}\subseteq M$ be an $\mathsf{N}$-horizontal leaf, and set $\mathfrak{a}_0:=\mathfrak{a}_x$ for any $x\in\widehat{\Sigma}$.
By Corollary~\ref{cor:horizontal-complement-subalgebra}, this does not depend on the choice of $x$.
In \cite[Corollary~9.8]{BL23}, the authors prove that the image of the map $\hat{h}|_{\widehat{\Sigma}}:\widehat{\Sigma}\rightarrow\mathrm{Sym}^2_+(\mathfrak{n}^*)$ lies in a single $\mathsf{Aut}(\mathfrak{n})$-orbit.
In Corollary~\ref{cor:centralizer-orbit}, we obtained a refinement of this result: the image lies in a single $\mathsf{J}$-orbit, where
$$
    \mathsf{J}:=\{\varphi\in\mathsf{Aut}(\mathfrak{n}):[\varphi,\mathrm{ad}_A|_{\mathfrak{n}}]=0\text{ for every }A\in\mathfrak{a}_0\}.
$$
We now refine this result further.
Set $\mathsf{A}_0:=\exp(\mathfrak{a}_0)\subseteq\mathsf{G}$ and notice that $\mathrm{Ad}(\mathsf{A}_0)|_{\mathfrak{n}}\subseteq\mathsf{J}$.

\begin{theorem}\label{thm:a0-orbit-image}
Under Assumption~\ref{ass:main}, the image of the map $\hat{h}|_{\widehat{\Sigma}}:\widehat{\Sigma}\rightarrow\mathrm{Sym}^2_+(\mathfrak{n}^*)$ lies in a single $\mathrm{Ad}(\mathsf{A}_0)|_{\mathfrak{n}}$-orbit.
\end{theorem}

By Remark~\ref{rem:affine-orbit-isometries}, this means that every pair of $\mathsf{N}$-orbits in $M$ is $c_a$-equivariantly isometric for some $a\in\mathsf{A}_0$, where $c_a(n):=ana^{-1}$.
We call an isometry between $\mathsf{N}$-orbits \emph{affine} if it is $\varphi$-equivariant for some $\varphi\in\mathsf{Aut}(\mathsf{N})$, and an \emph{affine conjugation} if $\varphi=c_g$ for some $g\in\mathsf{G}$.
This justifies the title of this section.

Notice that $\mathsf{A}_0\subseteq\mathsf{G}$ is a closed, connected, simply connected abelian subgroup, by Corollary~\ref{cor:solvable-exponential-decomposition}\,(ii).
By Corollary~\ref{cor:a0-action-on-section}, its action on $M$ restricts to an action on $\widehat{\Sigma}$ with compact orbit space $\widehat{\Sigma}/\mathsf{A}_0\equiv M/\mathsf{G}$.
It also acts isometrically on $\mathrm{Sym}^2_+(\mathfrak{n}^*)$ through the adjoint representation $\mathrm{Ad}(\mathsf{A}_0)|_{\mathfrak{n}}\subseteq\mathsf{GL}(\mathfrak{n})$.
For these actions, the map $\hat{h}|_{\widehat{\Sigma}}$ is $\mathsf{A}_0$-equivariant, as follows by restricting the identity in Lemma~\ref{lem:orbit-metric-equivariance} to $\mathfrak{n}$.

We divide the argument into two subsections.
In the first, we show that the $\mathsf{A}_0$-action on $\mathrm{Sym}^2_+(\mathfrak{n}^*)$ is free, proper and isometric, and that every orbit meeting the image of $\hat{h}|_{\widehat{\Sigma}}$ is totally geodesic.
In the second, we compose the $\mathsf{A}_0$-equivariant harmonic map $\hat{h}|_{\widehat{\Sigma}}$  (see Theorem~\ref{thm:nilradical-metric-identities}\,(ii)) with the convex squared-distance function to one such orbit to construct an $\mathsf{A}_0$-invariant subharmonic function on $\widehat{\Sigma}$. Compactness of $\widehat{\Sigma}/\mathsf{A}_0$ and the strong maximum principle force this function to be constant and hence zero, proving Theorem~\ref{thm:a0-orbit-image}.

\subsection{The $\mathsf{A}_0$-action on $\mathrm{Sym}^2_+(\mathfrak{n}^*)$}\label{subsec:a0-inner-product-action}

\textcolor{white}{aaa}

\begin{lemma}\label{lem:a0-action-free-proper}
The action of $\mathsf{A}_0$ on $\mathrm{Sym}^2_+(\mathfrak{n}^*)$ is free, proper and isometric.
\end{lemma}

\begin{proof}
The action is isometric by the $\mathsf{GL}(\mathfrak{n})$-invariance of the metric introduced in Appendix~\ref{app:inner-products}.
Fix $x_0\in\widehat{\Sigma}$.
By Theorem~\ref{thm:solvable-algebraic-rigidity}\,(iv), for every $a\in\mathsf{A}_0$ we have
$$
    \mathrm{Ad}(a)|_{\mathfrak{n}}\in\mathsf{O}(\mathfrak{n},\hat{h}_{x_0})
    \quad\Longleftrightarrow\quad a=e,
$$
since $\mathsf{A}_0\cap Z(\mathsf{N})=\{e\}$ by Corollary~\ref{cor:solvable-exponential-decomposition}\,(iii).
Thus the stabilizer of $\hat{h}_{x_0}$ is trivial.

We next show that the orbit $\mathsf{A}_0\cdot\hat{h}_{x_0}$ is closed.
For every $A\in\mathfrak{a}_0$, Corollary~\ref{cor:normal-adjoint-derivations} gives
$$
    \exp(A)\cdot\hat{h}_{x_0}
    =e^{\mathrm{ad}_A|_{\mathfrak{n}}}\cdot\hat{h}_{x_0}
    =e^{\mathrm{S}_{x_0}(\mathrm{ad}_A|_{\mathfrak{n}})}\cdot e^{\mathrm{sk}_{x_0}(\mathrm{ad}_A|_{\mathfrak{n}})}\cdot\hat{h}_{x_0}
    =e^{\mathrm{S}_{x_0}(\mathrm{ad}_A|_{\mathfrak{n}})}\cdot\hat{h}_{x_0}.
$$
Therefore, defining the linear subspace
$$
    D:=\{\mathrm{S}_{x_0}(\mathrm{ad}_A|_{\mathfrak{n}}):A\in\mathfrak{a}_0\}\subseteq\mathrm{End}(\mathfrak{n}),
$$
we have
$$
    \mathsf{A}_0\cdot\hat{h}_{x_0}=\exp(D)\cdot\hat{h}_{x_0}.
$$
By Corollary~\ref{cor:normal-adjoint-derivations}\,(iii), $D$ consists of commuting $\hat{h}_{x_0}$-symmetric endomorphisms.
Thus there is a $\hat{h}_{x_0}$-orthonormal basis of $\mathfrak{n}$ in which all elements of $D$ are diagonal.
In this basis, $\exp(D)$ is a connected Lie subgroup of the group of positive diagonal matrices.
Since the latter is simply connected and abelian, $\exp(D)$ is closed in it.
The positive diagonal group is closed in $\mathsf{GL}(\mathfrak{n})$, so $\exp(D)$ is also closed in $\mathsf{GL}(\mathfrak{n})$.
The action of $\mathsf{GL}(\mathfrak{n})$ on $\mathrm{Sym}^2_+(\mathfrak{n}^*)$ is proper (see Appendix~\ref{app:inner-products}), so $\exp(D)\cdot\hat{h}_{x_0}$ is a closed embedded submanifold.

The orbit $\mathsf{A}_0\cdot\hat{h}_{x_0}$ is therefore closed, and the stabilizer of $\hat{h}_{x_0}$ is trivial.
By Lemma~\ref{lem:closed-orbit-compact-isotropy-proper-action}, the $\mathsf{A}_0$-action is proper.
Since $\mathsf{A}_0$ is connected, simply connected and abelian, Lemma~\ref{lem:simply-connected-solvable-free-action} implies that this action is free.
\end{proof}

We now show that the $\mathsf{A}_0$-orbits meeting the image of $\hat{h}|_{\widehat{\Sigma}}$ are totally geodesic.

\begin{lemma}\label{lem:a0-orbits-totally-geodesic}
For every $x\in\widehat{\Sigma}$, the orbit $\mathsf{A}_0\cdot\hat{h}_x\subseteq\mathrm{Sym}^2_+(\mathfrak{n}^*)$ is a totally geodesic submanifold.
\end{lemma}

\begin{proof}
Fix $x_0\in\widehat{\Sigma}$ and set $\bar{h}:=\hat{h}_{x_0}$.
Let $X_h=\rho(E_h)h$ be an $\mathsf{A}_0$-basic vector field on $\mathrm{Sym}^2_+(\mathfrak{n}^*)$, represented by a smooth map $E :\mathrm{Sym}^2_+(\mathfrak{n}^*)\rightarrow \mathrm{End}(\mathfrak{n})$ such that $E_h$ is $h$-symmetric for every $h$.
Let $U_h=\rho(\mathrm{ad}_A|_{\mathfrak{n}})h$ be the $\mathsf{A}_0$-fundamental vector field associated with $A\in\mathfrak{a}_0$.
We denote the Levi-Civita connection of $\mathrm{Sym}^2_+(\mathfrak{n}^*)$ by $\nabla^{\mathrm{sym}}$.
To prove that the shape operator vanishes at $\bar{h}$, it suffices to show that $(\nabla^{\mathrm{sym}}_UX)_{\bar{h}}$ is $\mathsf{A}_0$-horizontal.
Since the action is isometric, the same will then hold along the whole orbit.

Since $X$ is $\mathsf{A}_0$-horizontal, we have
$$
    \mathrm{Tr}\bigl(E_h\cdot\mathrm{S}_h(\mathrm{ad}_{\tilde{A}}|_{\mathfrak{n}})\bigr)=0
    \qquad\text{for every }\tilde{A}\in\mathfrak{a}_0.
$$
Differentiating at $h=\bar{h}$ in the direction $U_{\bar{h}}$ gives
\begin{equation}\label{eq:a0-horizontal-orthogonality-derivative}
    \mathrm{Tr}\bigl(dE_{\bar{h}}(U_{\bar{h}})\cdot\mathrm{S}_{x_0}(\mathrm{ad}_{\tilde{A}}|_{\mathfrak{n}})\bigr)
    +\mathrm{Tr}\bigl(E_{\bar{h}}\cdot d(\mathrm{S}_*(\mathrm{ad}_{\tilde{A}}|_{\mathfrak{n}}))_{\bar{h}}(U_{\bar{h}})\bigr)=0.
\end{equation}
By Remark~\ref{rem:inner-product-fixed-endomorphism} and Corollary~\ref{cor:normal-adjoint-derivations}\,(iii),
$$
    d(\mathrm{S}_*(\mathrm{ad}_{\tilde{A}}|_{\mathfrak{n}}))_{\bar{h}}(U_{\bar{h}})
    =[\mathrm{S}_{x_0}(\mathrm{ad}_A|_{\mathfrak{n}}),(\mathrm{ad}_{\tilde{A}}|_{\mathfrak{n}})^{t_{x_0}}]
    =0.
$$
Substituting into \eqref{eq:a0-horizontal-orthogonality-derivative}, and using that the trace pairing of a symmetric and a skew-symmetric endomorphism vanishes, we obtain
$$
    \mathrm{Tr}\bigl(\mathrm{S}_{x_0}(dE_{\bar{h}}(U_{\bar{h}}))\cdot\mathrm{S}_{x_0}(\mathrm{ad}_{\tilde{A}}|_{\mathfrak{n}})\bigr)=0
    \qquad\text{for every }\tilde{A}\in\mathfrak{a}_0.
$$
By Proposition~\ref{lem:inner-product-connection}, this means precisely that $(\nabla^{\mathrm{sym}}_UX)_{\bar{h}}$ is $\mathsf{A}_0$-horizontal, as desired.
\end{proof}

\subsection{Proof of Theorem~\ref{thm:a0-orbit-image}}\label{subsec:proof-a0-orbit-image}

\textcolor{white}{aaa}

We will use the following result.

\begin{prop}[{\cite[Theorem~4.1(1)]{BishopONeill69}}]\label{prop:convex-squared-distance}
Let $M$ be a complete simply connected Riemannian manifold with non-positive sectional curvature, and let $S\subseteq M$ be a nonempty, closed, connected, totally geodesic submanifold.
Then the function $\varphi:M\rightarrow\mathbb{R}$ defined by $\varphi(x):=d(x,S)^2$ is smooth and convex.
\end{prop}

We also record the following elementary fact.

\begin{lemma}\label{lem:harmonic-convex-subharmonic}
Let $f:M_1\rightarrow M_2$ be a harmonic map between Riemannian manifolds, and let $\varphi:M_2\rightarrow\mathbb{R}$ be a smooth convex function.
Then $\varphi\circ f$ is subharmonic.
\end{lemma}

\begin{proof}
Let $\{\xi_1,\ldots,\xi_k\}$ be a local orthonormal frame on $M_1$.
The chain rule for covariant derivatives gives
$$
    \mathrm{Hess}(\varphi\circ f)(\xi_i,\xi_i)
    =\mathrm{Hess}\,\varphi(df(\xi_i),df(\xi_i))
    +d\varphi\bigl((\nabla df)(\xi_i,\xi_i)\bigr).
$$
Taking the trace, we obtain
\begin{equation}\label{eq:harmonic-convex-composition}
    \Delta(\varphi\circ f)
    =\sum_i\mathrm{Hess}\,\varphi(df(\xi_i),df(\xi_i))
    +d\varphi(\tau(f)),
\end{equation}
where $\tau(f):=\sum_i(\nabla df)(\xi_i,\xi_i)$ is the tension field of $f$.
Since $f$ is harmonic, $\tau(f)=0$.
Since $\varphi$ is convex, its Hessian is positive semi-definite.
Thus \eqref{eq:harmonic-convex-composition} gives $\Delta(\varphi\circ f)\geq0$.
\end{proof}

We can now prove the theorem.

\begin{proof}[Proof of Theorem~\ref{thm:a0-orbit-image}]
Fix $x_0\in\widehat{\Sigma}$ and consider the orbit
$$
    O:=\mathrm{Ad}(\mathsf{A}_0)|_{\mathfrak{n}}\cdot\hat{h}_{x_0}
    =\mathsf{A}_0\cdot\hat{h}_{x_0}
    \subseteq\mathrm{Sym}^2_+(\mathfrak{n}^*).
$$
By Lemmas~\ref{lem:a0-action-free-proper} and \ref{lem:a0-orbits-totally-geodesic}, $O$ is a closed, connected, totally geodesic submanifold.
The symmetric space $\mathrm{Sym}^2_+(\mathfrak{n}^*)$ is complete, simply connected and has non-positive sectional curvature (see \cite[Section~3.1]{Heber98}).
Therefore, by Proposition~\ref{prop:convex-squared-distance}, the function
$$
    \varphi:\mathrm{Sym}^2_+(\mathfrak{n}^*)\rightarrow\mathbb{R},
    \qquad\varphi(h):=d(h,O)^2,
$$
is smooth and convex.
Moreover, $\varphi$ is $\mathsf{A}_0$-invariant, since the action is isometric and preserves $O$.

The map $\hat{h}|_{\widehat{\Sigma}}$ is $\mathsf{A}_0$-equivariant by Lemma~\ref{lem:orbit-metric-equivariance}, and harmonic by Theorem~\ref{thm:nilradical-metric-identities}\,(ii).
Hence the composition $\varphi\circ\hat{h}:\widehat{\Sigma}\rightarrow\mathbb{R}$ is $\mathsf{A}_0$-invariant and, by Lemma~\ref{lem:harmonic-convex-subharmonic}, subharmonic.
Its invariance implies that it descends to a continuous function on $\widehat{\Sigma}/\mathsf{A}_0$.
This quotient is compact by Corollary~\ref{cor:a0-action-on-section}, so the descended function attains a maximum.
Consequently, $\varphi\circ\hat{h}$ attains a maximum on $\widehat{\Sigma}$.

Since $\widehat{\Sigma}$ is connected, the strong maximum principle implies that $\varphi\circ\hat{h}$ is constant (see \cite[Theorem~3.5]{GilbargTrudinger01}).
But $(\varphi\circ\hat{h})(x_0)=0$, so $\varphi\circ\hat{h}\equiv0$.
Thus
$$
    d(\hat{h}_x,O)^2=0
    \qquad\text{for every }x\in\widehat{\Sigma}.
$$
Since $O$ is closed, we conclude that 
$$
    \hat{h}_x\in O=\mathrm{Ad}(\mathsf{A}_0)|_{\mathfrak{n}}\cdot\hat{h}_{x_0}
    \qquad\text{for every }x\in\widehat{\Sigma},
$$
as desired.
\end{proof}

\section{Proof of Theorem A: The solvable case}\label{sec:solvable-splitting}

Combining the results of the previous sections, we prove the following special case of Theorem~\ref{thm:main-splitting}.

\begin{theorem}\label{thm:solvable-splitting}
Let $M$ be a complete connected negative-Einstein manifold, and let $\mathsf{G}$ be a connected solvable Lie group acting freely, properly, isometrically and cocompactly on $M$.
Then every orbit is equivariantly isometric to a fixed $\mathsf{G}$-homogeneous manifold $H$, and $M$ is equivariantly isometric to $H\times M/\mathsf{G}$.
\end{theorem}

The main step is to prove that $M/\mathsf{N}$ has non-positive Ricci curvature.
We obtain this from a maximum principle applied to harmonic one-forms on the $\mathsf{N}$-horizontal leaves.
The resulting estimate forces the $\mathsf{G}$-action to be polar and its orbits to be totally geodesic.
Finally, Theorem~\ref{thm:trivial-n-holonomy} and the freeness of the $\mathsf{A}_0$-action on $\mathrm{Sym}^2_+(\mathfrak{n}^*)$ imply that the $\mathsf{G}$-holonomy is trivial, yielding the global splitting.

\subsection{Horizontal vector fields associated with vertical tensors}\label{subsec:fields-vertical-tensors}

\textcolor{white}{aaa}

In this subsection, we associate horizontal vector fields with vertical endomorphisms by taking the pointwise adjoint of the shape operator.

\begin{definition}\label{def:shape-adjoint-field}
For every vertical $(1,1)$-tensor $E$, we denote by $L^*E$ the horizontal vector field characterized by
$$
    \langle L^*E,X\rangle=\mathrm{Tr}(E\cdot L_X)
    \qquad\text{for every horizontal vector field }X.
$$
Similarly, for every $\mathsf{N}$-vertical $(1,1)$-tensor $\widehat{E}$, we denote by $\hat{L}^*\widehat{E}$ the corresponding $\mathsf{N}$-horizontal vector field.
\end{definition}

In fact, $L^*$ is the pointwise adjoint of the map $\mathcal{H}\rightarrow\mathrm{End}\,\mathcal{V}$, $X\mapsto L_X$, where $\mathcal{H}$ carries the Riemannian metric and $\mathrm{End}\,\mathcal{V}$ carries the induced inner product $\langle E_1,E_2\rangle=\mathrm{Tr}(E_1E_2^t)$.
Since $L_X$ is self-adjoint, this gives precisely the formula in Definition~\ref{def:shape-adjoint-field}.
The same interpretation applies to $\hat{L}^*$.

Since $\mathsf{N}$ is normal in $\mathsf{G}$, the $\mathsf{G}$-action preserves both $\widehat{\mathcal{V}}$ and $\widehat{\mathcal{H}}$.
Thus, it acts on the bundles $\mathrm{End}\,\widehat{\mathcal{V}}$ and $\mathrm{End}\,\widehat{\mathcal{H}}$, and it makes sense to speak of $\mathsf{G}$-equivariant sections of these bundles.

\begin{remark}\label{rem:shape-adjoint-invariance}
The following properties hold:
\begin{enumerate}
\renewcommand{\labelenumi}{(\roman{enumi})}
\item If $E$ is a $\mathsf{G}$-equivariant vertical $(1,1)$-tensor, then $L^*E$ is basic.
\item If $\widehat{E}$ is a $\mathsf{G}$-equivariant $\mathsf{N}$-vertical $(1,1)$-tensor, then $\hat{L}^*\widehat{E}$ is $\mathsf{N}$-horizontal and $\mathsf{G}$-invariant.
\end{enumerate}
Indeed, the $\mathsf{G}$-action preserves the metric, the Levi-Civita connection and the corresponding distributions.
Consequently, for every $g\in\mathsf{G}$,
\begin{equation}\label{eq:shape-operator-equivariance}
    d(m_g)(L_XU)
    =L_{d(m_g)X}\bigl(d(m_g)U\bigr),
\end{equation}
where $X$ is horizontal and $U$ is vertical.
If $E$ is equivariant, invariance of the metric and the trace under conjugation shows that $L^*E$ is invariant, and hence basic.
The same argument with $\hat{L}$ and $\widehat{E}$ proves item~(ii).
\end{remark}

In Section~\ref{sec:shape-commutativity}, we studied the vector field $\mathrm{W}_A=L^*[\mathrm{ad}(A),\mathrm{ad}(A)^t]$ when $A$ is invariant, horizontally Killing and $\mathcal{A}$-valued.
The invariance of $A$ ensures that $[\mathrm{ad}(A),\mathrm{ad}(A)^t]$ is equivariant, so $\mathrm{W}_A$ is basic by Remark~\ref{rem:shape-adjoint-invariance}.

The section $\hat{L}^*\hat{L}$ was defined in Theorem~\ref{thm:nilradical-metric-identities}\,(iii).
By Definition~\ref{def:shape-adjoint-field}, it satisfies $\hat{L}^*\hat{L}(\widehat{X})=\hat{L}^*(\hat{L}_{\widehat{X}})$.
If $\widehat{X}$ is $\mathsf{G}$-invariant, then $\hat{L}_{\widehat{X}}$ is $\mathsf{G}$-equivariant by \eqref{eq:shape-operator-equivariance}, since $d(m_g)\widehat{X}=\widehat{X}$ for every $g\in\mathsf{G}$.
Remark~\ref{rem:shape-adjoint-invariance} then implies that $\hat{L}^*\hat{L}(\widehat{X})$ is invariant.
Since the action is free and proper, we obtain the following.

\begin{remark}\label{obs:shape-square-equivariance}
The section $\hat{L}^*\hat{L}$ of $\mathrm{End}\,\widehat{\mathcal{H}}\rightarrow M$ is $\mathsf{G}$-equivariant, pointwise self-adjoint and positive semidefinite.
\end{remark}

\subsection{The quotient $M/\mathsf{N}$ has non-positive Ricci curvature}\label{subsec:nilradical-quotient-ricci}

\textcolor{white}{aaa}

Our goal in this subsection is to prove that $M/\mathsf{N}$ has non-positive Ricci curvature.
By Theorem~\ref{thm:nilradical-metric-identities}\,(iii), this is equivalent to every eigenvalue of $\hat{L}^*\hat{L}$ being at most $1$.
We will express the shape operators in terms of $\mathcal{A}$-valued fields and use invariant horizontally Killing fields to construct harmonic one-forms $\omega_A$ on the $\mathsf{N}$-horizontal leaves.
Choosing one of these forms suitably will allow us to apply the Bochner formula at a maximum of its norm and obtain the desired eigenvalue bound.

We first use Theorem~\ref{thm:a0-orbit-image} to express every shape operator $\hat{L}_{\widehat{X}}$ as one associated with an $\mathcal{A}$-valued field.

\begin{prop}\label{prop:shape-alpha}
There exists a unique smooth $\mathsf{G}$-equivariant bundle map $\alpha:\widehat{\mathcal{H}}\rightarrow\mathcal{A}$ characterized by
$$
    \hat{L}_{\widehat{X}}=\hat{L}_{\alpha(\widehat{X})}
    \qquad\text{for every }\widehat{X}\in\widehat{\mathcal{H}}.
$$
\end{prop}

\begin{proof}
Fix a connected $\mathsf{N}$-horizontal leaf $\widehat{\Sigma}$ and write $\mathfrak{a}_0:=\mathfrak{a}_x$ for $x\in\widehat{\Sigma}$.
By Theorem~\ref{thm:a0-orbit-image}, the image of $\hat{h}|_{\widehat{\Sigma}}$ lies in an $\mathrm{Ad}(\mathsf{A}_0)|_{\mathfrak{n}}$-orbit.
Thus, $d\hat{h}_x(\widehat{X}_x)$ is tangent to this orbit for every $\widehat{X}_x\in\widehat{\mathcal{H}}_x$.
By Lemma~\ref{lem:orbit-metric-differential} applied to the $\mathsf{N}$-action, there exists $A_x\in\mathcal{A}_x$ such that
$$
    -\rho({}^{\mathfrak{n}}\hat{L}_{\widehat{X}_x})\hat{h}_x
    =d\hat{h}_x(\widehat{X}_x)
    =\rho\bigl(\mathrm{ad}_{{}^{\raisebox{0.15em}{$\scriptstyle\mathfrak{g}$}}\! A_x}|_{\mathfrak{n}}\bigr)\hat{h}_x.
$$
The kernel of the map $E\mapsto\rho(E)\hat{h}_x$ consists of the $\hat{h}_x$-skew-symmetric endomorphisms.
Since ${}^{\mathfrak{n}}\hat{L}_{\widehat{X}_x}$ is self-adjoint, taking symmetric parts gives
$$
    {}^{\mathfrak{n}}\hat{L}_{\widehat{X}_x}
    =-\mathrm{S}_x\bigl(\mathrm{ad}_{{}^{\raisebox{0.15em}{$\scriptstyle\mathfrak{g}$}}\! A_x}|_{\mathfrak{n}}\bigr)
    ={}^{\mathfrak{n}}\hat{L}_{A_x}.
$$
Here the last equality is Lemma~\ref{lem:shape-symmetric-adjoint}.

The map $A_x\mapsto\hat{L}_{A_x}$ is injective.
Indeed, if $\hat{L}_{A_x}=0$, then Lemma~\ref{lem:shape-symmetric-adjoint} implies that $\mathrm{ad}_{{}^{\raisebox{0.15em}{$\scriptstyle\mathfrak{g}$}}\! A_x}|_{\mathfrak{n}}$ is skew-symmetric.
Since $\mathrm{ad}_{{}^{\raisebox{0.15em}{$\scriptstyle\mathfrak{g}$}}\! A_x}$ vanishes on $\mathfrak{a}_x$, it is $h_x$-skew-symmetric on $\mathfrak{g}$.
Theorem~\ref{thm:solvable-algebraic-rigidity}\,(i) gives $A_x=0$.

Consequently, the vector $A_x$ associated with $\widehat{X}_x$ is unique and depends linearly on $\widehat{X}_x$.
This defines $\alpha_x$.
The smooth bundle map $\mathcal{A}\rightarrow\mathrm{End}\,\widehat{\mathcal{V}}$, $A\mapsto\hat{L}_A$, is injective and has constant rank, so its inverse on its image is smooth.
Hence $\alpha$ is smooth.
To prove equivariance, fix $g\in\mathsf{G}$ and $\widehat{X}_x\in\widehat{\mathcal{H}}_x$.
Since the action preserves $\mathcal{A}$ and the shape operators, the vector $d(m_g)_x(\alpha_x(\widehat{X}_x))$ belongs to $\mathcal{A}_{g\cdot x}$ and has the same shape operator as $d(m_g)_x(\widehat{X}_x)$.
By uniqueness in the definition of $\alpha_{g\cdot x}$, we obtain
$$
    \alpha_{g\cdot x}\bigl(d(m_g)_x\widehat{X}_x\bigr)
    =d(m_g)_x(\alpha_x(\widehat{X}_x)).
$$
\end{proof}

We now use these shape operators to construct the one-forms needed in the maximum principle argument.
For every $\mathcal{A}$-valued vector field $A$, consider the one-form $\omega_A:\widehat{\mathcal{H}}\rightarrow\mathbb{R}$ defined by
$$
    \omega_A(\widehat{X})
    :=\mathrm{Tr}(\hat{L}_A\cdot\hat{L}_{\widehat{X}})
    =\langle\hat{L}^*\hat{L}(A),\widehat{X}\rangle.
$$
We regard $\omega_A$ as a one-form on each $\mathsf{N}$-horizontal leaf by restriction to its tangent bundle.

\begin{prop}\label{prop:harmonic-shape-one-form}
If $A$ is invariant, horizontally Killing and $\mathcal{A}$-valued, then $\omega_A$ is closed and coclosed, and hence harmonic, along every $\mathsf{N}$-horizontal leaf $\widehat{\Sigma}$.
\end{prop}

\begin{proof}
By Corollary~\ref{cor:parallel-a-shape}, we have $\nabla_{\widehat{X}}(\hat{L}_A)=0$ for every $\mathsf{N}$-horizontal vector field $\widehat{X}$.
Thus, for $\mathsf{N}$-horizontal vector fields $\widehat{X},\widehat{Y}$,
$$
    \widehat{X}\bigl(\omega_A(\widehat{Y})\bigr)
    =\mathrm{Tr}\bigl(\hat{L}_A\cdot\nabla_{\widehat{X}}(\hat{L}_{\widehat{Y}})\bigr).
$$
Since $\widehat{\mathrm{A}}=0$ by Theorem~\ref{thm:nilradical-polarity}, covariant differentiation in $\mathsf{N}$-horizontal directions preserves both $\widehat{\mathcal{V}}$ and $\widehat{\mathcal{H}}$.
In particular, no $\mathsf{N}$-vertical projection is needed in the preceding formula, and the induced connection on $\widehat{\mathcal{H}}|_{\widehat{\Sigma}}$ is the Levi-Civita connection of the leaf.
Consequently,
$$
    (\nabla^{\widehat{\Sigma}}_{\widehat{X}}\omega_A)(\widehat{Y})
    =\mathrm{Tr}\bigl(\hat{L}_A\cdot(\nabla_{\widehat{X}}\hat{L})_{\widehat{Y}}\bigr).
$$
It follows that
$$
    d\omega_A(\widehat{X},\widehat{Y})
    =\mathrm{Tr}\bigl(\hat{L}_A\cdot\bigl((\nabla_{\widehat{X}}\hat{L})_{\widehat{Y}}-(\nabla_{\widehat{Y}}\hat{L})_{\widehat{X}}\bigr)\bigr)
    =0,
$$
where the last equality follows from Proposition~\ref{prop:polar-curvature-shape} applied to the $\mathsf{N}$-action.
Thus, $\omega_A$ is closed.

To prove that it is coclosed, let $\{\widehat{X}_i\}$ be a local orthonormal frame on $\widehat{\Sigma}$.
Taking the trace of the formula for $\nabla^{\widehat{\Sigma}}\omega_A$ gives
$$
    d^*\omega_A
    =-\sum_i\mathrm{Tr}\bigl(\hat{L}_A\cdot(\nabla_{\widehat{X}_i}\hat{L})_{\widehat{X}_i}\bigr)
    =-\mathrm{Tr}(\hat{L}_A\cdot\hat{\tau})=0,
$$
since $\hat{\tau}=0$ by Theorem~\ref{thm:nilradical-metric-identities}\,(i).
\end{proof}

We can now combine the harmonicity of these one-forms with cocompactness to control the largest eigenvalue of $\hat{L}^*\hat{L}$.

\begin{theorem}\label{thm:orbit-ricci-lower-bound}
The quotient $M/\mathsf{N}$ has non-positive Ricci curvature.
Equivalently,
$$
    \|\hat{L}_{\widehat{X}}\|\leq\|\widehat{X}\|
    \qquad\text{for every }\widehat{X}\in\widehat{\mathcal{H}}.
$$
\end{theorem}

\begin{proof}
By Theorem~\ref{thm:nilradical-metric-identities}\,(iii), the equivalence follows, which in turn is equivalent to every eigenvalue of $\hat{L}^*\hat{L}$ being at most $1$. Denote $\lambda_{\max}$ the maximum eigenvalue of an endomorphism. 
Since $\hat{L}^*\hat{L}$ is $\mathsf{G}$-equivariant (Remark~\ref{obs:shape-square-equivariance}), the map $x\mapsto\lambda_{\max}(\hat{L}^*\hat{L}_x)$ defines a continuous $\mathsf{G}$-invariant function on $M$.
The quotient $M/\mathsf{G}$ is compact, so this function attains a maximum $\lambda_0$ at some $x_0\in M$.
Our goal is to prove that $\lambda_0\leq1$.

Proposition~\ref{prop:harmonic-shape-one-form} allows us to apply the Bochner formula to $\|\omega_A\|^2$ whenever $A$ is invariant, horizontally Killing and $\mathcal{A}$-valued.
To extract information about $\lambda_0$, we choose $A$ so that $\|\omega_A\|^2$ also attains its maximum at $x_0$.

If $\lambda_0=0$, there is nothing to prove.
Suppose that $\lambda_0>0$, and choose a unit eigenvector $\widehat{X}_{x_0}\in\widehat{\mathcal{H}}_{x_0}$ of $\hat{L}^*\hat{L}_{x_0}$ with eigenvalue $\lambda_0$.
Set
$$
    A_{x_0}:=\frac{\alpha(\widehat{X}_{x_0})}{\sqrt{\lambda_0}}\in\mathcal{A}_{x_0}.
$$
By Proposition~\ref{prop:shape-alpha}, we have
\begin{equation}\label{eq:shape-one-form-at-maximum}
    \hat{L}^*\hat{L}_{x_0}(A_{x_0})
    =\frac{1}{\sqrt{\lambda_0}}\hat{L}^*_{x_0}(\hat{L}_{\widehat{X}_{x_0}})
    =\sqrt{\lambda_0}\,\widehat{X}_{x_0}.
\end{equation}
Moreover,
$$
    \|\hat{L}_{A_{x_0}}\|^2
    =\frac{1}{\lambda_0}\|\hat{L}_{\widehat{X}_{x_0}}\|^2
    =\frac{1}{\lambda_0}\langle\hat{L}^*\hat{L}_{x_0}(\widehat{X}_{x_0}),\widehat{X}_{x_0}\rangle=1.
$$
Extend $A_{x_0}$ to an invariant horizontally Killing $\mathcal{A}$-valued vector field $A$ using Proposition~\ref{prop:invariant-a-extension}.
Corollary~\ref{cor:parallel-a-shape} implies that $\|\hat{L}_A\|$ is constant along $\mathsf{N}$-horizontal curves.
It is also $\mathsf{G}$-invariant, so connectedness of $M$ gives $\|\hat{L}_A\|\equiv1$.

Since $\hat{L}^*\hat{L}(A)$ is the metric dual of $\omega_A$ on each $\mathsf{N}$-horizontal leaf, for every $x\in M$ we have
$$
    \|\omega_A\|^2(x)=\|\hat{L}^*\hat{L}_x(A_x)\|^2
    \leq\|\hat{L}^*_x\|_{\mathrm{op}}^2\cdot\|\hat{L}_{A_x}\|^2
    =\lambda_{\max}(\hat{L}^*\hat{L}_x)\leq\lambda_0.
$$
Here $\|\cdot\|_{\mathrm{op}}$ denotes the operator norm, and we used the identity $\|\hat{L}^*_x\|_{\mathrm{op}}^2=\lambda_{\max}(\hat{L}^*\hat{L}_x)$.
On the other hand, equation~\eqref{eq:shape-one-form-at-maximum} gives
$$
    \|\omega_A\|^2(x_0)
    =\|\hat{L}^*\hat{L}_{x_0}(A_{x_0})\|^2
    =\lambda_0.
$$
Thus, $\|\omega_A\|^2$ attains its maximum at $x_0$, as required.

Let $\widehat{\Sigma}$ be the $\mathsf{N}$-horizontal leaf through $x_0$.
By Proposition~\ref{prop:harmonic-shape-one-form}, $\omega_A$ is harmonic on $\widehat{\Sigma}$.
The Bochner formula gives
\begin{equation}\label{eq:bochner-shape-one-form}
    \frac12\Delta^{\widehat{\Sigma}}\|\omega_A\|^2-\|\nabla^{\widehat{\Sigma}}\omega_A\|^2
    =\mathrm{ric}^{\widehat{\Sigma}}\bigl(\hat{L}^*\hat{L}(A),\hat{L}^*\hat{L}(A)\bigr).
\end{equation}
At $x_0$, we have $\Delta^{\widehat{\Sigma}}\|\omega_A\|^2\leq0$.
Evaluating~\eqref{eq:bochner-shape-one-form} at $x_0$ and using~\eqref{eq:shape-one-form-at-maximum}, we obtain
$$
    \lambda_0\,\mathrm{ric}^{\widehat{\Sigma}}(\widehat{X}_{x_0},\widehat{X}_{x_0})\leq0.
$$
Since $\lambda_0>0$, the Ricci curvature in this direction is non-positive.
The projection $\widehat{\Sigma}\rightarrow M/\mathsf{N}$ is a local isometry, so Theorem~\ref{thm:nilradical-metric-identities}\,(iii) gives
$$
    0\geq\mathrm{ric}^{\widehat{\Sigma}}(\widehat{X}_{x_0},\widehat{X}_{x_0})
    =-\|\widehat{X}_{x_0}\|^2
    +\langle\hat{L}^*\hat{L}_{x_0}(\widehat{X}_{x_0}),\widehat{X}_{x_0}\rangle
    =-1+\lambda_0.
$$
Hence $\lambda_0\leq1$, as desired.
\end{proof}

\subsection{Proof of Theorem~\ref{thm:solvable-splitting}}\label{subsec:proof-solvable-splitting}

\textcolor{white}{aaa}

We now use the preceding estimate to prove polarity and the vanishing of the shape operator of the $\mathsf{G}$-orbits.
We then eliminate the $\mathsf{G}$-holonomy to obtain the global splitting.

\begin{prop}\label{prop:polarity-a-directions}
The $\mathsf{G}$-action is polar.
Moreover, for every $\mathcal{A}$-valued vector field $A$ and every $\mathsf{G}$-horizontal vector field $X$,
$$
    \mathrm{ric}^{\mathsf{G}}(A,A)=-\|A\|^2,
    \qquad L_X(A)=0.
$$
\end{prop}

\begin{proof}
By Lemma~\ref{lem:orbit-ricci-a-directions} and Theorem~\ref{thm:orbit-ricci-lower-bound}, we have
$$
    \mathrm{ric}^{\mathsf{G}}(A,A)
    =-\mathrm{Tr}(\hat{L}_A^2)
    =-\|\hat{L}_A\|^2
    \geq-\|A\|^2.
$$

Since the claimed identities are pointwise, Proposition~\ref{prop:invariant-a-extension} allows us to assume that $A$ is invariant and horizontally Killing.
Consider the invariant function $f:=\|A\|^2/2$.
Lemma~\ref{lem:horizontally-killing-norm} gives
$$
    \Delta_Bf-\langle\nabla f,H\rangle
    =\mathrm{ric}^{\mathsf{G}}(A,A)+\|A\|^2
    +\|\mathrm{A}A\|^2+2\|LA\|^2\geq0.
$$
Proposition~\ref{prop:cocompact-maximum-principle}\,(ii) implies that $f$ is constant, so every nonnegative term on the right-hand side vanishes.
Thus, for every $\mathsf{G}$-horizontal vector field $X$,
$$
    \mathrm{ric}^{\mathsf{G}}(A,A)=-\|A\|^2,
    \qquad\mathrm{A}_X(A)=0,
    \qquad L_X(A)=0.
$$
By Lemma~\ref{lem:vertical-tensor-blocks}, $\mathrm{A}_X$ also vanishes on $\widehat{\mathcal{V}}$.
Since $\mathcal{V}=\mathcal{A}\oplus\widehat{\mathcal{V}}$, we conclude that $\mathrm{A}=0$.
Hence the $\mathsf{G}$-action is polar.
\end{proof}

We next combine the equality in $\mathcal{A}$-directions with Theorem~\ref{thm:orbit-ricci-lower-bound} to show that the shape operator also vanishes on $\widehat{\mathcal{V}}$.

\begin{prop}\label{prop:totally-geodesic-orbits}
The map $\hat{h}:M\rightarrow\mathrm{Sym}^2_+(\mathfrak{n}^*)$ is constant along $\mathsf{G}$-horizontal curves, and the $\mathsf{G}$-orbits in $M$ are totally geodesic submanifolds.
\end{prop}

\begin{proof}
Polarizing the identity in Proposition~\ref{prop:polarity-a-directions} and using Lemma~\ref{lem:orbit-ricci-a-directions}, we obtain
$$
    \langle\hat{L}^*\hat{L}(A),\tilde{A}\rangle
    =\mathrm{Tr}(\hat{L}_A\cdot\hat{L}_{\tilde{A}})
    =-\mathrm{ric}^{\mathsf{G}}(A,\tilde{A})
    =\langle A,\tilde{A}\rangle
$$
for every pair of $\mathcal{A}$-valued vector fields $A,\tilde{A}$.
Thus, the $\mathcal{A}$-component of $\hat{L}^*\hat{L}(A)$ is $A$.
The estimate in Theorem~\ref{thm:orbit-ricci-lower-bound} also gives $\|\hat{L}^*\|_{\mathrm{op}}\leq1$.
It follows that
$$
    \|A\|^2+\|\mathcal{H}(\hat{L}^*\hat{L}(A))\|^2
    =\|\hat{L}^*\hat{L}(A)\|^2
    \leq\|\hat{L}_A\|^2
    =\|A\|^2.
$$
Therefore, $\mathcal{H}(\hat{L}^*\hat{L}(A))=0$, and hence
$$
    \hat{L}^*\hat{L}(A)=A
    \qquad\text{for every $\mathcal{A}$-valued vector field }A.
$$

Let $X$ be a $\mathsf{G}$-horizontal vector field.
By Proposition~\ref{prop:shape-alpha},
$$
    \|\hat{L}_X\|^2
    =\mathrm{Tr}(\hat{L}_X\cdot\hat{L}_{\alpha(X)})
    =\langle X,\hat{L}^*\hat{L}(\alpha(X))\rangle
    =\langle X,\alpha(X)\rangle=0.
$$
Thus, $\hat{L}_X=0$ for every $\mathsf{G}$-horizontal vector field $X$.
Lemma~\ref{lem:orbit-metric-differential} for the $\mathsf{N}$-action implies that $\hat{h}$ is constant along $\mathsf{G}$-horizontal curves.
Finally, Lemma~\ref{lem:vertical-tensor-blocks} and Proposition~\ref{prop:polarity-a-directions} give $L_X|_{\widehat{\mathcal{V}}}=\hat{L}_X=0$ and $L_X|_{\mathcal{A}}=0$.
Hence $L=0$, so the $\mathsf{G}$-orbits are totally geodesic.
\end{proof}

It remains to show that the local product structure has trivial holonomy, so that it gives a global splitting.

\begin{prop}\label{prop:trivial-g-holonomy}
The holonomy of the principal bundle $M\rightarrow M/\mathsf{G}$, with respect to the $\mathsf{G}$-horizontal distribution, is trivial.
\end{prop}

\begin{proof}
Fix $x\in M$, and let $\widehat{\Sigma}$ be the $\mathsf{N}$-horizontal section through $x$, with $\mathfrak{a}_0=\mathfrak{a}_x$ and $\mathsf{A}_0=\exp(\mathfrak{a}_0)$.
Let $g_{\mathrm{h}}\in\mathsf{Hol}_x$.
By Theorem~\ref{thm:trivial-n-holonomy}, we have $g_{\mathrm{h}}\in\mathsf{A}_0$.

By Proposition~\ref{prop:totally-geodesic-orbits}, $\hat{h}$ is constant along $\mathsf{G}$-horizontal curves.
Its equivariance therefore gives
$$
    \mathrm{Ad}(g_{\mathrm{h}})|_{\mathfrak{n}}\cdot\hat{h}_x
    =\hat{h}_{g_{\mathrm{h}}\cdot x}
    =\hat{h}_x.
$$
Thus, $g_{\mathrm{h}}$ fixes $\hat{h}_x$ under the $\mathsf{A}_0$-action on $\mathrm{Sym}^2_+(\mathfrak{n}^*)$.
This action is free by Lemma~\ref{lem:a0-action-free-proper}, so $g_{\mathrm{h}}=e$.
\end{proof}

\begin{proof}[Proof of Theorem~\ref{thm:solvable-splitting}]
This follows from Propositions~\ref{prop:totally-geodesic-orbits} and~\ref{prop:trivial-g-holonomy}, together with Proposition~\ref{prop:splitting-criterion-outline}.
\end{proof}

\section{Proof of Theorem A: the general case}\label{sec:general-splitting}

We will now prove the following equivariant refinement of Theorem~\ref{thm:main-splitting}.

\begin{theorem}\label{thm:equivariant-main-splitting}
Let $M$ be a complete, connected negative-Einstein manifold, and let $\mathsf{G}\subseteq\mathsf{Isom}(M)$ be a closed connected subgroup acting cocompactly on $M$.
Then the following assertions hold:
\begin{enumerate}
\renewcommand{\labelenumi}{(\roman{enumi})}
    \item Every orbit is principal, and all orbits are equivariantly isometric to a fixed $\mathsf{G}$-homogeneous manifold $H$.
    \item $M$ is equivariantly isometric to $H\times M/\mathsf{G}$, where $\mathsf{G}$ acts trivially on the second factor.
    \item There exists a closed, connected, simply connected, solvable subgroup $\mathsf{S}\subseteq\mathsf{G}$ such that $\mathsf{S}\cdot x=\mathsf{G}\cdot x$ for every $x\in M$.
\end{enumerate}
\end{theorem}

The proof has two steps.
In the first step, we assume that $\mathsf{G}$ has finite center and use the generalized Iwasawa decomposition of Theorem~\ref{thm:generalized-iwasawa} to find a solvable subgroup to which Theorem~\ref{thm:solvable-splitting} applies.
In the second step, we show that $\mathsf{G}$ always has finite center.
After proving its discreteness, we use that the center is finitely generated and hence decomposes as the product of a torsion-free subgroup and a finite subgroup.
We quotient by the torsion-free subgroup and apply the finite-center case.
The resulting simply connected orbits force this subgroup to be trivial, showing that the center is finite.
We conclude this section with consequences of the splitting.

\subsection{The finite-center case}\label{subsec:general-splitting-finite-center}

\textcolor{white}{aaa}

\begin{theorem}\label{thm:finite-center-splitting}
Let $M$ be a complete, connected negative-Einstein manifold, and let $\mathsf{G}\subseteq\mathsf{Isom}(M)$ be a closed connected subgroup with finite center acting cocompactly on $M$.
Then the assertions of Theorem~\ref{thm:equivariant-main-splitting} hold.
\end{theorem}

\begin{proof}
By Theorem~\ref{thm:generalized-iwasawa}, we can write $\mathsf{G}=\mathsf{K}\mathsf{S}$, where $\mathsf{K}$ is compact and $\mathsf{S}$ is a closed, connected, simply connected solvable subgroup.
We first verify that the $\mathsf{S}$-action satisfies the hypotheses of Theorem~\ref{thm:solvable-splitting}.
Since $\mathsf{S}$ is closed in $\mathsf{Isom}(M)$, its action is proper and hence free by Lemma~\ref{lem:simply-connected-solvable-free-action}.
Moreover, $\mathsf{G}/\mathsf{S}$ is compact, being the image of $\mathsf{K}$ under the quotient map, so Lemma~\ref{lem:cocompact-subgroup} implies that the $\mathsf{S}$-action is cocompact.
Theorem~\ref{thm:solvable-splitting} therefore gives an $\mathsf{S}$-equivariant isometric splitting $M\cong H\times M/\mathsf{S}$, under which the $\mathsf{S}$-orbits are the submanifolds $H\times\{b\}$, for $b\in M/\mathsf{S}$.
Both factors are negative-Einstein.

We next show that $\mathsf{S}\cdot x=\mathsf{G}\cdot x$ for every $x\in M$.
Since $\mathsf{G}$ is connected, it suffices to prove that $U^*$ is $\mathsf{S}$-vertical for every $U\in\mathfrak{g}$.
Under the product identification $M\cong H\times M/\mathsf{S}$, the $\mathsf{S}$-horizontal component of $U^*$ restricts, for each $p\in H$, to a Killing field on $\{p\}\times M/\mathsf{S}$.
Since $M/\mathsf{S}$ is compact and negative-Einstein, Bochner's theorem implies that this component vanishes; see \cite{bochner1946vector}.
Thus, the orbits coincide, so $M/\mathsf{S}=M/\mathsf{G}$, and the two actions have the same vertical and horizontal distributions.

Fix $x_0\in M$, and let $\Sigma\subseteq M$ be the connected horizontal section through $x_0$ arising from the product decomposition.
Identifying $\mathsf{S}\cdot x_0$ with $H$ and $\Sigma$ with $M/\mathsf{S}$, the $\mathsf{S}$-equivariant isometry given by Theorem~\ref{thm:solvable-splitting} takes the form
$$
    \Phi:(\mathsf{S}\cdot x_0)\times\Sigma\longrightarrow M,
    \qquad \Phi(s\cdot x_0,x)=s\cdot x.
$$
We claim that $\Phi$ is $\mathsf{G}$-equivariant, where $\mathsf{G}$ acts on the first factor by its action on $\mathsf{G}\cdot x_0=\mathsf{S}\cdot x_0$.

Since the $\mathsf{G}$- and $\mathsf{S}$-orbits coincide and $\mathsf{S}$ acts freely, for every $(g,x)\in \mathsf{G}\times M$ there is a unique $\mu(g,x)\in \mathsf{S}$ such that $g\cdot x=\mu(g,x)\cdot x$.
The map $\mu:\mathsf{G}\times M\to \mathsf{S}$ is smooth: in a local trivialization of $M\to M/\mathsf{S}$, it is obtained by multiplying the $\mathsf{S}$-coordinate of $g\cdot x$ by the inverse of the $\mathsf{S}$-coordinate of $x$.

To prove that $\Phi$ is $\mathsf{G}$-equivariant, it is enough to show that $\mu$ is constant along horizontal curves in the second variable.
Indeed, this would give $\mu(g,s\cdot x)=\mu(g,s\cdot x_0)$ for every $s\in \mathsf{S}$ and $x\in\Sigma$, since $s\cdot\Sigma$ is connected and horizontal.
We would then have
\begin{align*}
    g\cdot\Phi(s\cdot x_0,x)
    &=g\cdot s\cdot x
    =\mu(g,s\cdot x)\cdot s\cdot x\\
    &=\mu(g,s\cdot x_0)\cdot s\cdot x
    =\Phi(\mu(g,s\cdot x_0)\cdot s\cdot x_0,x)
    =\Phi(g\cdot s\cdot x_0,x),
\end{align*}
as required.
We now prove this constancy.

Recall that $m_g:M\to M$ denotes the map $x\mapsto g\cdot x$, and that $r_x:\mathsf{G}\to M$ denotes the orbit map $g\mapsto g\cdot x$.
Fix $g\in \mathsf{G}$ and differentiate $g\cdot x=\mu(g,x)\cdot x$ in a horizontal direction $X_x$.
This gives
\begin{equation}\label{eq:horizontal-derivative-action-coefficient}
    d(m_g)_x(X_x)
    =
    d(m_{\mu(g,x)})_x(X_x)
    +
    d(r_x)_{\mu(g,x)}\left(\partial_{(0,X_x)}\mu(g,x)\right).
\end{equation}
Notice that $d(m_g)_x(\mathcal{H}_x)\subseteq\mathcal{H}_{g\cdot x}$, while $d(r_x)_g$ always takes values in the vertical space.
Hence, the first two terms in \eqref{eq:horizontal-derivative-action-coefficient} are horizontal and the last term is vertical, so the last term must vanish.
Since $\partial_{(0,X_x)}\mu(g,x)\in T_{\mu(g,x)}\mathsf{S}$ and the $\mathsf{S}$-action is free, this implies $\partial_{(0,X_x)}\mu(g,x)=0$.

Thus, $\mu$ is constant along horizontal curves in the second variable, and $\Phi$ is $\mathsf{G}$-equivariant.
In particular, the restrictions of $\Phi$ to the submanifolds $(\mathsf{S}\cdot x_0)\times\{x\}$ give $\mathsf{G}$-equivariant isometries between all orbits.

Finally, let $y=\Phi(s\cdot x_0,x)\in M$.
By the $\mathsf{G}$-equivariance and injectivity of $\Phi$, an element $g\in \mathsf{G}$ fixes $y$ if and only if it fixes $s\cdot x_0$.
Thus, $\mathsf{G}_y=s\mathsf{G}_{x_0}s^{-1}$.
All isotropy subgroups are therefore conjugate, and every $\mathsf{G}$-orbit is principal.
Since $\Sigma$ is isometric to $M/\mathsf{S}=M/\mathsf{G}$, this completes the proof.
\end{proof}

\subsection{Reduction to finite center}\label{subsec:general-splitting-center-reduction}

\textcolor{white}{aaa}

We first prove that the center is discrete.
In the homogeneous case, this is due to Jablonski and Petersen \cite[Corollary~1.8]{jablonski2017step}.

\begin{lemma}\label{lem:discrete-center}
Let $\mathsf{G}$ be a connected Lie group acting effectively, isometrically and cocompactly on a Riemannian manifold $M$ with negative Ricci curvature.
Then $\mathfrak{z}(\mathfrak{g})=\{0\}$, or equivalently, $Z(\mathsf{G})$ is discrete.
\end{lemma}

\begin{proof}
Let $U\in\mathfrak{z}(\mathfrak{g})$. Connectedness of $\mathsf{G}$ implies that $\mathrm{Ad}(g)U=U$ for every $g\in\mathsf{G}$. Hence $U^*$ and $f:=\|U^*\|^2/2$ are $\mathsf{G}$-invariant. The function $f$ descends to the compact orbit space $M/\mathsf{G}$ and therefore attains a maximum at some $x_0\in M$.

Since $U^*$ is Killing, Bochner's formula gives
$$
\Delta f=\|\nabla U^*\|^2-\mathrm{ric}^M(U^*,U^*)\geq0;
$$
see \cite{bochner1946vector} or \cite[Lemma~1.5]{jablonski2017step}. At $x_0$ we also have $\Delta f(x_0)\leq0$ and hence $\mathrm{ric}^M(U^*_{x_0},U^*_{x_0})=0$. Negative definiteness of the Ricci tensor then implies $U^*_{x_0}=0$. Thus, the maximum of the nonnegative function $f$ is zero, and $U^*\equiv0$. Effectiveness gives $U=0$.
\end{proof}

Every discrete normal subgroup of a connected Lie group is finitely generated; see \cite[Corollary~8.A.23(2)]{cornulier2016metric}.
Thus, when $Z(\mathsf{G})$ is discrete, we may write $Z(\mathsf{G})=\Gamma\cdot\mathsf{F}$, where $\Gamma\cong\mathbb{Z}^k$ for some $k\geq0$, $\mathsf{F}$ is the finite torsion subgroup, and $\Gamma\cap\mathsf{F}=\{e\}$.
The next lemma allows us to apply the finite-center case to the quotient by $\Gamma$.
This construction is inspired by \cite[Lemma~5.1]{bohm2025non}.

\begin{lemma}\label{lem:finite-center-reduction}
Let $\mathsf{G}$ be a connected Lie group with discrete center, and write $Z(\mathsf{G})=\Gamma\cdot\mathsf{F}$ as above.
Then $\mathsf{G}/\Gamma$ is connected and has center isomorphic to $\mathsf{F}$.

Suppose, moreover, that $\mathsf{G}$ acts properly and isometrically on a connected Riemannian manifold $M$.
Then $\Gamma$ acts freely and properly on $M$, the quotient $M/\Gamma$ is a connected Riemannian manifold, and the projection $M\to M/\Gamma$ is a Riemannian covering.
The induced action of $\mathsf{G}/\Gamma$ on $M/\Gamma$ is well defined, proper and isometric.
In addition:
\begin{enumerate}
\renewcommand{\labelenumi}{(\roman{enumi})}
    \item If the $\mathsf{G}$-action is effective and cocompact, then so is the induced $\mathsf{G}/\Gamma$-action.
    \item If the $\mathsf{G}/\Gamma$-orbits in $M/\Gamma$ are simply connected, then $\Gamma=\{e\}$, and hence $\mathsf{G}$ has finite center.
\end{enumerate}
\end{lemma}

\begin{proof}
We first show that $Z(\mathsf{G}/\Gamma)=Z(\mathsf{G})/\Gamma$.
If $g_0\in\mathsf{G}$ represents a central element of $\mathsf{G}/\Gamma$, then $gg_0g^{-1}g_0^{-1}\in\Gamma$ for every $g\in\mathsf{G}$.
The commutator map $g\mapsto gg_0g^{-1}g_0^{-1}$ is continuous, with connected domain and discrete image, and takes the value $e$ at the identity.
It is therefore identically $e$, so $g_0\in Z(\mathsf{G})$.
The reverse inclusion is immediate.
Hence $Z(\mathsf{G}/\Gamma)=Z(\mathsf{G})/\Gamma\cong\mathsf{F}$.

Properness of the $\mathsf{G}$-action implies that each isotropy subgroup $\mathsf{G}_x$ is compact.
Since $\Gamma$ is closed and discrete in $\mathsf{G}$, the subgroup $\Gamma\cap \mathsf{G}_x$ is compact and discrete, and hence finite.
As $\Gamma$ is torsion-free, $\Gamma\cap \mathsf{G}_x=\{e\}$, so the $\Gamma$-action is free.
It is also proper, since $\Gamma$ is closed in $\mathsf{G}$.
Thus, $M/\Gamma$ carries the quotient Riemannian metric for which the projection $\pi:M\to M/\Gamma$ is a Riemannian covering.

Since $\Gamma$ is central, the formula $\bar g\cdot\bar x:=\overline{g\cdot x}$ defines a smooth action of $\mathsf{G}/\Gamma$ on $M/\Gamma$.
We verify the remaining properties separately.

\emph{Isometric.}
For $g\in \mathsf{G}$, denote its action on $M$ by $m_g$ and the action of $\bar g\in \mathsf{G}/\Gamma$ on $M/\Gamma$ by $m_{\bar g}$.
The identity $m_{\bar g}\circ\pi=\pi\circ m_g$, together with the fact that $\pi$ is a local isometry and $m_g$ is an isometry, gives
$$
    \left\langle d(m_{\bar g})_{\bar x}\bigl(d\pi_x(v)\bigr),d(m_{\bar g})_{\bar x}\bigl(d\pi_x(w)\bigr)\right\rangle
    =\langle v,w\rangle
    =\left\langle d\pi_x(v),d\pi_x(w)\right\rangle
$$
for every $x\in M$ and $v,w\in T_xM$.
Since $d\pi_x$ is an isomorphism, the induced action is isometric.

\emph{Proper.}
Suppose that $\bar x_n\to\bar x$ in $M/\Gamma$ and that $\bar g_n\cdot\bar x_n\to\bar y$ for a sequence $\bar g_n\in \mathsf{G}/\Gamma$.
Using local sections of the covering $\pi$, choose representatives $x_n\to x$ and $y_n\to y$ in $M$ such that $\pi(x_n)=\bar x_n$ and $\pi(y_n)=\bar g_n\cdot\bar x_n$.
Choose any representatives $g_n\in \mathsf{G}$ of $\bar g_n$.
There exists $\gamma_n\in\Gamma$ such that $y_n=\gamma_ng_n\cdot x_n$.
Properness of the $\mathsf{G}$-action implies that $(\gamma_ng_n)$ has a convergent subsequence.
Passing to the quotient, the corresponding subsequence of $\bar g_n=\overline{\gamma_ng_n}$ also converges.
This proves properness.

\emph{Effective.}
Assume that the $\mathsf{G}$-action is effective, and suppose that $\bar g\in \mathsf{G}/\Gamma$ acts trivially on $M/\Gamma$.
Then $m_g$ is a deck transformation of the covering $\pi$.
Choose $x_0\in M$ and $\gamma\in\Gamma$ such that $g\cdot x_0=\gamma\cdot x_0$.
The maps $m_g,m_\gamma:M\to M$ are lifts of the same map $\pi:M\to M/\Gamma$ and agree at $x_0$.
Since $M$ is connected, uniqueness of lifts gives $m_g=m_\gamma$.
Effectiveness of the $\mathsf{G}$-action implies $g=\gamma$, so $\bar g=\bar e$.

\emph{Cocompact.}
The projection $\pi$ induces a homeomorphism $M/\mathsf{G}\to(M/\Gamma)/(\mathsf{G}/\Gamma)$ given by $\mathsf{G}\cdot x\mapsto(\mathsf{G}/\Gamma)\cdot\bar x$.
Thus, compactness of $M/\mathsf{G}$ implies compactness of $(M/\Gamma)/(\mathsf{G}/\Gamma)$.

Finally, suppose that the $\mathsf{G}/\Gamma$-orbits are simply connected.
Let $\gamma\in\Gamma$.
Since $\mathsf{G}$ is connected, choose a path $g:[0,1]\to \mathsf{G}$ with $g(0)=e$ and $g(1)=\gamma$.
For any $x\in M$, the path $\sigma(t):=g(t)\cdot x$ projects to a loop in $(\mathsf{G}/\Gamma)\cdot\bar x$.
This loop is null-homotopic in that orbit, and hence in $M/\Gamma$.
Its lift $\sigma$ is therefore closed, so $\gamma\cdot x=x$.
Freeness of the $\Gamma$-action implies $\gamma=e$, proving $\Gamma=\{e\}$.
\end{proof}

\begin{proof}[Proof of Theorem~\ref{thm:equivariant-main-splitting}]
By Lemma~\ref{lem:discrete-center}, $\mathsf{G}$ has discrete center.
Write $Z(\mathsf{G})=\Gamma\cdot\mathsf{F}$ as in Lemma~\ref{lem:finite-center-reduction}.
That lemma shows that $\Gamma$ acts freely, properly and isometrically on $M$.
Since $M\to M/\Gamma$ is a Riemannian covering and $M$ is complete, $M/\Gamma$ is a complete, connected negative-Einstein manifold.
Moreover, $\mathsf{G}/\Gamma$ is connected, has finite center, and acts effectively, properly, isometrically and cocompactly on $M/\Gamma$.
By effectiveness and properness, we may identify $\mathsf{G}/\Gamma$ with a closed connected subgroup of $\mathsf{Isom}(M/\Gamma)$.
Theorem~\ref{thm:finite-center-splitting} therefore applies to this action.

In particular, there is a closed, connected, simply connected solvable subgroup $\bar{\mathsf{S}}\subseteq\mathsf{G}/\Gamma$ acting transitively on each $\mathsf{G}/\Gamma$-orbit.
The $\bar{\mathsf{S}}$-action is proper, since $\bar{\mathsf{S}}$ is closed, and hence free by Lemma~\ref{lem:simply-connected-solvable-free-action}.
Since $\bar{\mathsf{S}}$ acts simply transitively on each $\mathsf{G}/\Gamma$-orbit, each such orbit is diffeomorphic to $\bar{\mathsf{S}}$ and is, in particular, simply connected.
Lemma~\ref{lem:finite-center-reduction}\,(ii) now gives $\Gamma=\{e\}$, so $\mathsf{G}$ has finite center.
A final application of Theorem~\ref{thm:finite-center-splitting} completes the proof.
\end{proof}

Theorem~\ref{thm:main-splitting} follows by taking $\mathsf{G}=\mathsf{Isom}(M)^o$: as recalled in Subsection~\ref{subsec:properties-isometric-actions}, cocompact symmetries imply that $M$ is complete and that $M/\mathsf{Isom}(M)^o$ is compact.

\begin{remark}\label{rem:trivial-center}
For $\mathsf{G}$ as in Theorem~\ref{thm:equivariant-main-splitting}, the preceding proof shows that $Z(\mathsf{G})$ is finite.
In fact, $\mathsf{G}$ is centerless.
Let $z\in Z(\mathsf{G})$ have order $m$, and fix $x\in M$.
The subgroup $\mathsf{S}$ given by Theorem~\ref{thm:equivariant-main-splitting} acts freely and transitively on each $\mathsf{G}$-orbit, so there exists $s\in\mathsf{S}$ such that $z\cdot x=s\cdot x$.
Since $z$ is central, we have
$$
    x=z^m\cdot x=s^m\cdot x.
$$
Freeness of the $\mathsf{S}$-action gives $s^m=e$.
By Lemma~\ref{lem:simply-connected-solvable-free-action}, $\mathsf{S}$ has no nontrivial compact subgroups, so $s=e$.
Thus, $z\cdot x=x$ for every $x\in M$, and effectiveness implies $z=e$.
\end{remark}

\begin{remark}
For $\mathsf{G}$ as in Theorem~\ref{thm:equivariant-main-splitting}, we obtain an Iwasawa-like decomposition $\mathsf{G}=\mathsf{K}\mathsf{A}_0\mathsf{N}_0$.
Indeed, Theorem~\ref{thm:generalized-iwasawa} gives a decomposition $\mathsf{G}=\mathsf{K}\mathsf{S}$, where $\mathsf{K}$ is connected and maximal compact and $\mathsf{S}$ is closed, connected, simply connected and solvable.
The proof of Theorem~\ref{thm:finite-center-splitting} shows that the $\mathsf{S}$-action is free, proper and cocompact.
Applying Corollary~\ref{cor:solvable-exponential-decomposition}\,(iii) to this action gives $\mathsf{S}=\mathsf{A}_0\ltimes\mathsf{N}_0$, where $\mathsf{N}_0$ is the connected nilradical of $\mathsf{S}$, which need not coincide with the connected nilradical of $\mathsf{G}$.
The subgroup $\mathsf{A}_0$ is closed, connected, simply connected and abelian, while $\mathsf{N}_0$ is closed, simply connected and nilpotent; see Corollary~\ref{cor:solvable-exponential-decomposition}\,(i)--(ii).
Both subgroups are closed in $\mathsf{G}$ because $\mathsf{S}$ is closed in $\mathsf{G}$.
Consequently, the multiplication map
$$
    \mathsf{K}\times\mathsf{A}_0\times\mathsf{N}_0\longrightarrow\mathsf{G},
    \qquad (k,a,n)\longmapsto kan
$$
is a diffeomorphism.
\end{remark}

\subsection{Consequences of the splitting}\label{subsec:general-splitting-consequences}

\textcolor{white}{aaa}

We first show that the orbit decomposition is independent of the choice of closed connected cocompact symmetry group.
We then prove Corollary~\ref{cor:parallel-ricci-splitting}.
All other corollaries stated in the Introduction have already been proved there.

\begin{corollary}\label{cor:canonical-orbit-decomposition}
Let $M$ be a complete, connected negative-Einstein manifold, and let $\mathsf{G}_1,\mathsf{G}_2\subseteq\mathsf{Isom}(M)$ be closed connected subgroups acting cocompactly on $M$.
Then $\mathsf{G}_1\cdot x=\mathsf{G}_2\cdot x$ for every $x\in M$.
\end{corollary}

\begin{proof}
Since $\mathsf{G}_1,\mathsf{G}_2\subseteq\mathsf{Isom}(M)^o$, it suffices to show that $\mathsf{G}_1\cdot x=\mathsf{Isom}(M)^o\cdot x$ for every $x\in M$.
By Theorem~\ref{thm:equivariant-main-splitting}, there is a $\mathsf{G}_1$-equivariant isometric splitting $M\cong H\times M/\mathsf{G}_1$, under which the $\mathsf{G}_1$-orbits are the submanifolds $H\times\{b\}$, for $b\in M/\mathsf{G}_1$.
The second factor is compact and negative-Einstein.
As in the proof of Theorem~\ref{thm:finite-center-splitting}, the component of any Killing field tangent to the second factor restricts to a Killing field on $\{p\}\times M/\mathsf{G}_1$ for each $p\in H$, and hence vanishes by Bochner's theorem.
Thus, every fundamental field of $\mathsf{Isom}(M)^o$ is tangent to the $\mathsf{G}_1$-orbits.
Connectedness of $\mathsf{Isom}(M)^o$ yields the desired equality of orbits.
\end{proof}

\begin{proof}[Proof of Corollary~\ref{cor:parallel-ricci-splitting}]
Since $M$ has cocompact symmetries, it is complete.
The Ricci endomorphism is self-adjoint and parallel.
Connectedness therefore implies that its distinct eigenvalues $\lambda_1,\ldots,\lambda_k$ and their multiplicities are constant.
The corresponding smooth eigendistributions give a parallel orthogonal decomposition $TM=E_1\oplus\cdots\oplus E_k$.
Completeness and simple connectedness allow us to apply the de Rham decomposition theorem, obtaining an isometric splitting $M\cong M_1\times\cdots\times M_k$; see \cite{derham1952reductibilite,Pet16}.
Each $M_i$ is complete, connected, simply connected and Einstein, with Einstein constant $\lambda_i$; see \cite[Theorem~1.100]{Besse87}.

We next show that each factor has cocompact symmetries.
Every isometry $f$ of $M$ preserves the Ricci endomorphism.
Since the eigenvalues are distinct, $df(E_i)=E_i$ for every $i$.
This implies that
\begin{equation}\label{eq:ricci-product-isometry-group}
\mathsf{Isom}(M)\cong\mathsf{Isom}(M_1)\times\cdots\times\mathsf{Isom}(M_k).
\end{equation}
Indeed, let $\pi_i:M\to M_i$ be the projection, choose base points in the other factors, and let $\iota_i:M_i\to M$ be the corresponding inclusion.
Set $f_i:=\pi_i\circ f\circ\iota_i$.
For $j\neq i$, we have $df(E_j)=E_j$ and $d\pi_i(E_j)=0$, so $d(\pi_i\circ f)(E_j)=0$.
Connectedness of the factors implies that $\pi_i\circ f$ depends only on the $M_i$-coordinate.
Hence
\begin{equation}\label{eq:product-isometry-components}
f(x_1,\ldots,x_k)=(f_1(x_1),\ldots,f_k(x_k)).
\end{equation}
Applying the same argument to $f^{-1}$ shows that each $f_i$ is a diffeomorphism.
Since $f$ preserves the product metric, each $f_i$ is an isometry.
Conversely, every product of isometries of the factors is an isometry of $M$, proving \eqref{eq:ricci-product-isometry-group}.

Taking identity components in \eqref{eq:ricci-product-isometry-group} gives $\mathsf{Isom}(M)^o\cong\mathsf{Isom}(M_1)^o\times\cdots\times\mathsf{Isom}(M_k)^o$.
Thus, each $\pi_i$ induces a continuous surjection $M/\mathsf{Isom}(M)^o\to M_i/\mathsf{Isom}(M_i)^o$.
Its domain is compact, so each $M_i$ has cocompact symmetries.

We now distinguish the signs of the Einstein constants.
If $\lambda_i>0$, then $M_i$ is compact by Myers's theorem \cite{myers1941riemannian}.
If $\lambda_i=0$, Theorem~9.2 of \cite{bohm2025non} and simple connectedness give an isometric splitting $M_i\cong K_i\times\mathbb{R}^{n_i}$, with $K_i$ compact and Ricci-flat.
If $\lambda_i<0$, Theorem~\ref{thm:main-splitting} gives $M_i\cong H_i\times K_i$, where $H_i$ is homogeneous and both factors are negative-Einstein, with $K_i$ compact.
Taking the product of all compact factors gives $K$, and taking the product of all Euclidean and homogeneous negative-Einstein factors gives $H$.
Thus, $M\cong H\times K$ with the asserted decompositions.
In particular, $H$ is homogeneous and has non-positive Ricci curvature.
\end{proof}

\begin{appendix}

\section{The symmetric space of inner products}\label{app:inner-products}
In this appendix, we recall the homogeneous-space description of the space of inner products on a finite-dimensional real vector space.
We then derive differentiation formulas for metric adjoints and Lie brackets in terms of pointwise self-adjoint representatives of vector fields.
Finally, we introduce the invariant Riemannian metric and compute its Levi-Civita connection.

\subsection{The homogeneous-space description}\label{subsec:inner-products-homogeneous-model}
\textcolor{white}{aaa}

Let $\mathfrak{g}$ be a finite-dimensional real vector space.
The space $\mathrm{Sym}^2_+(\mathfrak{g}^*)$ of inner products on $\mathfrak{g}$ is an open convex subset of the vector space $\mathrm{Sym}^2(\mathfrak{g}^*)$ of symmetric bilinear forms.
It therefore carries a natural smooth structure, with $T_h(\mathrm{Sym}^2_+(\mathfrak{g}^*))\equiv\mathrm{Sym}^2(\mathfrak{g}^*)$ for every $h\in\mathrm{Sym}^2_+(\mathfrak{g}^*)$.

The group $\mathsf{GL}(\mathfrak g)$ acts smoothly and transitively on $\mathrm{Sym}^2_+(\mathfrak g^*)$ by
\begin{equation}\label{eq:action}
    (q\cdot h)(U,V):=h(q^{-1}U,q^{-1}V),
    \qquad q\in\mathsf{GL}(\mathfrak g),\quad U,V\in\mathfrak g.
\end{equation}
The action is not effective: its ineffective kernel is $\{\pm \mathrm{Id}_{\mathfrak{g}}\}$.
The isotropy subgroup at $h$ is the compact orthogonal group $\mathsf{O}(\mathfrak{g},h)$.
Since the action is transitive with compact isotropy, it is proper by \cite[Corollary~3.1.1]{Field07}.
Thus, after fixing a background inner product $\bar{h}\in \mathrm{Sym}^2_+(\mathfrak{g}^*)$, the orbit map induces a $\mathsf{GL}(\mathfrak{g})$-equivariant diffeomorphism
$$
    \mathsf{GL}(\mathfrak g)/\mathsf{O}(\mathfrak g,\bar h)
    \longrightarrow \mathrm{Sym}^2_+(\mathfrak{g}^*),
    \qquad q\mathsf{O}(\mathfrak g,\bar h)\longmapsto q\cdot\bar h.
$$
To avoid confusion, we denote the Lie algebra of $\mathsf{GL}(\mathfrak{g})$ by $\mathrm{End}(\mathfrak{g})$ rather than $\mathfrak{gl}(\mathfrak{g})$.

For $h\in\mathrm{Sym}^2_+(\mathfrak{g}^*)$, denote by $E^{t_h}$ the $h$-adjoint of an endomorphism $E$, and set $\mathrm{S}_h(E):=(E+E^{t_h})/2$.
The corresponding Cartan decomposition is
$$
\mathrm{End}(\mathfrak{g})=\mathfrak{so}(\mathfrak{g},h)\oplus\mathrm{sym}(\mathfrak{g},h),
$$
where the two summands consist of the $h$-skew-adjoint and $h$-self-adjoint endomorphisms, respectively.

Differentiating \eqref{eq:action}, the fundamental vector field associated with $E\in\mathrm{End}(\mathfrak{g})$ is
$$
\rho(E)h:=\frac{d}{dt}\big|_{t=0} \ e^{tE}\cdot h
=-h(E\cdot,\cdot)-h(\cdot,E\cdot)
=-2h(\mathrm{S}_h(E)\cdot,\cdot).
$$
Thus, $\rho(E)h$ depends only on $\mathrm{S}_h(E)$, and the map
\begin{equation}\label{eq:symmetric-representative}
\mathrm{sym}(\mathfrak{g},h)\longrightarrow T_h(\mathrm{Sym}^2_+(\mathfrak{g}^*)),
\qquad E\longmapsto\rho(E)h,
\end{equation}
is a linear isomorphism.
Consequently, every smooth vector field $X$ on $\mathrm{Sym}^2_+(\mathfrak{g}^*)$ can be uniquely represented as $X_h=\rho(\mathrm{E}^X_h)h$ for a smooth function $\mathrm{E}^X:\mathrm{Sym}^2_+(\mathfrak{g}^*)\to\mathrm{End}(\mathfrak g)$ satisfying $\mathrm{E}^X_h\in\mathrm{sym}(\mathfrak g,h)$ for each $h$.
If $X_h=\rho(E_h)h$ for an arbitrary smooth map $E:\mathrm{Sym}^2_+(\mathfrak{g}^*)\to\mathrm{End}(\mathfrak g)$, its self-adjoint representative is $\mathrm{E}^X_h=\mathrm{S}_h(E_h)$.
\vskip2mm

{\bf{\noindent Convention.}}
Throughout this appendix, $\mathrm{E}^X$ always denotes the unique pointwise self-adjoint representative of the vector field $X$ given by \eqref{eq:symmetric-representative}.

\subsection{Differentiation formulas}\label{subsec:inner-products-differentiation}
\textcolor{white}{aaa}

The following proposition describes the variation of metric adjoints.
We write $\mathrm{sk}_h(E):=(E-E^{t_h})/2$ for the $h$-skew-adjoint part of an endomorphism $E$.

\begin{prop}\label{lem:inner-product-adjoint-differential}
Let $E:\mathrm{Sym}^2_+(\mathfrak{g}^*)\to\mathrm{End}(\mathfrak{g})$ be a smooth map, and define
$E^{t_*}:\mathrm{Sym}^2_+(\mathfrak{g}^*)\rightarrow\mathrm{End}(\mathfrak{g})$ by $h\mapsto E_h^{t_h}$. For every smooth vector field $X$ on $\mathrm{Sym}^2_+(\mathfrak{g}^*)$ and every $h\in\mathrm{Sym}^2_+(\mathfrak{g}^*)$, we have
$$
    d(E^{t_*})_h(X_h)
    =dE_h(X_h)^{t_h}+2[\mathrm{E}_h^X,E_h^{t_h}].
$$
In particular, if $E_h$ is $h$-self-adjoint for every $h$, then
$$
    \mathrm{sk}_h\bigl(dE_h(X_h)\bigr)=[\mathrm{E}_h^X,E_h].
$$
\end{prop}

\begin{proof}
By the definition of the metric adjoint,
$$
    h(E_hU,V)=h(U,E_h^{t_h}V)
    \qquad\text{for every }U,V\in\mathfrak{g}.
$$
Differentiating this identity in the direction $X_h$ gives
$$
    X_h(E_h\cdot,\cdot)+h\bigl(dE_h(X_h)\cdot,\cdot\bigr)
    =X_h(\cdot,E_h^{t_h}\cdot)+h\bigl(\cdot,d(E^{t_*})_h(X_h)\cdot\bigr).
$$
Substituting $X_h=-2h(\mathrm{E}_h^X\cdot,\cdot)$ gives
$$
    -2h(\mathrm{E}_h^XE_h\cdot,\cdot)+h\bigl(dE_h(X_h)\cdot,\cdot\bigr)
    =-2h(\mathrm{E}_h^X\cdot,E_h^{t_h}\cdot)+h\bigl(\cdot,d(E^{t_*})_h(X_h)\cdot\bigr).
$$
Moving all endomorphisms to the second argument of $h$ and using its nondegeneracy, we obtain
$$
    -2E_h^{t_h}\mathrm{E}_h^X+dE_h(X_h)^{t_h}
    =-2\mathrm{E}_h^XE_h^{t_h}+d(E^{t_*})_h(X_h).
$$
Rearranging proves the first formula.
If $E$ is pointwise self-adjoint, then $E^{t_*}=E$, so the first formula becomes
$$
    dE_h(X_h)-dE_h(X_h)^{t_h}=2[\mathrm{E}_h^X,E_h],
$$
which proves the second formula.
\end{proof}

The following consequence will be used in Section~\ref{sec:affine-conjugation-orbits}.

\begin{remark}\label{rem:inner-product-fixed-endomorphism}
For a fixed $E_0\in\mathrm{End}(\mathfrak{g})$, consider the map
$\mathrm{S}_*(E_0):\mathrm{Sym}^2_+(\mathfrak{g}^*)\rightarrow\mathrm{End}(\mathfrak{g})$ defined by $h\mapsto\mathrm{S}_h(E_0)$.
Applying Proposition~\ref{lem:inner-product-adjoint-differential} to the constant map with value $E_0$ gives, for every smooth vector field $X$ on $\mathrm{Sym}^2_+(\mathfrak{g}^*)$,
$$
    d(\mathrm{S}_*(E_0))_h(X_h)
    =[\mathrm{E}^X_h,E_0^{t_h}]
    \qquad\text{for every }h\in\mathrm{Sym}^2_+(\mathfrak{g}^*).
$$
\end{remark}

We now express the Lie bracket of vector fields in terms of their representatives.

\begin{lemma}\label{lem:inner-product-lie-bracket}
Let $X,Y$ be smooth vector fields on $\mathrm{Sym}^2_+(\mathfrak{g}^*)$.
Then
$$
    [X,Y]_h=\rho\bigl(d\mathrm{E}_h^Y(X_h)-d\mathrm{E}_h^X(Y_h)\bigr)h.
$$
\end{lemma}

\begin{proof}
Using the identification $T_h(\mathrm{Sym}^2_+(\mathfrak{g}^*))\equiv\mathrm{Sym}^2(\mathfrak{g}^*)$, we regard $X$ and $Y$ as maps from $\mathrm{Sym}^2_+(\mathfrak{g}^*)$ into the vector space $\mathrm{Sym}^2(\mathfrak{g}^*)$.
Their Lie bracket is therefore
$$
    [X,Y]_h=dY_h(X_h)-dX_h(Y_h).
$$
Differentiating $X_h=-2h(\mathrm{E}_h^X\cdot,\cdot)$ in the direction $Y_h=-2h(\mathrm{E}_h^Y\cdot,\cdot)$ gives
$$
    dX_h(Y_h)
    =-2Y_h(\mathrm{E}_h^X\cdot,\cdot)-2h\bigl(d\mathrm{E}_h^X(Y_h)\cdot,\cdot\bigr)
    =4h(\mathrm{E}_h^Y\mathrm{E}_h^X\cdot,\cdot)-2h\bigl(d\mathrm{E}_h^X(Y_h)\cdot,\cdot\bigr).
$$
The analogous expression for $dY_h(X_h)$ yields
\begin{equation}\label{eq:inner-product-bracket-expansion}
    [X,Y]_h
    =4h([\mathrm{E}_h^X,\mathrm{E}_h^Y]\cdot,\cdot)
    -2h\bigl((d\mathrm{E}_h^Y(X_h)-d\mathrm{E}_h^X(Y_h))\cdot,\cdot\bigr).
\end{equation}
By Proposition~\ref{lem:inner-product-adjoint-differential},
$$
    \mathrm{sk}_h\bigl(d\mathrm{E}_h^Y(X_h)\bigr)
    =[\mathrm{E}_h^X,\mathrm{E}_h^Y]
    =-\mathrm{sk}_h\bigl(d\mathrm{E}_h^X(Y_h)\bigr).
$$
Consequently,
$$
    \mathrm{S}_h\bigl(d\mathrm{E}_h^Y(X_h)-d\mathrm{E}_h^X(Y_h)\bigr)
    =d\mathrm{E}_h^Y(X_h)-d\mathrm{E}_h^X(Y_h)-2[\mathrm{E}_h^X,\mathrm{E}_h^Y].
$$
Since $\rho(F)h=-2h(\mathrm{S}_h(F)\cdot,\cdot)$ for every $F\in\mathrm{End}(\mathfrak{g})$, equation~\eqref{eq:inner-product-bracket-expansion} now gives the asserted formula.
\end{proof}

\subsection{The Riemannian metric and Levi-Civita connection}\label{subsec:inner-products-metric}

\textcolor{white}{aaa}

There is a canonical $\mathsf{GL}(\mathfrak{g})$-invariant metric on $\mathrm{Sym}^2_+(\mathfrak{g}^*)$.
Given smooth vector fields $X,Y$, set
$$
    \ip{X}{Y}_h
    :=
    \mathrm{Tr}(\mathrm{E}^X_h\cdot \mathrm{E}^Y_h).
$$

Equipped with this metric, $\mathrm{Sym}^2_+(\mathfrak{g}^*)$ is a Riemannian symmetric space of nonpositive sectional curvature.
In particular, it is complete, and it is also simply connected, hence a Hadamard manifold (see \cite[Section~3.1]{Heber98}).

\begin{prop}\label{lem:inner-product-connection}
Let $X,Y$ be smooth vector fields on $\mathrm{Sym}^2_+(\mathfrak{g}^*)$.
Then
$$
    (\nabla_XY)_h=\rho\bigl(d\mathrm{E}_h^Y(X_h)\bigr)h.
$$
\end{prop}

\begin{proof}
Define
$$
    (\widetilde\nabla_XY)_h:=\rho\bigl(d\mathrm{E}_h^Y(X_h)\bigr)h.
$$
This expression is $C^\infty$-linear in $X$.
For any smooth function $f$, we have $\mathrm{E}^{fY}=f\mathrm{E}^Y$, so
$$
    d\mathrm{E}_h^{fY}(X_h)=X_h(f)\mathrm{E}_h^Y+f(h)d\mathrm{E}_h^Y(X_h).
$$
Applying $\rho(\,\cdot\,)h$ gives the Leibniz rule
$$
    \widetilde\nabla_X(fY)=X(f)Y+f\widetilde\nabla_XY,
$$
and hence $\widetilde\nabla$ is a connection.
By Lemma~\ref{lem:inner-product-lie-bracket},
$$
    \widetilde\nabla_XY-\widetilde\nabla_YX=[X,Y],
$$
so this connection is torsion-free.

To verify compatibility with the metric, let $Z$ be a smooth vector field.
Differentiating the definition of the metric gives
$$
    X_h\ip{Y}{Z}
    =\mathrm{Tr}\!\left(d\mathrm{E}_h^Y(X_h)\cdot \mathrm{E}_h^Z\right)
    +\mathrm{Tr}\!\left(\mathrm{E}_h^Y\cdot d\mathrm{E}_h^Z(X_h)\right).
$$
The pointwise self-adjoint representative of $\widetilde\nabla_XY$ is $\mathrm{S}_h(d\mathrm{E}_h^Y(X_h))$.
Since $\mathrm{E}_h^Y$ and $\mathrm{E}_h^Z$ are $h$-self-adjoint, taking the $h$-symmetric part of either derivative does not change the corresponding trace.
Therefore,
$$
    X_h\ip{Y}{Z}
    =\ip{\widetilde\nabla_XY}{Z}_h+\ip{Y}{\widetilde\nabla_XZ}_h.
$$
Thus $\widetilde\nabla$ is compatible with the metric and torsion-free, so it agrees with the Levi-Civita connection.
\end{proof}

\section{The generalized Iwasawa decomposition}\label{app:generalized-iwasawa}

The goal of this appendix is to prove the following extension of the Iwasawa decomposition to connected Lie groups with finite center.

\begin{theorem}[Generalized Iwasawa decomposition]\label{thm:generalized-iwasawa}
Let $\mathsf{G}$ be a connected Lie group with finite center, and let $\mathsf{N}\subseteq\mathsf{G}$ be its connected nilradical.
For every maximal compact subgroup $\mathsf{K}\subseteq\mathsf{G}$, there exists a closed, connected, simply connected solvable subgroup $\mathsf{S}\subseteq\mathsf{G}$ containing $\mathsf{N}$ such that the multiplication map $\mathsf{K}\times\mathsf{S}\to\mathsf{G}$, $(k,s)\mapsto ks$, is a diffeomorphism.
\end{theorem}

When $\mathsf{G}$ is semisimple, the theorem is the usual Iwasawa decomposition $\mathsf{G}=\mathsf{K}\mathsf{A}\mathsf{N}_0$, with $\mathsf{S}:=\mathsf{A}\mathsf{N}_0$; see \cite[Theorem~6.46]{Knapp02}.
The finite-center assumption cannot be replaced by discreteness of the center.
Indeed, the universal covering group of $\mathsf{SL}_2(\mathbb{R})$ has infinite discrete center and no nontrivial compact subgroups.
The asserted decomposition would therefore force the group itself to be solvable.

Theorem~\ref{thm:generalized-iwasawa} can be deduced from the proof of \cite[Lemma~6.7]{Cornulier08}; we explain this deduction in Remark~\ref{rem:cornulier-iwasawa} at the end of this section.
W\"ustner~\cite{Wustner00} proves a related factorization with a solvable factor containing the nilradical and a compactly embedded factor, which need not itself be compact.

The theorem allows us to reduce the finite-center case of Theorem~\ref{thm:main-splitting} to the solvable case.
Indeed, let $M$ be a complete Riemannian manifold on which $\mathsf{G}$ acts effectively, properly, isometrically and cocompactly, and write $\mathsf{G}=\mathsf{K}\mathsf{S}$ as above.
The restricted isometric $\mathsf{S}$-action is proper because $\mathsf{S}$ is closed, and it is free by Lemma~\ref{lem:simply-connected-solvable-free-action}.
Moreover, $\mathsf{G}/\mathsf{S}$ is compact, so Lemma~\ref{lem:cocompact-subgroup} implies that this action is also cocompact.
Thus $M\to M/\mathsf{S}$ is a principal $\mathsf{S}$-bundle with compact base.
If $M$ is additionally connected and negative-Einstein, the $\mathsf{S}$-action satisfies all the hypotheses of Theorem~\ref{thm:solvable-splitting}.
We carry out this reduction in Section~\ref{sec:general-splitting}; see Theorem~\ref{thm:finite-center-splitting}.

To prove the theorem, we first show that the nilradical is simply connected when $Z(\mathsf{G})$ is discrete, equivalently $\mathfrak{z}(\mathfrak{g})=0$, and then obtain a decomposition of the solvable radical analogous to that in Theorem~\ref{thm:generalized-iwasawa}.
Finally, we combine this decomposition with an Iwasawa decomposition of a suitable Levi factor.
\subsection{The nilradical and the solradical}\label{subsec:iwasawa-radical}

\textcolor{white}{aaa}

We will use the Cartan--Iwasawa--Malcev theorem: a connected Lie group has a maximal compact subgroup, all maximal compact subgroups are conjugate, and the inclusion of a maximal compact subgroup is a homotopy equivalence.
In particular, maximal compact subgroups of connected solvable Lie groups are tori.

We also recall that a connected, simply connected solvable Lie group is diffeomorphic to Euclidean space and has no nontrivial compact subgroups by Lemma~\ref{lem:simply-connected-solvable-free-action}.

\begin{lemma}\label{lem:simply-connected-nilradical}
Let $\mathsf{G}$ be a connected Lie group with discrete center.
Then its connected nilradical $\mathsf{N}$ is simply connected.
\end{lemma}

\begin{proof}
Let $\mathsf{T}$ be a maximal compact subgroup of $\mathsf{N}$.
Since $\mathsf{N}$ is connected and nilpotent, $\mathsf{T}$ is a torus, and \cite[Lemma~3.7]{BL23} shows that it is central in $\mathsf{N}$.
Conjugacy of maximal compact subgroups therefore implies that $\mathsf{T}$ is the unique maximal compact subgroup of $\mathsf{N}$.

Normality of $\mathsf{N}$ in $\mathsf{G}$ implies that $g\mathsf{T}g^{-1}$ is a maximal compact subgroup of $\mathsf{N}$ for every $g\in\mathsf{G}$.
By uniqueness, $g\mathsf{T}g^{-1}=\mathsf{T}$, so conjugation gives a continuous homomorphism $\mathsf{G}\to\mathsf{Aut}(\mathsf{T})$.
Write $\mathsf{T}\cong\mathbb{T}^m$, where $\mathbb{T}^m:=(\mathbb{R}/\mathbb{Z})^m$ is the $m$-dimensional torus.
Then $\mathsf{Aut}(\mathsf{T})\cong\mathsf{GL}_m(\mathbb{Z})$ is discrete; see \cite[Chapter~IV, equation~(1.2)]{BrockerTomDieck85}.
Connectedness of $\mathsf{G}$ forces this homomorphism to be trivial.
Hence $\mathsf{T}\subseteq Z(\mathsf{G})$, and discreteness of $Z(\mathsf{G})$ together with connectedness of $\mathsf{T}$ gives $\mathsf{T}=\{e\}$.
The Cartan--Iwasawa--Malcev theorem then implies that $\mathsf{N}$ is simply connected.
\end{proof}

\begin{remark}\label{rem:nilradical-free-action}
If a connected Lie group $\mathsf{G}$ with discrete center acts properly and isometrically on a Riemannian manifold $M$, then its connected nilradical $\mathsf{N}$ acts freely, properly and isometrically.
Indeed, $\mathsf{N}$ is closed in $\mathsf{G}$, so its action is proper, and it is simply connected by Lemma~\ref{lem:simply-connected-nilradical}.
Freeness follows from Lemma~\ref{lem:simply-connected-solvable-free-action}.
\end{remark}

Let $\mathsf{R}\subseteq\mathsf{G}$ denote the connected solradical, that is, the largest connected solvable normal subgroup.
It is closed, and its connected nilradical agrees with that of $\mathsf{G}$; see \cite[Section~3.18]{Varadarajan84}.
The next proposition gives the required decomposition of $\mathsf{R}$, assuming only that $\mathsf{N}$ is simply connected.
The preceding lemma guarantees this assumption when $Z(\mathsf{G})$ is discrete.

\begin{prop}\label{prop:solvable-radical-decomposition}
Let $\mathsf{G}$ be a connected Lie group with connected solradical $\mathsf{R}$ and simply connected nilradical $\mathsf{N}$.
There exists a closed, connected, simply connected solvable subgroup $\mathsf{S}_{\mathsf{R}}\subseteq\mathsf{R}$ containing $\mathsf{N}$, normal in $\mathsf{G}$, such that
$$
    \mathsf{R}=\mathsf{K}_{\mathsf{R}}\ltimes\mathsf{S}_{\mathsf{R}}
$$
for every maximal compact subgroup $\mathsf{K}_{\mathsf{R}}\subseteq\mathsf{R}$.
\end{prop}

\begin{proof}
Since $[\mathfrak{r},\mathfrak{r}]\subseteq\mathfrak{n}$ by \cite[Corollary~1.41]{Knapp02}, the quotient $\mathsf{R}/\mathsf{N}$ is a connected abelian Lie group.
Choose an isomorphism $\mathsf{R}/\mathsf{N}\cong\mathbb{T}^m\times\mathbb{R}^n$, let $\pi:\mathsf{R}\to\mathsf{R}/\mathsf{N}$ be the projection, and define
$$
    \mathsf{S}_{\mathsf{R}}:=\pi^{-1}(\{e\}\times\mathbb{R}^n).
$$
This is a closed solvable subgroup containing $\mathsf{N}$, and $\mathsf{S}_{\mathsf{R}}/\mathsf{N}\cong\mathbb{R}^n$.
Since $\mathsf{N}$ is connected and simply connected by assumption, the homotopy exact sequence of the bundle $\mathsf{N}\to\mathsf{S}_{\mathsf{R}}\to\mathsf{S}_{\mathsf{R}}/\mathsf{N}$ shows that $\mathsf{S}_{\mathsf{R}}$ is connected and simply connected.

Let $\mathfrak{s}_{\mathsf{R}}$ be its Lie algebra.
By \cite[Corollary~1.41]{Knapp02},
$$
    [\mathfrak{g},\mathfrak{s}_{\mathsf{R}}]
    \subseteq[\mathfrak{g},\mathfrak{r}]
    \subseteq\mathfrak{n}
    \subseteq\mathfrak{s}_{\mathsf{R}}.
$$
Thus $\mathfrak{s}_{\mathsf{R}}$ is an ideal in $\mathfrak{g}$.
Since both groups are connected, $\mathsf{S}_{\mathsf{R}}$ is normal in $\mathsf{G}$.

Now let $\mathsf{K}_{\mathsf{R}}$ be any maximal compact subgroup of $\mathsf{R}$, and write $\mathsf{K}_{\mathsf{R}}\cong\mathbb{T}^k$.
The Cartan--Iwasawa--Malcev theorem gives $\pi_1(\mathsf{R})\cong\mathbb{Z}^k$.
On the other hand, the homotopy exact sequence of $\mathsf{N}\to\mathsf{R}\to\mathsf{R}/\mathsf{N}$ gives
$$
    \pi_1(\mathsf{R})\cong\pi_1(\mathsf{R}/\mathsf{N})\cong\mathbb{Z}^m,
$$
because $\mathsf{N}$ is connected and simply connected.
Hence $k=m$.

By Lemma~\ref{lem:simply-connected-solvable-free-action}, $\mathsf{K}_{\mathsf{R}}\cap\mathsf{N}=\{e\}$, so $\pi|_{\mathsf{K}_{\mathsf{R}}}$ is injective.
Its image is compact, and its projection to $\mathbb{R}^n$ is trivial by the same lemma applied to the additive group $\mathbb{R}^n$.
Thus $\pi(\mathsf{K}_{\mathsf{R}})\subseteq\mathbb{T}^m\times\{0\}$.
Both groups are connected and have dimension $m$, so
$$
    \pi(\mathsf{K}_{\mathsf{R}})=\mathbb{T}^m\times\{0\}.
$$
Together with the definition of $\mathsf{S}_{\mathsf{R}}$, this gives $\mathsf{R}/\mathsf{N}=\pi(\mathsf{K}_{\mathsf{R}})\pi(\mathsf{S}_{\mathsf{R}})$.
Since $\mathsf{N}\subseteq\mathsf{S}_{\mathsf{R}}$, it follows that $\mathsf{R}=\mathsf{K}_{\mathsf{R}}\mathsf{S}_{\mathsf{R}}$.

Finally, $\mathsf{K}_{\mathsf{R}}\cap\mathsf{S}_{\mathsf{R}}=\{e\}$ by Lemma~\ref{lem:simply-connected-solvable-free-action}.
Normality of $\mathsf{S}_{\mathsf{R}}$ now gives the asserted semidirect-product decomposition.
\end{proof}

\subsection{Proof of the generalized Iwasawa decomposition}\label{subsec:proof-generalized-iwasawa}

\textcolor{white}{aaa}

Recall that every connected Lie group admits a Levi decomposition $\mathsf{G}=\mathsf{L}\mathsf{R}$, where $\mathsf{R}$ is its connected solradical and $\mathsf{L}$ is a connected semisimple subgroup; see \cite[Theorem~3.18.13]{Varadarajan84}.
In general, $\mathsf{L}$ need not be closed, and $\mathsf{L}\cap\mathsf{R}$ is discrete but need not be trivial.

We divide the proof of Theorem~\ref{thm:generalized-iwasawa} into five steps, combining the decomposition $\mathsf{R}=\mathsf{K}_{\mathsf{R}}\ltimes\mathsf{S}_{\mathsf{R}}$ from Proposition~\ref{prop:solvable-radical-decomposition} with the Iwasawa decomposition of a suitable Levi factor $\mathsf{L}$.
First, we choose $\mathsf{L}$ to commute with $\mathsf{K}_{\mathsf{R}}$.
Second, we show that $\mathsf{L}$ is closed and has finite center, so the Iwasawa decomposition gives $\mathsf{L}=\mathsf{K}_{\mathsf{L}}\mathsf{S}_{\mathsf{L}}$, where $\mathsf{S}_{\mathsf{L}}:=\mathsf{A}_{\mathsf{L}}\mathsf{N}_{\mathsf{L}}$.
Third, we prove that $\mathsf{S}:=\mathsf{S}_{\mathsf{L}}\mathsf{S}_{\mathsf{R}}$ is closed, connected, simply connected and solvable.
Fourth, we set $\mathsf{K}:=\mathsf{K}_{\mathsf{L}}\mathsf{K}_{\mathsf{R}}$ and show that multiplication $\mathsf{K}\times\mathsf{S}\to\mathsf{G}$ is a diffeomorphism.
Finally, we prove that $\mathsf{K}$ is maximal compact.

\begin{proof}[Proof of Theorem~\ref{thm:generalized-iwasawa}]
It suffices to construct a maximal compact subgroup $\mathsf{K}\subseteq\mathsf{G}$ and a closed, connected, simply connected solvable subgroup $\mathsf{S}\subseteq\mathsf{G}$ containing $\mathsf{N}$ such that multiplication $\mathsf{K}\times\mathsf{S}\to\mathsf{G}$ is a diffeomorphism.
Indeed, suppose such subgroups have been constructed, and let $\mathsf{K}_0\subseteq\mathsf{G}$ be any prescribed maximal compact subgroup.
By conjugacy, $\mathsf{K}_0=g\mathsf{K}g^{-1}$ for some $g\in\mathsf{G}$.
Conjugating the decomposition gives a diffeomorphism $\mathsf{K}_0\times(g\mathsf{S}g^{-1})\to\mathsf{G}$ by multiplication.
The subgroup $g\mathsf{S}g^{-1}$ is closed, connected, simply connected and solvable, and contains $\mathsf{N}$ because $\mathsf{N}$ is normal in $\mathsf{G}$.

Fix a maximal compact subgroup $\mathsf{K}_{\mathsf{R}}\subseteq\mathsf{R}$, and let $\mathsf{R}=\mathsf{K}_{\mathsf{R}}\ltimes\mathsf{S}_{\mathsf{R}}$ be the decomposition from Proposition~\ref{prop:solvable-radical-decomposition}.

\textcolor{white}{aaa}

\noindent\emph{(i) A Levi factor commuting with $\mathsf{K}_{\mathsf{R}}$.}
Let $\pi:\mathsf{G}\to\mathsf{G}/\mathsf{R}$ be the projection and $\pi_*:\mathfrak{g}\to\mathfrak{g}/\mathfrak{r}$ its differential.
Consider the Lie algebra of the centralizer of $\mathsf{K}_{\mathsf{R}}$,
$$
    \mathfrak{c}:=\mathrm{Lie}(Z_{\mathsf{G}}(\mathsf{K}_{\mathsf{R}}))
    =\{U\in\mathfrak{g}:\mathrm{Ad}_{\mathsf{G}}(k)U=U\text{ for every }k\in\mathsf{K}_{\mathsf{R}}\}.
$$
We want to choose a Levi subalgebra inside $\mathfrak{c}$, so that the corresponding Levi factor commutes with $\mathsf{K}_{\mathsf{R}}$.
To ensure that it is also a Levi subalgebra of $\mathfrak{g}$, we need to show that $\mathfrak{c}\cap\mathfrak{r}$ is the solradical of $\mathfrak{c}$.
We begin by establishing $\pi_*(\mathfrak{c})=\mathfrak{g}/\mathfrak{r}$.

Let $\mu$ be the normalized Haar measure on $\mathsf{K}_{\mathsf{R}}$, so $\mu(\mathsf{K}_{\mathsf{R}})=1$.
For $U\in\mathfrak{g}$, define
$$
    \overline{U}:=\int_{\mathsf{K}_{\mathsf{R}}}\mathrm{Ad}_{\mathsf{G}}(k)U\,d\mu(k).
$$
For $k\in\mathsf{K}_{\mathsf{R}}$, left invariance of $\mu$ gives
$$
    \mathrm{Ad}_{\mathsf{G}}(k)\overline{U}
    =\int_{\mathsf{K}_{\mathsf{R}}}\mathrm{Ad}_{\mathsf{G}}(k\tilde{k})U\,d\mu(\tilde{k})
    =\int_{\mathsf{K}_{\mathsf{R}}}\mathrm{Ad}_{\mathsf{G}}(\tilde{k})U\,d\mu(\tilde{k})
    =\overline{U}.
$$
Thus $\overline{U}\in\mathfrak{c}$.
Since $k\in\mathsf{K}_{\mathsf{R}}\subseteq\ker\pi$, we have $\pi(kgk^{-1})=\pi(g)$ for every $g\in\mathsf{G}$.
Differentiating at the identity gives $\pi_*(\mathrm{Ad}_{\mathsf{G}}(k)U)=\pi_*(U)$.
Therefore
$$
    \pi_*(\overline{U})
    =\int_{\mathsf{K}_{\mathsf{R}}}\pi_*(\mathrm{Ad}_{\mathsf{G}}(k)U)\,d\mu(k)
    =\int_{\mathsf{K}_{\mathsf{R}}}\pi_*(U)\,d\mu(k)
    =\pi_*(U).
$$
This proves $\pi_*(\mathfrak{c})=\mathfrak{g}/\mathfrak{r}$, or equivalently $\mathfrak{g}=\mathfrak{c}+\mathfrak{r}$.

The intersection $\mathfrak{c}\cap\mathfrak{r}$ is a solvable ideal in $\mathfrak{c}$.
Since $\pi_*(\mathfrak{c})=\mathfrak{g}/\mathfrak{r}$, the quotient $\mathfrak{c}/(\mathfrak{c}\cap\mathfrak{r})\cong\mathfrak{g}/\mathfrak{r}$ is semisimple.
Thus $\mathfrak{c}\cap\mathfrak{r}$ is the solradical of $\mathfrak{c}$.
Choose a Levi decomposition $\mathfrak{c}=\mathfrak{l}\ltimes(\mathfrak{c}\cap\mathfrak{r})$.
Then $\mathfrak{l}\cap\mathfrak{r}=0$ and $\mathfrak{g}=\mathfrak{c}+\mathfrak{r}=\mathfrak{l}+\mathfrak{r}$, so $\mathfrak{g}=\mathfrak{l}\ltimes\mathfrak{r}$ is a Levi decomposition.

Let $\mathsf{L}$ be the connected immersed subgroup with Lie algebra $\mathfrak{l}$.
It is a Levi factor, and $\mathfrak{l}\subseteq\mathfrak{c}$ implies $\mathsf{L}\subseteq Z_{\mathsf{G}}(\mathsf{K}_{\mathsf{R}})$.
Thus $\mathsf{L}$ commutes with $\mathsf{K}_{\mathsf{R}}$.

\textcolor{white}{aaa}

\noindent\emph{(ii) The Levi factor is closed and has finite center.}
Consider the adjoint homomorphism $\mathrm{Ad}_{\mathsf{G}}:\mathsf{G}\to\mathsf{GL}(\mathfrak{g})$ and the connected immersed subgroup $\mathrm{Ad}_{\mathsf{G}}(\mathsf{L})\subseteq\mathsf{GL}(\mathfrak{g})$, whose Lie algebra is $\mathrm{ad}_{\mathfrak{g}}(\mathfrak{l})\subseteq\mathrm{End}(\mathfrak{g})$.
Since $\mathfrak{z}(\mathfrak{g})=0$, the restriction $\mathrm{ad}_{\mathfrak{g}}:\mathfrak{l}\to\mathrm{ad}_{\mathfrak{g}}(\mathfrak{l})$ is a Lie algebra isomorphism.
Thus $\mathrm{Ad}_{\mathsf{G}}(\mathsf{L})$ is semisimple, and \cite[Proposition~7.9]{Knapp02} implies that it is closed in $\mathsf{GL}(\mathfrak{g})$ and has finite center.

Consequently,
$$
    \mathsf{L}Z(\mathsf{G})
    =\mathrm{Ad}_{\mathsf{G}}^{-1}\bigl(\mathrm{Ad}_{\mathsf{G}}(\mathsf{L})\bigr)
$$
is a closed subgroup of $\mathsf{G}$.
Since $Z(\mathsf{G})$ is finite, $\mathsf{L}$ is the identity component of the closed subgroup $\mathsf{L}Z(\mathsf{G})$, and hence $\mathsf{L}$ is closed in $\mathsf{G}$.

To see that $Z(\mathsf{L})$ is finite, consider the restriction $\mathrm{Ad}_{\mathsf{G}}:Z(\mathsf{L})\to Z(\mathrm{Ad}_{\mathsf{G}}(\mathsf{L}))$.
It has finite kernel $Z(\mathsf{L})\cap Z(\mathsf{G})\subseteq Z(\mathsf{G})$ and finite image contained in $Z(\mathrm{Ad}_{\mathsf{G}}(\mathsf{L}))$.
In particular,
$$
    |Z(\mathsf{L})|
    \leq |Z(\mathsf{G})|\cdot|Z(\mathrm{Ad}_{\mathsf{G}}(\mathsf{L}))|
    <\infty.
$$

We may therefore apply the Iwasawa decomposition \cite[Theorem~6.46]{Knapp02} to write $\mathsf{L}=\mathsf{K}_{\mathsf{L}}\mathsf{S}_{\mathsf{L}}$, where $\mathsf{K}_{\mathsf{L}}$ is maximal compact and $\mathsf{S}_{\mathsf{L}}:=\mathsf{A}_{\mathsf{L}}\mathsf{N}_{\mathsf{L}}$ is closed, connected, simply connected and solvable.
Since $\mathsf{L}$ is closed, $\mathsf{S}_{\mathsf{L}}$ is also closed in $\mathsf{G}$.

We will also need $\mathsf{L}\cap\mathsf{R}\subseteq\mathsf{K}_{\mathsf{L}}$.
Indeed, $\mathsf{L}\cap\mathsf{R}$ is a discrete normal subgroup of the connected group $\mathsf{L}$, so it lies in $Z(\mathsf{L})$ and is therefore finite.
Hence $(\mathsf{L}\cap\mathsf{R})\mathsf{K}_{\mathsf{L}}$ is a compact subgroup containing $\mathsf{K}_{\mathsf{L}}$.
Maximality gives the claimed inclusion.

\textcolor{white}{aaa}

\noindent\emph{(iii) Construction of the solvable subgroup.}
Define $\mathsf{S}:=\mathsf{S}_{\mathsf{L}}\mathsf{S}_{\mathsf{R}}$.
This is a subgroup because $\mathsf{S}_{\mathsf{R}}$ is normal in $\mathsf{G}$ by Proposition~\ref{prop:solvable-radical-decomposition}.
To establish closedness, consider the external semidirect product $\mathsf{L}\ltimes\mathsf{R}$ defined by the conjugation action of $\mathsf{L}$ on $\mathsf{R}$.
The multiplication map $\Phi:\mathsf{L}\ltimes\mathsf{R}\to\mathsf{G}$, $\Phi(l,r):=lr$, is a surjective Lie-group homomorphism with kernel
$$
    \ker\Phi=\{(d,d^{-1}):d\in\mathsf{L}\cap\mathsf{R}\}.
$$
By Step~(ii), this kernel is finite.
Thus $\Phi$ is a finite-sheeted covering, and in particular a closed map.
The subset $\mathsf{S}_{\mathsf{L}}\times\mathsf{S}_{\mathsf{R}}$ is closed and connected in $\mathsf{L}\ltimes\mathsf{R}$, so its image $\mathsf{S}$ is closed and connected in $\mathsf{G}$.

Furthermore,
$$
    \mathsf{S}_{\mathsf{L}}\cap\mathsf{S}_{\mathsf{R}}
    \subseteq\mathsf{S}_{\mathsf{L}}\cap(\mathsf{L}\cap\mathsf{R})
    \subseteq\mathsf{S}_{\mathsf{L}}\cap\mathsf{K}_{\mathsf{L}}
    =\{e\}.
$$
Together with normality of $\mathsf{S}_{\mathsf{R}}$, this gives $\mathsf{S}=\mathsf{S}_{\mathsf{L}}\ltimes\mathsf{S}_{\mathsf{R}}$.
Both factors are solvable and simply connected, so $\mathsf{S}$ has these properties as well.
Finally, $\mathsf{N}\subseteq\mathsf{S}_{\mathsf{R}}\subseteq\mathsf{S}$.

\textcolor{white}{aaa}

\noindent\emph{(iv) The product decomposition.}
Set $\mathsf{K}:=\mathsf{K}_{\mathsf{L}}\mathsf{K}_{\mathsf{R}}$.
The factors are compact and commute by Step~(i), so $\mathsf{K}$ is a compact subgroup.
Since $\mathsf{S}_{\mathsf{L}}$ also commutes with $\mathsf{K}_{\mathsf{R}}$,
$$
    \mathsf{G}
    =\mathsf{L}\mathsf{R}
    =\mathsf{K}_{\mathsf{L}}\mathsf{S}_{\mathsf{L}}\mathsf{K}_{\mathsf{R}}\mathsf{S}_{\mathsf{R}}
    =\mathsf{K}_{\mathsf{L}}\mathsf{K}_{\mathsf{R}}\mathsf{S}_{\mathsf{L}}\mathsf{S}_{\mathsf{R}}
    =\mathsf{K}\mathsf{S}.
$$
The intersection $\mathsf{K}\cap\mathsf{S}$ is compact because $\mathsf{S}$ is closed, and hence trivial by Lemma~\ref{lem:simply-connected-solvable-free-action}.

Since $\mathsf{G}=\mathsf{K}\mathsf{S}$ and $\mathsf{K}\cap\mathsf{S}=\{e\}$, their Lie algebras are complementary in $\mathfrak{g}$ and multiplication $\mathsf{K}\times\mathsf{S}\to\mathsf{G}$ is bijective.
By \cite[Lemma~6.44]{Knapp02}, this map is also a local diffeomorphism, hence a diffeomorphism.

\textcolor{white}{aaa}

\noindent\emph{(v) Maximality of the compact subgroup.}
Let $\mathsf{C}$ be a compact subgroup containing $\mathsf{K}$.
For $c\in\mathsf{C}$, write $c=ks$ using Step~(iv).
Then $s=k^{-1}c\in\mathsf{C}\cap\mathsf{S}$.
This intersection is trivial by Lemma~\ref{lem:simply-connected-solvable-free-action}.
Hence $c=k\in\mathsf{K}$, proving that $\mathsf{K}$ is maximal compact.
\end{proof}

\begin{remark}\label{rem:cornulier-iwasawa}
We now give an alternative proof of Theorem~\ref{thm:generalized-iwasawa} using the construction in the proof of \cite[Lemma~6.7]{Cornulier08}.
Cornulier constructs a closed, connected solvable subgroup $\mathsf{S}$ containing $\mathsf{N}$, with $\pi_1(\mathsf{S})\cong\pi_1(\mathsf{N})$, and a cocompact subgroup $\mathsf{S}\times\Gamma$, where $\Gamma$ is free abelian and central in $\mathsf{G}$.

Since we assume that $Z(\mathsf{G})$ is finite, $\Gamma=\{e\}$, so $\mathsf{S}$ is cocompact in $\mathsf{G}$.
Moreover, Lemma~\ref{lem:simply-connected-nilradical} implies that $\mathsf{N}$ is simply connected, and hence so is $\mathsf{S}$.
By Lemma~\ref{lem:simply-connected-solvable-free-action}, $\mathsf{S}$ has no nontrivial compact subgroups and is diffeomorphic to $\mathbb{R}^{\dim\mathsf{S}}$, hence contractible.

The homotopy exact sequence of the fibre bundle $\mathsf{S}\to\mathsf{G}\to\mathsf{G}/\mathsf{S}$, together with Whitehead's theorem, therefore shows that $\mathsf{G}/\mathsf{S}$ is homotopy equivalent to $\mathsf{G}$.
For any maximal compact subgroup $\mathsf{K}\subseteq\mathsf{G}$, the Cartan--Iwasawa--Malcev theorem also gives $\mathsf{G}\simeq\mathsf{K}$.
Thus $\mathsf{G}/\mathsf{S}$ and $\mathsf{K}$ are homotopy equivalent compact manifolds, and comparing their top homology with $\mathbb{Z}_2$-coefficients gives $\dim(\mathsf{G}/\mathsf{S})=\dim\mathsf{K}$.

For the $\mathsf{K}$-action on $\mathsf{G}/\mathsf{S}$, the isotropy subgroup at $g\mathsf{S}$ is the compact subgroup $\mathsf{K}\cap g\mathsf{S}g^{-1}$.
Conjugation by $g^{-1}$ identifies it with a compact subgroup of $\mathsf{S}$, which is trivial by Lemma~\ref{lem:simply-connected-solvable-free-action}.
Hence the action is free.

Since $\dim\mathsf{K}=\dim(\mathsf{G}/\mathsf{S})$, its orbits are open, and compactness of $\mathsf{K}$ makes them closed.
Connectedness of $\mathsf{G}/\mathsf{S}$ therefore implies that the action is transitive.
Consequently, $\mathsf{G}/\mathsf{S}$ is diffeomorphic to $\mathsf{K}$, and $\mathsf{G}=\mathsf{K}\mathsf{S}$ with $\mathsf{K}\cap\mathsf{S}=\{e\}$.
Their Lie algebras are complementary, so \cite[Lemma~6.44]{Knapp02} implies that the bijective multiplication map $\mathsf{K}\times\mathsf{S}\to\mathsf{G}$ is a diffeomorphism.
\end{remark}

\end{appendix}

\bibliographystyle{alpha}
\bibliography{mybib}

\section*{Notation index}
\phantomsection\makeatletter\def\@currentlabel{Notation index}\makeatother\label{sec:notation-index}

This index includes notation used in more than one section of the paper.
Notation confined to a single section or appendix is omitted.

We recall the standing convention (see Subsection~\ref{subsec:nilradical-action}): we use a hat to distinguish the tensors and distributions associated with the Riemannian submersion $M\rightarrow M/\mathsf{N}$.
These hatted tensors and distributions are not included in this index.

\subsection*{The action, its distributions and holonomy}

\begin{description}
\item[$\mathrm{g}$, $\langle\cdot,\cdot\rangle$]
The Riemannian metric on $M$.

\item[$\mathsf{G}$, $\mathfrak{g}$]
The connected symmetry group and its Lie algebra.

\item[$\mathsf{N}$, $\mathfrak{n}$]
The connected nilradical of $\mathsf{G}$ and its Lie algebra.

\item[$m_g$, $r_x$]
The action map $m_g:M\to M$, $m_g(x)=g\cdot x$, and the orbit map $r_x:\mathsf{G}\to M$, $r_x(g)=g\cdot x$.

\item[$\mathcal{V}$, $\mathcal{H}$]
The vertical and horizontal distributions, $\mathcal{V}_x=T_x(\mathsf{G}\cdot x)$ and $\mathcal{H}_x=\mathcal{V}_x^\perp$.

\item[$\mathfrak{a}_x$, $\mathcal{A}_x$]
The $h_x$-orthogonal complement of $\mathfrak{n}$ in $\mathfrak{g}$ and the corresponding distribution $\mathcal{A}_x=\mathcal{V}_x\cap\widehat{\mathcal{H}}_x$.

\item[$U,V,W$;\\ $X,Y,Z$; $\xi,\tilde{\xi}$]
Vertical fields; horizontal fields; and arbitrary vector fields, respectively.

\par\vspace{\baselineskip}

\item[$\Gamma(E)$]
The space of smooth sections of a vector bundle $E$, over its specified base.

\item[$\mathsf{Hol}_x$]
The holonomy group of $M\to M/\mathsf{G}$, regarded as an immersed Lie subgroup of $\mathsf{G}$; see Subsection~\ref{subsec:holonomy}.

\item[$\widehat{\Sigma}$]
A connected $\mathsf{N}$-horizontal leaf.

\item[$\mathfrak{a}_0$, $\mathsf{A}_0$]
For a fixed $\mathsf{N}$-horizontal leaf $\widehat{\Sigma}$, the constant subalgebra $\mathfrak{a}_0=\mathfrak{a}_x$, $x\in\widehat{\Sigma}$, and the subgroup $\mathsf{A}_0=\exp(\mathfrak{a}_0)$.

\end{description}

\subsection*{Killing coordinates and inner products}

\begin{description}
\item[$U^*$]
The fundamental Killing field induced by $U\in\mathfrak{g}$.

\item[$\theta$, $\theta_x$]
The trivialization $\theta:M\times\mathfrak{g}\to\mathcal{V}$ defined by $\theta(x,U)=U_x^*$ and its fiberwise isomorphisms $\theta_x:\mathfrak{g}\rightarrow\mathcal{V}_x$.

\item[${}^{\mathfrak{g}}U$, ${}^{\mathfrak{g}}E$]
The Killing coordinates of a vertical field and a vertical endomorphism; see Subsection~\ref{subsec:trivialization}.

\item[$\mathrm{ad}(A)$]
For an $\mathcal{A}$-valued field $A$, the vertical endomorphism with Killing coordinates
${}^{\mathfrak{g}}\mathrm{ad}(A)_x=\mathrm{ad}_{{}^{\raisebox{0.15em}{$\scriptstyle\mathfrak{g}$}}\! A_x}$.

\item[$\mathrm{Sym}^2_+(\mathfrak{g}^*)$]
The space of inner products on $\mathfrak{g}$, with the Riemannian metric described in Appendix~\ref{app:inner-products}.

\item[$h_x$, $\hat{h}_x$]
The inner products on $\mathfrak{g}$ and $\mathfrak{n}$ induced by the orbit metrics; see Subsection~\ref{subsec:trivialization}.

\item[$\rho(E)h$]
The infinitesimal action on symmetric bilinear forms,
$\rho(E)h=-h(E\cdot,\cdot)-h(\cdot,E\cdot)$;
see Appendix~\ref{app:inner-products}.

\item[$E^{t_h}$, $E^{t_x}$]
The adjoints with respect to the inner product $h$ and $h_x$, respectively.

\item[$\mathrm{S}_x(E),\mathrm{sk}_x(E)$]
The $h_x$-symmetric and $h_x$-skew-symmetric parts of $E$.
\end{description}

\subsection*{Connections, curvature and differential operators}

\begin{description}
\item[$\nabla$]
The Levi-Civita connection of $M$.

\item[$\nabla^{\mathrm{ad}}$]
The connection induced in the adjoint bundle; see Definition~\ref{def:adjoint-connection}.

\item[$\mathrm{A}$, $L$]
The O'Neill tensors for the Riemannian submersion $M\to M/\mathsf{G}$; see Subsection~\ref{subsec:oneill-tensors}.

\item[$H$]
The mean curvature field of the $\mathsf{G}$-orbits; see Subsection~\ref{subsec:oneill-curvature}.

\item[$\mathrm{A}^*\mathrm{A}$, $\tau$]
The self-adjoint vertical endomorphisms involved in O'Neill's vertical Ricci formula; see Subsection~\ref{subsec:oneill-curvature}.

\item[$\beta^+$, ${}^{\mathfrak{n}}\beta^+$]
The positive definite $\mathsf{N}$-vertical endomorphism arising from GIT and its Killing coordinates; see Subsection~\ref{subsec:git}.

\item[$\mathrm{ric}$, $\mathrm{Ric}$]
The Ricci tensor and Ricci endomorphism of $M$, related by
$\mathrm{ric}(\xi,\tilde{\xi})=\langle\mathrm{Ric}\,\xi,\tilde{\xi}\rangle$.
A superscript with another Riemannian manifold refers to its Ricci tensor or endomorphism.
A superscript with a Lie group refers to the corresponding orbits.

\item[$\mathrm{Ric}^{\mathcal{V}}$]
The vertical compression $U\mapsto\mathcal{V}(\mathrm{Ric}\,U)$, for $U\in\mathcal{V}$.

\item[$\mathrm{div}_B$, $\Delta_B$]
The divergence on $M/\mathsf{G}$, pulled back to basic vector fields on $M$, and the operator $\Delta_Bf=\mathrm{div}_B(\nabla f)$ on invariant functions; see Subsection~\ref{subsec:divergence-maximum-principle}.
\end{description}

\end{document}